\documentclass[11pt]{amsart}
\usepackage[T1]{fontenc}
\usepackage{lmodern}
\usepackage{amssymb,mathtools}
\usepackage[letterpaper,hmargin=1.15in,top=1.1in,bottom=1.1in]{geometry}
\usepackage[hidelinks]{hyperref}

\newtheorem{theorem}{Theorem}[section]
\newtheorem{lemma}[theorem]{Lemma}
\newtheorem{proposition}[theorem]{Proposition}
\newtheorem{corollary}[theorem]{Corollary}
\theoremstyle{remark}
\newtheorem{remark}[theorem]{Remark}
\numberwithin{equation}{section}
\newcommand{\N}{\mathbb{N}}
\newcommand{\Z}{\mathbb{Z}}
\newcommand{\Q}{\mathbb{Q}}
\newcommand{\kk}{\Bbbk}
\newcommand{\T}{\mathsf{T}}

\newcommand{\Cas}{\operatorname{Cas}}
\newcommand{\crd}{\operatorname{crd}}
\newcommand{\detR}{\det\nolimits_R}
\newcommand{\detE}{\det\nolimits_E}
\newcommand{\Mod}{\mathfrak{M}}
\newcommand{\detM}{\det\nolimits_{\Mod}}
\newcommand{\detx}[1]{\det\nolimits_{#1}}
\newcommand{\eps}{\varepsilon}
\newcommand{\cW}{\mathcal{W}}
\newcommand{\cH}{\mathcal{H}}
\newcommand{\cE}{\mathcal{E}}
\newcommand{\cO}{\mathcal{O}}
\newcommand{\cN}{\mathcal{N}}
\newcommand{\cL}{\mathcal{L}}
\newcommand{\fn}{\mathfrak{n}}
\newcommand{\Nar}{\mathrm{Nar}}
\newcommand{\symp}{\mathrm{sp}}
\newcommand{\sorth}{\mathrm{so}}
\newcommand{\Sp}{\mathrm{Sp}}
\newcommand{\SO}{\mathrm{SO}}
\newcommand{\Rt}{\widetilde{R}}
\newcommand{\Ft}{\widetilde{F}}

\begin{document}

\title[Hankel determinants of Motzkin paths and the Hirota--Miwa equation]{Hankel determinants of Motzkin path generating functions and the Hirota--Miwa equation}

\author{Haiming Sun}
\address{Peking University, Beijing, China}
\email{2400012896@stu.pku.edu.cn}

\date{November 2026}

\begin{abstract}
For $k\ge0$ let $a_{n,k}(x,c)$ count Motzkin paths from $(0,0)$ to $(n,k)$ that never go below height $0$, where level steps have weight $c+x$ at height $0$ and $c$ above it, and let $D_{k,n}=\det\bigl(a_{i+j,k}(x,c)\bigr)_{0\le i,j\le n-1}$. T.~Yu proved Cigler's conjecture $D_{k,n}D_{k,n-2k-2}-x^2D_{k,n-1}D_{k,n-2k-1}=(1-x^2)D_{k,n-k-1}^2$ $(n\ge2k+2)$ by writing $D_{k,n}$ as a $k\times k$ determinant of Chebyshev polynomials modulo $U_k-xU_{k-1}$. Using a linear map found by Yu, we show that the recurrence is a reduction of the Hirota--Miwa equation for Casorati determinants. This gives a second proof, and the recurrence holds for $n\ge2k+2+\max(0,N-\lfloor k/2\rfloor)$ for the Hankel determinants of $n\mapsto L(X^n(U_k-xU_{k-1}))$ whenever the moment functional $L$ has the recurrence coefficients of the Chebyshev polynomials above height $N$. For any monic polynomial in place of $U_k-xU_{k-1}$ the Hirota--Miwa equation gives an identity with a free parameter between symplectic characters, and Yu's identity for rectangular characters follows from it by fixing the parameter and restricting to a line. On the curve $xc+x^2+1=0$ we express the Hankel determinants of the sequences shifted by $2$ and $-1$ through the $D_{k,n}$. There the recurrence becomes the Desnanot--Jacobi identity, which gives a third proof, and a recurrence for the Hankel determinants of odd convolution powers of Narayana polynomials, shifted backward by $k$. We also prove identities that Cigler stated as observations, an analogue of odd orthogonal type, and a factorization for paths with a weighted return step.
\end{abstract}

\maketitle
\section{Introduction}

\subsection{Cigler's recurrence}\label{ss:cigler}
Let $x,c$ be indeterminates and let $k\ge0$ be an integer. Put $s_0=c+x$ and $s_j=c$ for $j\ge1$, and define $a_{n,j}=a_{n,j}(x,c)\in\Z[x,c]$ for $n\ge0$ and $j\in\Z$ by
\[
a_{0,j}=\delta_{j,0},\qquad a_{n,j}=0\quad(j<0),\qquad a_{n,j}=a_{n-1,j-1}+s_j\,a_{n-1,j}+a_{n-1,j+1}\quad(n\ge1,\ j\ge0).
\]
Thus $a_{n,j}$ is the weighted number of Motzkin paths from $(0,0)$ to $(n,j)$ that never go below height $0$, where up and down steps have weight $1$, level steps at height $0$ have weight $c+x$, and level steps at positive height have weight $c$. For $n\ge0$ let
\[
D_{k,n}=D_{k,n}(x,c)=\det\bigl(a_{i+j,k}\bigr)_{0\le i,j\le n-1},\qquad D_{k,0}=1.
\]
In \cite{CMO} and \cite[eq.~(18)]{C23}, Cigler conjectured a quadratic recurrence for the sequence $(D_{k,n})_{n\ge0}$. (In \cite{C23} the determinant $D_{k,n}(x,c)$ is denoted $D_{0,k,n}(b,c)$, with $b$ in place of $x$.) The conjecture was recently proved by Yu.

\begin{theorem}[Yu {\cite[Theorem~B]{Yu}}; conjectured by Cigler]\label{thm:B}
For all $k\ge0$ and all $n\ge2k+2$,
\begin{equation}\label{eq:main}
D_{k,n}\,D_{k,n-2k-2}-x^2\,D_{k,n-1}\,D_{k,n-2k-1}=(1-x^2)\bigl(D_{k,n-k-1}\bigr)^2\quad\text{in }\Z[x,c].
\end{equation}
\end{theorem}

The determinants $D_{k,n}(x,c)$ do not depend on $c$. This is an instance of the invariance of Hankel determinants under the binomial transform \cite{L}, \cite[Theorem~6.1]{SS} (see Corollary~\ref{cor:cindep}), and it was conjectured by Cigler \cite[Conjecture~2.2, eq.~(2.4)]{C22}. We also write $D_{k,n}(x)$ for $D_{k,n}(x,c)$. Yu's proof of Theorem~\ref{thm:B} rests on the following result.

\begin{theorem}[Yu {\cite[Theorem~C]{Yu}}]\label{thm:C}
Let $k\ge1$, let $U_m\in\Z[X]$ $(m\in\Z)$ be defined by $U_{-1}=0$, $U_0=1$ and $U_{m+1}=XU_m-U_{m-1}$ for all $m\in\Z$, and let $P_k=U_k-xU_{k-1}$. For $m\in\Z$ let $w(m)$ be the determinant of the $k\times k$ matrix whose rows are the coordinate vectors, with respect to the basis of $\Q(x)[X]/(P_k)$ formed by the residue classes of $U_0,\dots,U_{k-1}$, of the residue classes of $U_{m-1},U_m,\dots,U_{m+k-2}$. Then $w(m)\in\Z[x]$, and
\begin{gather}
D_{k,n}=(-1)^{nk+k(k+1)/2}\,w(n-k)\qquad\text{for all } n\ge0,\label{eq:C1}\\
w(\ell+k+1)\,w(\ell-k-1)-x^2\,w(\ell+k)\,w(\ell-k)=(1-x^2)\,w(\ell)^2\qquad\text{for all }\ell\in\Z.\label{eq:C2}
\end{gather}
\end{theorem}

Theorem~\ref{thm:B} follows from \eqref{eq:C1} and \eqref{eq:C2}. The formula \eqref{eq:C1} is a special case of a formula of Christoffel type for Hankel determinants of linear combinations of moments \cite[Theorem~1]{K23}, combined with the relation $F_{m+1}-xF_m=\Lambda F_{m-k}$ of \eqref{eq:introLambda} below and the value $\det\Lambda=(-1)^{k(k+1)/2}$ (see the proof of Theorem~\ref{thm:C} in \S\ref{ss:proofs}). The identity \eqref{eq:C2} is proved in \cite{Yu} by means of a model of $\Q(x)[X]/(P_k)$ by antisymmetric Laurent polynomials, a three-term Pl\"ucker relation for a basis of the two-dimensional kernel of a window of $k+2$ consecutive residue classes, and a divisibility argument.

In this paper we show that \eqref{eq:C2} is a reduction of the Hirota--Miwa equation. This gives a second proof of Theorems~\ref{thm:B} and~\ref{thm:C}. A third proof of Theorem~\ref{thm:B} rests on the Desnanot--Jacobi identity. Yu \cite[Remark~4.4]{Yu} interprets $w(m)$ as a symplectic character of rectangular shape and asks whether \eqref{eq:C2} is a specialization of a known identity for symplectic characters. Remark~\ref{rem:symplectic} gives a partial answer. It shows that \eqref{eq:C2} is obtained, by fixing a parameter and restricting to a line, from an identity between symplectic characters that holds at every point and is a form of the three-term Pl\"ucker relation.

\subsection{Results}\label{ss:results}
Let $F_j$ $(j\in\Z)$ be the residue class of $U_{j-1}$ modulo $P_k$, so that $w(m)$ is the determinant of the $k$ consecutive terms $F_m,\dots,F_{m+k-1}$ (a \emph{Casorati determinant}) of the sequence $(F_j)_{j\in\Z}$. The recurrence comes from an invertible linear map $\Lambda$ of $\Q(x)[X]/(P_k)$, found by Yu, such that
\begin{equation}\label{eq:introLambda}
F_{j+1}-xF_j=\Lambda F_{j-k}\qquad\text{and}\qquad F_j-xF_{j+1}=\Lambda F_{j+k+1}\qquad\text{for all }j\in\Z.
\end{equation}
These relations are \cite[Lemma~3.6(ii)]{Yu}, proved there in the Laurent polynomial model. The map $\Lambda$ is the multiplication by $-P_{k-1}$ \cite[Remark~4.3]{Yu}, and in Lemma~\ref{lem:Lambda} we verify \eqref{eq:introLambda} directly in the quotient ring, over an arbitrary field.

\emph{The Hirota--Miwa equation $($Section~$\ref{sec:HM})$.} For any sequence of vectors $(F_j)_{j\in\Z}$ in a space of dimension $r$, the $r\times r$ determinants formed from $r$ consecutive terms of the sequences $F$, $(F_{j+1}-aF_j)_j$ and $(F_{j+2}-(a+b)F_{j+1}+abF_j)_j$ satisfy a three-term relation, the Hirota--Miwa (discrete KP) equation for Casorati determinants \cite{Hir,Miw,OHTI} (Lemma~\ref{lem:HM}). By \eqref{eq:introLambda}, the operations $F_j\mapsto F_{j+1}-xF_j$ and $F_j\mapsto F_{j+1}-x^{-1}F_j$ are, up to the fixed map $\Lambda$ and a scalar, translations of the index by $-k$ and by $k+1$. Under this identification the Hirota--Miwa equation becomes \eqref{eq:C2}, and the coefficient $1-x^2=x(x^{-1}-x)$ comes from the difference $x^{-1}-x$ of the two parameters. More generally, if each of three operations $F_j\mapsto F_{j+1}-a_iF_j$ is a translation of the index followed by a fixed invertible linear map, then the Casorati determinants satisfy a three-term bilinear recurrence of Gale--Robinson type (Proposition~\ref{thm:reduction}).

\emph{General moment functionals $($\S\ref{ss:general}$)$.} For sequences $s=(s_j)_{j\ge0}$ and $\lambda=(\lambda_j)_{j\ge1}$ in a commutative ring $A$ let $L=L^{(s,\lambda)}$ be the $A$-linear functional on $A[X]$ for which $L(X^n)$ is the weighted number of Motzkin paths from $(0,0)$ to $(n,0)$ that never go below height $0$, where level steps at height $j$ have weight $s_j$ and down steps from height $j$ have weight $\lambda_j$ (\S\ref{ss:matrices}). Its monic orthogonal polynomials satisfy $p_{m+1}=(X-s_m)p_m-\lambda_mp_{m-1}$, and for $s_j=0$, $\lambda_j=1$ they are the polynomials $U_m$. Theorem~\ref{thm:hankelgen} states the following. Let $k\ge1$, $N\ge0$ and $x\in A$, suppose that $s_n=0$ and $\lambda_n=1$ for all $n>N$, and let $\cH_n=\det(b_{i+j})_{0\le i,j\le n-1}$ with $b_n=L\bigl(X^n(U_k-xU_{k-1})\bigr)$. Then
\[
\cH_n\,\cH_{n-2k-2}-x^2\,\cH_{n-1}\,\cH_{n-2k-1}=(1-x^2)\,\cH_{n-k-1}^{\,2}\qquad\text{for all }n\ge2k+2+\max\bigl(0,N-\lfloor k/2\rfloor\bigr).
\]
Theorem~\ref{thm:B} is the case $N=0$, $s_0=x$ (with $c=0$, which suffices by Corollary~\ref{cor:cindep}). The bound can be lowered slightly if $s_N=0$, and in the examples we computed the resulting bounds are attained (Remark~\ref{rem:range}). The theorem follows from Corollary~\ref{thm:module}, by which \eqref{eq:C2} holds for the Casorati determinants, of any size, of every solution of $f_{j+1}=Xf_j-f_{j-1}$ in a finite-dimensional module over $\kk[X]/(U_k-xU_{k-1})$, for any field $\kk$ and $x\in\kk$.

\emph{A spectral parameter $($\S\ref{ss:spectral}$)$.} For every monic polynomial $P$ of degree $k$ in place of $P_k$, the Hirota--Miwa equation gives an identity with a free parameter $y$ (Proposition~\ref{thm:spectral}). In terms of the roots of $P$ it is an identity between symplectic characters of rectangular and nearly rectangular shape that holds at every point. The identity \eqref{eq:C2} is obtained from it by putting $y=x$ and restricting to the line on which all fundamental characters except the first vanish (Remark~\ref{rem:symplectic}). The identity with the parameter $y$ also follows from the Desnanot--Jacobi identity (Remark~\ref{rem:DJ}).

\emph{Shifted Hankel determinants $($Section~$\ref{sec:shifted})$.} For $m\in\Z$ let
\[
D_{m,k,n}(x,c)=\det\bigl(a_{i+j+m,k}(x,c)\bigr)_{0\le i,j\le n-1},\qquad a_{r,k}:=0\ \text{ for }r<0,
\]
as in \cite{C23}. For all $(x,c)$ the determinants with backward shifts are forward shifted Hankel determinants of the sequence $r\mapsto a_{r,k}(0,c)-x\,a_{r,k-1}(0,c)$. Namely, for $k\ge1$ and $m,n\ge0$ (Proposition~\ref{prop:backward})
\[
D_{-m,k,n+m+k+1}(x,c)=(-1)^{\binom{m+k+1}2}\det\bigl(a_{i+j+m,k}(0,c)-x\,a_{i+j+m,k-1}(0,c)\bigr)_{0\le i,j\le n-1}.
\]
For $x=0$ this is Cigler's duality \cite[Theorem~2]{C24}. On the curve $xc+x^2+1=0$, that is, for $c=-x-x^{-1}$, the determinants with the shifts $2$ and $-1$ can be expressed through $d_n:=D_{k,n}(x)$. For $k\ge1$ we show (Theorem~\ref{thm:curve})
\begin{align*}
(x^2-1)\,D_{2,k,n}\bigl(x,-x-x^{-1}\bigr)&=x^{k+2}\,d_{n+1}-x^{-2n}\,d_{n+2k+2}&&(n\ge0),\\
(x^2-1)\,D_{-1,k,n}\bigl(x,-x-x^{-1}\bigr)&=(-1)^{n-1+k(k+1)/2}\bigl(x^{n}\,d_{n-k-1}-x^{k+2-n}\,d_{n+k}\bigr)&&(n\ge k+1).
\end{align*}
The corresponding formula for the shift $1$ is a case of \cite[Theorem~4]{C24a} (see \S\ref{ss:related}). With these formulas the Desnanot--Jacobi identity for the Hankel matrices $(a_{i+j,k})_{0\le i,j\le n}$ becomes \eqref{eq:main}, which gives a third proof of Theorem~\ref{thm:B} (Corollary~\ref{cor:DJB}).

The curve is the curve of the Narayana weights, and Theorem~\ref{thm:B} can be restated without the parameter $c$ and with positive weights. Let $q$ be an indeterminate and let $\Nar_{n,j}(q)\in\Z[q]$ be the weighted number of Motzkin paths from $(0,0)$ to $(n,j)$ that never go below height $0$, where up steps have weight $1$, down steps have weight $q$, level steps at height $0$ have weight $1$, and level steps at positive height have weight $1+q$. The numbers $\Nar_{n,0}(q)$ are the Narayana polynomials $1,\ 1,\ 1+q,\ 1+3q+q^2,\dots$, and for every $k\ge0$ the sequence $\bigl(\Nar_{n+k,k}(q)\bigr)_{n\ge0}$ is the sequence of the $(2k+1)$-th convolution powers of the Narayana polynomials in the sense of Cigler \cite[Section~3]{C24a} (Corollary~\ref{cor:narayana}(i)).

\begin{corollary}[Narayana form of Theorem~\ref{thm:B}]\label{cor:narintro}
Let $k\ge0$ and $\fn_n:=\det\bigl(\Nar_{i+j,k}(q)\bigr)_{0\le i,j\le n-1}$ for $n\ge0$. Then for all $n\ge2k+2$
\[
\fn_n\,\fn_{n-2k-2}-q^{2k+2}\,\fn_{n-1}\,\fn_{n-2k-1}=(1-q)\,q^{(k+1)^2}\,\fn_{n-k-1}^{\,2}\qquad\text{in }\Z[q].
\]
\end{corollary}

This is the first identity of Corollary~\ref{cor:narayana}(ii). It is equivalent to Theorem~\ref{thm:B}, because $\Nar_{n,j}(x^2)=(-x)^{n-j}a_{n,j}(x,-x-x^{-1})$. Since $\Nar_{i+j,k}(q)$ is the coefficient of $z^{i+j-k}$ in the generating function of the $(2k+1)$-th convolution power of the Narayana polynomials in the sense of \cite[Section~3]{C24a}, $\fn_n$ is a Hankel determinant of this convolution power shifted backward by $k$, the determinant $D_{2k+1,-k,q}(n)$ of \cite{C24a}. The same recurrence holds for the shift $1-k$, but in general not for the unshifted Hankel determinants (\S\ref{ss:narayana}). For $q=1$ the numbers $\Nar_{n+k,k}(1)$ $(n\ge0)$ are the coefficients of the $(2k+1)$-th power of the generating function of the Catalan numbers, and the right-hand side vanishes.

For arbitrary $(x,c)$ we give formulas of Christoffel type for the shifts $1$ and $-1$ (Proposition~\ref{prop:shift1}, Corollary~\ref{cor:back1}). For $x=1$ they give closed formulas and prove the identities (32) and (33) of \cite{C23}, which are stated there as observations (Corollary~\ref{rem:C32}).

\emph{An analogue of odd orthogonal type and paths with a weighted return step $($Section~$\ref{sec:typeB})$.} Let $k\ge1$, let $\cE_n$ be the Hankel determinant of order $n$ of the $k$-th column of the triangle of the weights $s=(c+x-1,c,c,\dots)$, $\lambda=(1+x,1,1,\dots)$, and let $\cO_n=\det\bigl(a_{i+j,k}(x,c)+a_{i+j,k-1}(x,c)\bigr)_{0\le i,j\le n-1}$. These determinants do not depend on $c$, and for all $n\ge2k+1$ (Theorem~\ref{thm:typeB})
\begin{align*}
(-1)^k\bigl(\cE_n\,\cE_{n-2k-1}-x^2\,\cE_{n-1}\,\cE_{n-2k}\bigr)&=(1-x)\,\cO_{n-k-1}^{\,2},\\
(-1)^k\bigl(\cO_n\,\cO_{n-2k-1}-x^2\,\cO_{n-1}\,\cO_{n-2k}\bigr)&=(1-x)(1+x)^2\,\cE_{n-k}^{\,2}.
\end{align*}
Here the role of the polynomials $U_m$ is played by the Chebyshev polynomials $W_m=U_m+U_{m-1}$ of the fourth kind, and for $n>k$ the determinants $\cE_n$ and $\cO_n$ are, up to sign, rectangular characters of $\Sp(2k)$ and of $\SO(2k+1)$ (Remark~\ref{rem:typeB}). Let $D^{[u]}_{k,n}$ be the Hankel determinant of order $n$ of the $k$-th column of the triangle of the weights $s=(c,c,\dots)$, $\lambda=(u,1,1,\dots)$, that is, of Motzkin paths in which the down steps ending at height $0$ have weight $u$. It does not depend on $c$ (Theorem~\ref{thm:return}). If $k=2k'+1$, then $D^{[u]}_{k,2n+1}=0$ and $D^{[u]}_{k,2n}=(-1)^nD_{k',n}(u-1)^2$. If $k=2k'$, then $D^{[u]}_{k,n}=\cE^{(k')}_{\lceil n/2\rceil}(u-1)\,\cO^{(k')}_{\lfloor n/2\rfloor}(u-1)$, where $\cE^{(k')}_n(x)$ and $\cO^{(k')}_n(x)$ are the determinants above for $k'$ in place of $k$.

\subsection{Related work}\label{ss:related}
\emph{Cigler's conjectures.} The recurrence \eqref{eq:main} appears as eq.~(18) in \cite{C23} (September 2023). In the MathOverflow question \cite{CMO}, posted on 18 October 2023, Cigler conjectured that $D_{k,n}(x,c)$ is the determinant of an explicit matrix $H_{k,n}(x)$ with entries $0$, $1$ and powers of $x$; on 22 October 2023 he added the recurrence, with the range $n\ge2k+2$, to the question. The determinant formula was proved by Jandhyala \cite{J}, who states that the quadratic recurrence remains open, and independently in \cite[Theorem~A]{Yu}. Before \cite{Yu}, the recurrence was known for $k\le1$, for $x=0$ and for $x=1$. The last two cases follow from the evaluations \cite[Theorems~1.1 and~1.6]{CS} of $D_{k,n}(0,c)$ (which goes back to \cite{CK11}) and of $D_{k,n}(1,c)$. We refer to \cite[\S1.2]{Yu} for a detailed account. The formula of Christoffel type behind \eqref{eq:C1} is \cite[Theorem~1]{K23} (see \cite[Section~2]{K23} for its history). The cases of three and four consecutive moments are treated in \cite{CK21}, with closed forms in terms of Chebyshev polynomials when the weights do not depend on the height except at height $0$ \cite[Corollaries~3 and~6]{CK21}. Lemma~\ref{lem:christoffel} below is a version of the formula that does not refer to the roots of the polynomial.

\emph{Gale--Robinson recurrences and the Hirota--Miwa equation.} Theorem~\ref{thm:B} has the shape of a three-term bilinear recurrence of Gale--Robinson type \cite{Gale}, with span $2k+2$ and shifts $1$ and $k+1$. Such recurrences, with arbitrary coefficients, are special cases of the octahedron recurrence \cite[Example~1.7 and the proof of Theorem~1.8]{FZ}, \cite[\S1.3]{Sp}, which is another name for the Hirota--Miwa equation \cite{HKW}. It is known that the Hirota--Miwa equation \cite{Hir,Miw}, which is satisfied by Casorati determinants \cite{OHTI}, reduces to recurrences of Somos and Gale--Robinson type when the three lattice directions are identified with translations in one variable \cite{HKW}. Proposition~\ref{thm:reduction} gives a sufficient condition for a Casorati determinant solution to be compatible with such a reduction. We have not found it in this form in the literature. The Desnanot--Jacobi identity for the shifted Hankel determinants of a sequence is the discrete-time Toda equation \cite{Hir77} in the form \cite[(2.8)]{SNvdK}, where its relation to the quotient-difference algorithm is also discussed. Thus the second proof of Theorem~\ref{thm:B} uses the Hirota--Miwa equation for Casorati determinants of size $k$, and the third proof uses the discrete-time Toda equation for the determinants $D_{m,k,n}$.

For the Hankel determinants of the $q$-deformed metallic numbers, Ovsienko and Pedon \cite[Conjecture~1.9, eq.~(1.22)]{OP} conjectured a recurrence of the same shape (for the $m$-th metallic number, span $2m+2$, shifts $1$ and $m+1$, coefficients $1$ and $-1$). They proved it for $m=1,2$, and Han and Pedon \cite[Theorem~C]{HP} proved it in general by means of Hankel continued fractions. Their sequences are different from ours, and we do not know a direct relation between the two results. In Remark~\ref{rem:HP} we show that the recurrence of \cite{HP} follows, for the shifts $1\le\ell\le m$, from another theorem of \cite{HP} and the Desnanot--Jacobi identity, by the argument of our third proof. Hankel determinant solutions of the Somos-4 recurrence are treated in \cite{Xin,CHX}, of the Somos-5 recurrence in \cite{CHX}, and of Somos-type recurrences attached to hyperelliptic curves in \cite{Hon,CL26}. Conditions under which a Hankel transform is an $(\alpha,\beta)$ Somos-4 sequence are given in \cite{WZ24,LWZ26}.

\emph{Shifted determinants.} For $x=0$ the determinants $D_{1,k,n}(0,c)$ and $D_{2,k,n}(0,c)$ are evaluated in \cite{CK11} and \cite{CS}, respectively. For $x=1$ the determinants $D_{1,k,n}(1,c)$ are evaluated in \cite[Theorem~1.7]{CS}, and the determinants $D_{2,k,n}(1,c)$, except for two residue classes of $n$ modulo $2k+1$, in \cite[Theorem~1.8]{CS}. The evaluations of \cite{CS} quoted here were conjectured in \cite{C23} (Conjectures~5, 8 and~9 there). The duality between backward and forward shifts for $x=0$ is \cite[Theorem~2]{C24}, and it is also a case of a duality of Liu, Wang, Zhang and Zhang \cite[Theorem~2.4]{LWZZ} for powers of power series that satisfy a quadratic equation. At the point $(x,c)=(-1,2)$ the $k$-th column of the triangle consists of the coefficients of the $(2k+1)$-th power of the generating function of the Catalan numbers. The corresponding duality, formula~(31) of \cite{C23}, is contained in \cite[Conjecture~10]{C23a}, which was proved by Fulmek \cite{F} and by Cigler \cite[Theorem~2]{C24a}. Cigler extended it to convolution powers of Narayana polynomials \cite[Theorem~4]{C24a} (see also \cite{H25}). In our notation \cite[Theorem~4]{C24a}, with $k+1$ in place of the $k$ there, is the identity \eqref{eq:D4} of Theorem~\ref{thm:curve} on the curve $xc+x^2+1=0$, and it contains \eqref{eq:D1}. For $c=1$ the numbers $a_{n+k,k}(0,1)$ are the coefficients of the $(k+1)$-th power of the generating function of the Motzkin numbers. Hankel determinants of these powers, that is, the determinants $D_{k,k,n}(0,1)$, are studied in \cite{WZ25}, and the duality between their backward and forward shifted Hankel determinants is \cite[Corollary~3.2]{LWZZ}. Mu and Wang \cite{MW17} evaluate the Hankel determinants of the shifted sequence $(0,b_0,b_1,\dots)$ for Catalan-like numbers $(b_n)$ in the sense of Aigner, as described in \cite[\S2.2]{HY}.

Theorems~\ref{thm:B} and~\ref{thm:C} are due to Yu \cite{Yu}. Our proof of \eqref{eq:C2}, and hence of Theorem~\ref{thm:B}, is different from his. For \eqref{eq:C1} we follow \cite{Yu}, except that the relations \eqref{eq:introLambda} are verified directly in the quotient ring. The identities \eqref{eq:D1} and \eqref{eq:D4} are \cite[Theorem~4]{C24a}, Proposition~\ref{prop:backward} is \cite[Theorem~2]{C24} for $x=0$, and the lemmas of Section~\ref{sec:prelim}, Lemma~\ref{lem:HM} and Lemma~\ref{lem:reciprocal} are known in the forms indicated there. To the best of our knowledge the other results are new. These are the results of Section~\ref{sec:general}, the identities \eqref{eq:D2} and \eqref{eq:D3}, Proposition~\ref{prop:backward} for $x\ne0$ off the curve $xc+x^2+1=0$, Proposition~\ref{prop:shift1}, Corollary~\ref{cor:back1}, a proof of the identities (32) and (33) of \cite{C23}, the reading of Theorem~\ref{thm:B} as a recurrence for convolution powers of Narayana polynomials, and Theorems~\ref{thm:typeB} and~\ref{thm:return} (the proof of Theorem~\ref{thm:return} uses a classical factorization, see \S\ref{ss:u}).

\subsection{Organization and notation}\label{ss:notation}
Section~\ref{sec:prelim} introduces the transfer matrix of weighted Motzkin paths and the Chebyshev polynomials, and it collects the lemmas on determinants that are used throughout. Section~\ref{sec:HM} contains the Hirota--Miwa equation for Casorati determinants and the second proof of Theorems~\ref{thm:B} and~\ref{thm:C}. Section~\ref{sec:general} contains the generalizations of Theorem~\ref{thm:B}, Section~\ref{sec:shifted} the results on shifted determinants and the third proof, and Section~\ref{sec:typeB} the analogue of odd orthogonal type and the paths with a weighted return step. Section~\ref{sec:remarks} contains remarks on the recurrence of Han and Pedon, on solitons and on numerical observations. Appendix~\ref{app:C32} contains the proofs of the identities (32) and (33) of \cite{C23}, and Appendix~\ref{ss:checks} describes the computer checks. Apart from the lemma on complementary minors, two identities for Chebyshev polynomials and one generating function, which we quote from \cite{Yu}, and some further facts cited in remarks, the paper is self-contained. We use the notation of \cite{Yu} where possible.

\emph{General conventions.} We put $\N=\{0,1,2,\dots\}$ and denote by $\Z$ and $\Q$ the integers and the rationals. Throughout, $A$ denotes a commutative ring and $\kk$ a field. A binomial coefficient $\binom nj$ is $0$ unless $0\le j\le n$. We write $\delta_{i,j}$ for the Kronecker delta and $B^{\T}$ for the transpose of a matrix or vector $B$. All matrices and vectors are indexed from $0$. In a list, a hat marks an omitted member. Unless stated otherwise, the letters $x$, $c$, $y$, $u$, $q$ denote indeterminates, $z$ is the variable of power series, and $t$ and $\tau$ are variables of Laurent polynomials. The symbols that are used in more than one section are listed in Table~\ref{tab:notation}.

\begin{table}[ht]
\caption{Notation}\label{tab:notation}
\small
\begin{tabular}{@{}lp{0.66\textwidth}@{}}
\hline
$U_m$ $(m\in\Z)$ & Chebyshev polynomials: $U_{-1}=0$, $U_0=1$, $U_{m+1}=XU_m-U_{m-1}$ (\S\ref{ss:cheb})\\
$P_{m,y}$, $P_m$, $W_m$ & $P_{m,y}=U_m-yU_{m-1}$;\ \ $P_m=P_{m,x}$;\ \ $W_m=P_{m,-1}=U_m+U_{m-1}$\\
$J$, $a^{(s,\lambda)}_{n,j}$, $p_m$, $h_m$, $L$ & objects attached to weights $s=(s_j)_{j\ge0}$, $\lambda=(\lambda_j)_{j\ge1}$ (\S\ref{ss:matrices})\\
$L_y$ & the functional $L$ for $s=(y,0,0,\dots)$, $\lambda=(1,1,\dots)$; its polynomials are $P_{m,y}$\\
$a_{n,j}(x,c)$ & the case $s=(c+x,c,c,\dots)$, $\lambda=(1,1,\dots)$\\
$D_{k,n}(x)$, $d_n$ & $D_{k,n}(x,c)$, which does not depend on $c$; $d_n=D_{k,n}(x)$ for fixed $k$ (Section~\ref{sec:shifted})\\
$D_{m,k,n}(x,c)$ & $\det(a_{i+j+m,k}(x,c))_{0\le i,j\le n-1}$, $m\in\Z$ (Section~\ref{sec:shifted})\\
$\detx{\chi}$, $\detR$ & determinant of coordinates modulo a monic polynomial $\chi$ (\S\ref{ss:det}); $\detR=\detx{P_k}$\\
$\partial_a$, $\Cas_m$ & $(\partial_aF)_j=F_{j+1}-aF_j$;\ \ $\Cas_m(F)=\det(F_m,\dots,F_{m+r-1})$ (\S\ref{ss:HM})\\
$K$, $R$, $F_j$, $w(m)$ & $K=\Q(x)$, $R=K[X]/(P_k)$, $F_j$ the class of $U_{j-1}$, $w(m)=\Cas_m(F)$ (\S\ref{ss:Lambda})\\
$\Lambda$, $\eps_k$ & multiplication by $-P_{k-1}$ on $R$;\ \ $\eps_k=(-1)^{k(k+1)/2}=\det\Lambda$ (Lemma~\ref{lem:Lambda})\\
$\cW_m(y)$, $\cW^*_m(y)$, $M_i(m)$ & $\cW_m(y)=\Cas_m(\partial_yF)$, $\cW^*_m(y)=y^k\cW_m(1/y)$, coefficients $M_i(m)$ (\S\ref{ss:spectral})\\
$Q(y)$ & $y^{2k+2}-xy^{2k+1}+xy-1$ (\S\ref{ss:spectral})\\
$\cH_n$, $\cH_n(x;y)$ & Hankel determinants of $n\mapsto L(X^nP_k)$; the case $L=L_y$ (\S\ref{ss:general}, \S\ref{ss:spectral})\\
$\zeta$, $\phi_m(c)$ & $\zeta=x+x^{-1}$;\ \ $\phi_m(c)=U_{m-1}(c)$, the Fibonacci polynomials (Section~\ref{sec:shifted})\\
$\Nar_{n,j}(q)$, $\cN_{m,n}$, $\fn_n$ & Narayana weights and their Hankel determinants (\S\ref{ss:narayana})\\
$r_m$, $\cE_n$, $\cO_n$, $D^{[u]}_{k,n}$ & $r_m=W_m-xW_{m-1}$; determinants of Section~\ref{sec:typeB}\\
\hline
\end{tabular}
\end{table}
\section{Preliminaries}\label{sec:prelim}

\subsection{Infinite matrices and weighted Motzkin paths}\label{ss:matrices}
Let $A$ be a commutative ring. Let $M(A)$ be the set of matrices $B=(B_{i,j})_{i,j\in\N}$ with entries in $A$ such that every row and every column of $B$ has only finitely many non-zero entries. With the usual matrix product, $M(A)$ is an associative ring with unit $I$. Hence $f(B)\in M(A)$ is defined for $B\in M(A)$ and $f\in A[X]$, and polynomials in the same $B$ commute. Let $e_m$ $(m\in\N)$ be the $m$-th standard basis column vector, so that $Be_m$ is the $m$-th column of $B$, $e_i^{\T}B$ is its $i$-th row, and $e_i^{\T}Be_j=B_{i,j}$.

Let $s=(s_j)_{j\ge0}$ and $\lambda=(\lambda_j)_{j\ge1}$ be sequences of elements of $A$ (the \emph{weights}). Let $J=J^{(s,\lambda)}\in M(A)$ be the tridiagonal matrix with
\[
J_{m,m}=s_m,\qquad J_{m,m+1}=1,\qquad J_{m+1,m}=\lambda_{m+1}\qquad(m\ge0),
\]
all other entries being $0$, and put $a^{(s,\lambda)}_{n,j}:=(J^n)_{0,j}$ for $n,j\ge0$. Define monic polynomials $p_m=p^{(s,\lambda)}_m\in A[X]$ by
\[
p_0=1,\qquad p_1=X-s_0,\qquad p_{m+1}=(X-s_m)\,p_m-\lambda_m\,p_{m-1}\quad(m\ge1),
\]
let $L=L^{(s,\lambda)}\colon A[X]\to A$ be the $A$-linear map $L(f)=\bigl(f(J)\bigr)_{0,0}=e_0^{\T}f(J)e_0$, and put $h_0:=1$ and $h_m:=\lambda_1\lambda_2\cdots\lambda_m$ for $m\ge1$. For $\lambda=(1,1,\dots)$ these are the objects of \cite[Lemma~3.1]{Yu}. The following facts are standard. We include the short proof, which extends that of \cite[Lemma~3.1 and Corollary~3.2]{Yu} to weighted down steps.

\begin{lemma}\label{lem:jacobi}
With the above notation, for all $i,j,k,m,n\ge0$:
\begin{enumerate}
\item[(i)] $a^{(s,\lambda)}_{0,j}=\delta_{j,0}$ and, for $n\ge1$,
\[
a^{(s,\lambda)}_{n,j}=a^{(s,\lambda)}_{n-1,j-1}+s_j\,a^{(s,\lambda)}_{n-1,j}+\lambda_{j+1}\,a^{(s,\lambda)}_{n-1,j+1}\qquad(a^{(s,\lambda)}_{n-1,-1}:=0);
\]
moreover $a^{(s,\lambda)}_{n,j}=0$ for $j>n$ and $a^{(s,\lambda)}_{n,n}=1$.
\item[(ii)] $e_0^{\T}\,p_m(J)=e_m^{\T}$ and $p_m(J)\,e_0=h_m\,e_m$.
\item[(iii)] $h_k\,a^{(s,\lambda)}_{n,k}=L(X^np_k)$ and $L(p_ip_j)=h_j\,\delta_{i,j}$.
\item[(iv)] Let $c\in A$, let $s'=(s_j+c)_{j\ge0}$, and let $J',p'_m,L'$ be the objects defined as above for the weights $(s',\lambda)$. Then $J'=J+cI$, $p'_m(X)=p_m(X-c)$, and $L'(f)=L\bigl(f(X+c)\bigr)$ for all $f\in A[X]$. Moreover, for every $\chi\in A[X]$,
\[
\det\bigl(L((X+c)^{i+j}\chi)\bigr)_{0\le i,j\le n-1}=\det\bigl(L(X^{i+j}\chi)\bigr)_{0\le i,j\le n-1}.
\]
\end{enumerate}
\end{lemma}

By (i), $a^{(s,\lambda)}_{n,j}$ is the weighted number of Motzkin paths from $(0,0)$ to $(n,j)$ that never go below height $0$, in which up steps have weight $1$, level steps at height $j$ have weight $s_j$, and down steps from height $j$ to height $j-1$ have weight $\lambda_j$; in particular $L(X^n)=a^{(s,\lambda)}_{n,0}$. By (iii), $L$ is the moment functional for which the $p_m$ are the monic orthogonal polynomials.

\begin{proof}
(i) Column $j$ of $J$ has the entries $J_{j-1,j}=1$ (if $j\ge1$), $J_{j,j}=s_j$ and $J_{j+1,j}=\lambda_{j+1}$, so the recurrence is the expansion of $(J^n)_{0,j}=\sum_i(J^{n-1})_{0,i}J_{i,j}$. The triangularity follows by induction on $n$, since $j-1,j,j+1>n-1$ if $j>n$, and $a^{(s,\lambda)}_{n,n}=a^{(s,\lambda)}_{n-1,n-1}$.

(ii) The rows of $J$ give $e_0^{\T}J=s_0e_0^{\T}+e_1^{\T}$ and $e_m^{\T}J=\lambda_me_{m-1}^{\T}+s_me_m^{\T}+e_{m+1}^{\T}$ for $m\ge1$. Hence $e_0^{\T}p_1(J)=e_1^{\T}$ and, inductively for $m\ge1$, $e_0^{\T}p_{m+1}(J)=e_m^{\T}(J-s_mI)-\lambda_me_{m-1}^{\T}=e_{m+1}^{\T}$. The columns of $J$ give $Je_0=s_0e_0+\lambda_1e_1$ and $Je_m=e_{m-1}+s_me_m+\lambda_{m+1}e_{m+1}$ for $m\ge1$. Hence $p_1(J)e_0=h_1e_1$ and, inductively for $m\ge1$ (using $\lambda_mh_{m-1}=h_m$), $p_{m+1}(J)e_0=(J-s_mI)h_me_m-\lambda_mh_{m-1}e_{m-1}=h_m(e_{m-1}+\lambda_{m+1}e_{m+1})-h_me_{m-1}=h_{m+1}e_{m+1}$.

(iii) By (ii),
\begin{align*}
h_k\,a^{(s,\lambda)}_{n,k}&=e_0^{\T}J^n(h_ke_k)=e_0^{\T}J^np_k(J)e_0=L(X^np_k),\\
L(p_ip_j)&=e_0^{\T}p_i(J)\,p_j(J)e_0=e_i^{\T}(h_je_j)=h_j\delta_{i,j}.
\end{align*}

(iv) Clearly $J'=J+cI$, and $p'_m(X)=p_m(X-c)$ follows by induction from the defining recurrences. For $f\in A[X]$ we have $f(J')=g(J)$ with $g(X)=f(X+c)$, so $L'(f)=L(f(X+c))$. Finally
\[
L\bigl((X+c)^{i+j}\chi\bigr)=\sum_{l,l'}\binom il c^{\,i-l}\,L(X^{l+l'}\chi)\,\binom j{l'}c^{\,j-l'},
\]
that is, the matrix on the left of the asserted identity is $T\cdot(L(X^{l+l'}\chi))_{l,l'<n}\cdot T^{\T}$ with the lower unitriangular matrix $T=\bigl(\binom il c^{\,i-l}\bigr)_{0\le i,l\le n-1}$.
\end{proof}

Two families of weights occur throughout, both with $\lambda=(1,1,\dots)$. For elements $y,c$ of $A$ we write $a_{n,j}(y,c)$ for the numbers $a^{(s,\lambda)}_{n,j}$ of the weights $s=(c+y,c,c,\dots)$. By Lemma~\ref{lem:jacobi}(i) this agrees with the definition of $a_{n,j}(x,c)$ in \S\ref{ss:cigler}. We write $L_y$ for the functional $L^{(s,\lambda)}$ of the weights $s=(y,0,0,\dots)$, whose polynomials $p_m$ are determined in Lemma~\ref{lem:cheb}(iii).

\subsection{Chebyshev polynomials}\label{ss:cheb}
Let $U_m\in\Z[X]$ $(m\in\Z)$ be defined, as in Theorem~\ref{thm:C}, by $U_{-1}=0$, $U_0=1$ and the two-sided recurrence $U_{m+1}=XU_m-U_{m-1}$ $(m\in\Z)$. Thus $U_1=X$, $U_{-2}=-1$, and $U_m(2\cos\theta)=\sin((m+1)\theta)/\sin\theta$. For an element $y$ of a commutative ring $A$ and $m\in\Z$ put
\[
P_{m,y}:=U_m-y\,U_{m-1}\in A[X],
\]
so that $P_{m,0}=U_m$. For $m\ge0$ these are co-recursive Chebyshev polynomials in the sense of \cite{Chi57}, that is, only the first coefficient of their recurrence differs from that of the $U_m$ (Lemma~\ref{lem:cheb}(iii)). Two cases have their own names. For the indeterminate $x$ we write $P_m:=P_{m,x}=U_m-xU_{m-1}$, as in Theorem~\ref{thm:C}, and we put $W_m:=P_{m,-1}=U_m+U_{m-1}$. Then $W_m(2\cos\theta)=\sin((m+\frac12)\theta)/\sin(\frac\theta2)$, so that $U_m(2X)$ and $W_m(2X)$ are the Chebyshev polynomials of the second and the fourth kind \cite[(18.5.2), (18.5.4)]{DLMF}.

\begin{lemma}\label{lem:cheb}
\begin{enumerate}
\item[(i)] In $\Z[t,t^{-1}]$ we have, for all $m\in\Z$,
\begin{equation}\label{eq:cheb}
(t-t^{-1})\,U_m(t+t^{-1})=t^{m+1}-t^{-m-1}.
\end{equation}
\item[(ii)] $U_{-m-2}=-U_m$ and $P_{m,y}\bigl(y+y^{-1}\bigr)=y^{-m}$ for all $m\in\Z$, the latter in $\Z[y,y^{-1}]$ for an indeterminate $y$.
\item[(iii)] $P_{0,y}=1$, $P_{1,y}=X-y$, and $P_{m+1,y}=XP_{m,y}-P_{m-1,y}$ for all $m\in\Z$. Consequently $(P_{m,y})_{m\ge0}$ are the polynomials $p_m$ of the weights $s=(y,0,0,\dots)$, $\lambda=(1,1,\dots)$, and for all $c\in A$ and $j,n\ge0$
\begin{equation}\label{eq:shiftm}
a_{n,j}(y,c)=L_y\bigl((X+c)^nP_{j,y}\bigr).
\end{equation}
\end{enumerate}
\end{lemma}

\begin{proof}
The identity {\rm(i)} and the first identity of {\rm(ii)} are \cite[Lemma~2.2]{Yu}. For the second one, by \eqref{eq:cheb} with $t=y$, $(y-y^{-1})P_{m,y}(y+y^{-1})=(y^{m+1}-y^{-m-1})-y\,(y^m-y^{-m})=y^{-m}(y-y^{-1})$.

(iii) The first assertions are clear from the recurrence of the $U_m$; they are the recurrence defining the $p_m$ for the weights $(y,0,0,\dots)$, $(1,1,\dots)$. For \eqref{eq:shiftm} apply Lemma~\ref{lem:jacobi}(iii),(iv) to these weights: the shifted weights are $(c+y,c,c,\dots)$, with $p'_j(X)=P_{j,y}(X-c)$, so $a_{n,j}(y,c)=L'(X^np'_j)=L_y\bigl((X+c)^nP_{j,y}\bigr)$.
\end{proof}

\begin{corollary}\label{cor:cindep}
For all $k,m,n\ge0$,
\begin{equation}\label{eq:Dshift}
\det\bigl(a_{i+j+m,k}(x,c)\bigr)_{0\le i,j\le n-1}=\det\bigl(L_x(X^{i+j}(X+c)^mP_k)\bigr)_{0\le i,j\le n-1}.
\end{equation}
In particular $D_{k,n}(x,c)=\det\bigl(L_x(X^{i+j}P_k)\bigr)_{0\le i,j\le n-1}$ does not depend on $c$ {\rm\cite[Corollary~3.2]{Yu}}. Moreover $D_{k,n}=1$ if $k=0$, and $D_{k,n}=0$ for $1\le n\le k$ if $k\ge1$.
\end{corollary}

\begin{proof}
By \eqref{eq:shiftm}, $a_{i+j+m,k}(x,c)=L_x\bigl((X+c)^{i+j}\chi\bigr)$ with $\chi=(X+c)^mP_k$, and Lemma~\ref{lem:jacobi}(iv) gives \eqref{eq:Dshift}. For $k=0$ the matrix $\bigl(L_x(X^{i+j})\bigr)_{i,j<n}$ is transformed into $\bigl(L_x(P_iP_j)\bigr)_{i,j<n}=I_n$ by a unitriangular change of basis applied to rows and columns (Lemma~\ref{lem:jacobi}(iii)). If $1\le n\le k$, the first row of the matrix $(a_{i+j,k})_{i,j<n}$ consists of the elements $a_{j,k}$ with $j<k$, which vanish by Lemma~\ref{lem:jacobi}(i).
\end{proof}

\subsection{Determinants modulo a monic polynomial}\label{ss:det}
Let $\kk$ be a field and let $\chi\in\kk[X]$ be monic of degree $d\ge1$. A \emph{monic basis} of $\kk[X]/(\chi)$ is a basis formed by the residue classes of monic polynomials $\beta_0,\dots,\beta_{d-1}$ with $\deg\beta_i=i$. For $f_1,\dots,f_d\in\kk[X]$ we denote by
\[
\detx{\chi}(f_1,\dots,f_d)
\]
the determinant of the $d\times d$ matrix whose $i$-th row is the coordinate vector of the residue class of $f_i$ with respect to a monic basis. Two monic bases differ by a unitriangular matrix, so this determinant does not depend on the choice of the monic basis; in particular $\detx{\chi}(\beta_0,\dots,\beta_{d-1})=1$ for every monic basis. We use the same notation for elements of $\kk[X]/(\chi)$. If $\Psi$ is a $\kk$-linear endomorphism of $\kk[X]/(\chi)$, then $\detx{\chi}(\Psi f_1,\dots,\Psi f_d)=\det(\Psi)\,\detx{\chi}(f_1,\dots,f_d)$.

\begin{lemma}[complementary minors {\cite[Lemma~3.3]{Yu}}]\label{lem:compl}
Let $r,r'\ge1$ be integers, $M\in\kk^{r\times(r+r')}$ and $Z\in\kk^{(r+r')\times r'}$ with $MZ=0$. For $S\subseteq\{0,\dots,r+r'-1\}$ with $|S|=r$ let $M_S$ be the $r\times r$ submatrix of $M$ formed by the columns with index in $S$, let $Z_S$ be the $r'\times r'$ submatrix of $Z$ formed by the rows with index not in $S$ $($rows and columns kept in increasing order$)$, and let $\|S\|$ denote the sum of the elements of $S$. Then for all such $S,S'$:
\[
(-1)^{\|S\|}\det(M_S)\det(Z_{S'})=(-1)^{\|S'\|}\det(M_{S'})\det(Z_S).
\]
\end{lemma}

This is Jacobi's theorem on complementary minors, that is, the duality between the Pl\"ucker coordinates of the row space of $M$ and those of its annihilator. In \cite{Yu} it is stated for the field $\Q(x)$, but the proof given there works over any field.

The next lemma is the formula of Christoffel type mentioned in \S\ref{ss:related}, in a form that does not refer to the roots of $f$ (compare \cite[Theorem~1]{K23}). For the weights $(x,0,0,\dots)$, $(1,1,\dots)$ and $f=P_k$ it is \cite[Lemma~3.4]{Yu}, and the proof is the same.

\begin{lemma}[Christoffel-type formula]\label{lem:christoffel}
Let $\cL\colon\kk[X]\to\kk$ be a $\kk$-linear map and $(q_m)_{m\ge0}$ monic polynomials with $\deg q_m=m$ and $\cL(q_iq_j)=\kappa_i\delta_{i,j}$ for all $i,j\ge0$, where $\kappa_i\in\kk$. Let $f\in\kk[X]$ be monic of degree $d\ge1$. Then for all $n\ge0$
\[
\det\bigl(\cL(X^{i+j}f)\bigr)_{0\le i,j\le n-1}=(-1)^{nd}\,\kappa_0\kappa_1\cdots\kappa_{n-1}\,\detx{f}\bigl(q_n,q_{n+1},\dots,q_{n+d-1}\bigr).
\]
In particular, for weights $(s,\lambda)$ with values in $\kk$,
\begin{equation}\label{eq:christoffel}
\det\bigl(L^{(s,\lambda)}(X^{i+j}f)\bigr)_{0\le i,j\le n-1}=(-1)^{nd}\,h_0h_1\cdots h_{n-1}\,\detx{f}\bigl(p_n,p_{n+1},\dots,p_{n+d-1}\bigr).
\end{equation}
\end{lemma}

\begin{proof}
For $n=0$ both sides are $1$; let $n\ge1$. The $q_m$ form a basis of $\kk[X]$, and $fq_j$ is monic of degree $j+d$; so for $j\ge0$ we can write $fq_j=\sum_{m=0}^{j+d}\Pi_{j,m}q_m$ with $\Pi_{j,m}\in\kk$ and $\Pi_{j,j+d}=1$, and we put $\Pi_{j,m}:=0$ for $m>j+d$. Multiplying by $q_i$ and applying $\cL$ gives $\cL(q_iq_jf)=\kappa_i\Pi_{j,i}$ for all $i,j\ge0$. Write $q_i=\sum_{l\le i}\xi_{i,l}X^l$ with $\xi_{i,i}=1$, and let $\Xi=(\xi_{i,l})_{0\le i,l\le n-1}$, a lower unitriangular matrix. Then $\bigl(\cL(q_iq_jf)\bigr)_{i,j<n}=\Xi\,\bigl(\cL(X^{l+l'}f)\bigr)_{l,l'<n}\,\Xi^{\T}$, and therefore
\[
\det\bigl(\cL(X^{i+j}f)\bigr)_{i,j<n}=\det\bigl(\kappa_i\Pi_{j,i}\bigr)_{i,j<n}=\kappa_0\cdots\kappa_{n-1}\,\det\bigl(\Pi_{j,m}\bigr)_{0\le j,m\le n-1}.
\]
It remains to show that $\det(\Pi_{j,m})_{j,m<n}=(-1)^{nd}\detx{f}(q_n,\dots,q_{n+d-1})$. For $g\in\kk[X]$ let $\crd(g)\in\kk^d$ be the coordinate row vector of the residue class of $g$ modulo $f$ with respect to the monic basis $q_0,\dots,q_{d-1}$. Since $\crd(fq_j)=0$, we have $\sum_{m=0}^{j+d}\Pi_{j,m}\crd(q_m)=0$ for all $j\ge0$. Put
\[
M=(\Pi_{j,m})_{0\le j\le n-1,\ 0\le m\le n+d-1}\in\kk^{n\times(n+d)},
\]
and let $Z$ be the $(n+d)\times d$ matrix with rows $\crd(q_0),\dots,\crd(q_{n+d-1})$. Then $MZ=0$ (for $j\le n-1$ all $m\le j+d$ are $\le n+d-1$). Apply Lemma~\ref{lem:compl} with $r=n$, $r'=d$, $S=\{0,\dots,n-1\}$ and $S'=\{d,\dots,n+d-1\}$. Here $\det M_S=\det(\Pi_{j,m})_{j,m<n}$; the matrix $Z_S$ consists of the rows $\crd(q_n),\dots,\crd(q_{n+d-1})$, so $\det Z_S=\detx{f}(q_n,\dots,q_{n+d-1})$; the matrix $M_{S'}=(\Pi_{j,d+i})_{0\le j,i\le n-1}$ is lower unitriangular, so $\det M_{S'}=1$; $Z_{S'}=I_d$; and $\|S\|=n(n-1)/2$, $\|S'\|=nd+n(n-1)/2$. Lemma~\ref{lem:compl} gives $(-1)^{n(n-1)/2}\det M_S=(-1)^{nd+n(n-1)/2}\det Z_S$. Finally, \eqref{eq:christoffel} is the case $\cL=L^{(s,\lambda)}$, $q_m=p_m$, $\kappa_m=h_m$, by Lemma~\ref{lem:jacobi}(iii).
\end{proof}

The last lemma of this section is Claim~1b in the proof of \cite[Proposition~3.7]{Yu}, stated there for particular vectors.

\begin{lemma}\label{lem:expansion}
Let $E$ be a $\kk$-vector space of finite dimension $r\ge1$ with a fixed basis, and for $u_1,\dots,u_r\in E$ let $\detE(u_1,\dots,u_r)$ be the determinant of the matrix whose $i$-th row is the coordinate vector of $u_i$. Let $u_0,\dots,u_r\in E$, let $z\in\kk$, and put $m_i:=\detE(u_0,\dots,\widehat{u_i},\dots,u_r)$ for $0\le i\le r$. Then
\begin{align*}
\detE\bigl(u_1-zu_0,\ u_2-zu_1,\ \dots,\ u_r-zu_{r-1}\bigr)&=\sum_{i=0}^r(-z)^i\,m_i,\\
\detE\bigl(u_0-zu_1,\ u_1-zu_2,\ \dots,\ u_{r-1}-zu_r\bigr)&=\sum_{i=0}^r(-z)^{r-i}\,m_i.
\end{align*}
\end{lemma}

\begin{proof}
Expand the first determinant by multilinearity into $2^r$ terms, one for each choice, in every argument $u_j-zu_{j-1}$ $(1\le j\le r)$, of its first or its second summand. If the argument $j$ chooses $u_j$ and the argument $j+1$ chooses $-zu_j$, the term vanishes. In the remaining terms there is $i\in\{0,\dots,r\}$ such that the arguments $j\le i$ choose $-zu_{j-1}$ and the arguments $j>i$ choose $u_j$. The chosen vectors are then $u_0,\dots,u_{i-1},u_{i+1},\dots,u_r$ in increasing order, and the term equals $(-z)^im_i$. In the second determinant, with the arguments $u_j-zu_{j+1}$ $(0\le j\le r-1)$, the remaining terms are those in which the arguments $j<i$ choose $u_j$ and the arguments $j\ge i$ choose $-zu_{j+1}$, and such a term equals $(-z)^{r-i}m_i$.
\end{proof}
\section{The Hirota--Miwa equation and Cigler's recurrence}\label{sec:HM}

\subsection{Casoratians and the Hirota--Miwa equation}\label{ss:HM}
In this subsection $\kk$ is an arbitrary field and $E$ is a $\kk$-vector space of finite dimension $r\ge1$ with a fixed basis; $\detE$ is defined as in Lemma~\ref{lem:expansion}. A \emph{sequence in $E$} is a family $F=(F_j)_{j\in\Z}$ of elements of $E$. For a sequence $F$ in $E$, $a\in\kk$ and $m\in\Z$ put
\[
(\partial_aF)_j:=F_{j+1}-aF_j\quad(j\in\Z),\qquad \Cas_m(F):=\detE(F_m,F_{m+1},\dots,F_{m+r-1}).
\]
The sequence $\partial_0F$ is $F$ translated by one step, so that $\Cas_m(\partial_0F)=\Cas_{m+1}(F)$, and $\partial_a\partial_b=\partial_b\partial_a$, with $(\partial_a\partial_bF)_j=F_{j+2}-(a+b)F_{j+1}+ab\,F_j$. The determinants $\Cas_m(F)$ are called Casorati determinants, or Casoratians.

\begin{lemma}[Hirota--Miwa equation for Casoratians]\label{lem:HM}
Let $F$ be a sequence in $E$, let $a_1,a_2,a_3\in\kk$ and $m\in\Z$. Then
\begin{equation}\label{eq:HM}
\begin{aligned}
&(a_2-a_3)\,\Cas_m(\partial_{a_1}F)\,\Cas_m(\partial_{a_2}\partial_{a_3}F)
+(a_3-a_1)\,\Cas_m(\partial_{a_2}F)\,\Cas_m(\partial_{a_1}\partial_{a_3}F)\\
&\qquad+(a_1-a_2)\,\Cas_m(\partial_{a_3}F)\,\Cas_m(\partial_{a_1}\partial_{a_2}F)=0.
\end{aligned}
\end{equation}
\end{lemma}

This is the bilinear difference equation of Hirota and Miwa \cite{Hir,Miw} in the form satisfied by Casorati determinants \cite{OHTI}, and it is a three-term Pl\"ucker relation. Both sides of \eqref{eq:HM} are polynomials with integer coefficients in $a_1,a_2,a_3$ and in the coordinates of $F_m,\dots,F_{m+r+1}$, so it would suffice to prove the identity over $\Q$ for vectors in general position. The short proof below works directly for all sequences over any field.

\begin{proof}
Replacing $F$ by the sequence $(F_{j+m})_{j\in\Z}$ we may assume $m=0$. Only the vectors $v_i:=F_i$, $0\le i\le r+1$, occur in \eqref{eq:HM}. Let $C$ be the $r\times(r+2)$ matrix whose $i$-th column is the coordinate vector of $v_i$ $(0\le i\le r+1)$. For row vectors $\mathbf u,\mathbf u'\in\kk^{r+2}$ let
\[
[\mathbf u,\mathbf u']:=\det\begin{pmatrix}C\\ \mathbf u\\ \mathbf u'\end{pmatrix}
\]
be the determinant of the $(r+2)\times(r+2)$ matrix obtained by writing $\mathbf u$ and $\mathbf u'$, in this order, below $C$. (Recall that $\detE$ of a family of vectors is the determinant of the matrix having their coordinate vectors as rows, or equivalently as columns.)

\emph{Step 1.} For all $\mathbf u_1,\mathbf u_2,\mathbf u_3,\mathbf u_4\in\kk^{r+2}$,
\begin{equation}\label{eq:pl3}
[\mathbf u_1,\mathbf u_2]\,[\mathbf u_3,\mathbf u_4]-[\mathbf u_1,\mathbf u_3]\,[\mathbf u_2,\mathbf u_4]+[\mathbf u_1,\mathbf u_4]\,[\mathbf u_2,\mathbf u_3]=0.
\end{equation}
Indeed, $(\mathbf u,\mathbf u')\mapsto[\mathbf u,\mathbf u']$ is bilinear and alternating, and $[\mathbf u,\mathbf u']=0$ whenever $\mathbf u$ or $\mathbf u'$ lies in the row space $\mathcal R\subseteq\kk^{r+2}$ of $C$. If the rows of $C$ are linearly dependent, then $[\mathbf u,\mathbf u']=0$ for all $\mathbf u,\mathbf u'$ and \eqref{eq:pl3} holds. Otherwise choose $\boldsymbol\zeta,\boldsymbol\zeta'\in\kk^{r+2}$ such that the rows of $C$ together with $\boldsymbol\zeta,\boldsymbol\zeta'$ form a basis of $\kk^{r+2}$, and write $\mathbf u_i=\mathbf c_i+\xi_i\boldsymbol\zeta+\eta_i\boldsymbol\zeta'$ with $\mathbf c_i\in\mathcal R$ and $\xi_i,\eta_i\in\kk$. Then $[\mathbf u_i,\mathbf u_j]=(\xi_i\eta_j-\xi_j\eta_i)\,[\boldsymbol\zeta,\boldsymbol\zeta']$, and \eqref{eq:pl3} follows from the identity
\[
(\xi_1\eta_2-\xi_2\eta_1)(\xi_3\eta_4-\xi_4\eta_3)-(\xi_1\eta_3-\xi_3\eta_1)(\xi_2\eta_4-\xi_4\eta_2)+(\xi_1\eta_4-\xi_4\eta_1)(\xi_2\eta_3-\xi_3\eta_2)=0,
\]
which is easily verified by expanding.

\emph{Step 2.} For $a\in\kk$ put $\mathbf e(a):=(1,a,a^2,\dots,a^{r+1})\in\kk^{r+2}$, and let $\mathbf e_\infty:=(0,\dots,0,1)\in\kk^{r+2}$. We claim that for all $a,b\in\kk$
\begin{equation}\label{eq:eval}
[\mathbf e(a),\mathbf e_\infty]=(-1)^r\,\Cas_0(\partial_aF),\qquad [\mathbf e(a),\mathbf e(b)]=(b-a)\,\Cas_0(\partial_a\partial_bF).
\end{equation}
For $n\ge2$ and $a\in\kk$ let $N^{(n)}_a$ be the $n\times n$ matrix with entries $1$ on the diagonal, $-a$ in the positions $(i,i+1)$, $0\le i\le n-2$, and $0$ elsewhere. Then $\det N^{(n)}_a=1$, and multiplying a matrix with columns $\mathbf c_0,\dots,\mathbf c_{n-1}$ from the right by $N^{(n)}_a$ replaces its columns by $\mathbf c_0,\ \mathbf c_1-a\mathbf c_0,\ \dots,\ \mathbf c_{n-1}-a\mathbf c_{n-2}$. In particular
\[
(1,b,\dots,b^{n-1})\,N^{(n)}_a=\bigl(1,\ b-a,\ (b-a)b,\ \dots,\ (b-a)b^{n-2}\bigr),\qquad (0,\dots,0,1)\,N^{(n)}_a=(0,\dots,0,1),
\]
and the first of these vectors is $(1,0,\dots,0)$ if $b=a$.

Fix $a$, put $v'_i:=v_{i+1}-av_i=(\partial_aF)_i$ for $0\le i\le r$, and let $C'$ be the $r\times(r+1)$ matrix whose columns are the coordinate vectors of $v'_0,\dots,v'_r$. Let $\mathbf u\in\kk^{r+2}$ and write $\mathbf u\,N^{(r+2)}_a=(u^\circ,\hat{\mathbf u})$ with $u^\circ\in\kk$ and $\hat{\mathbf u}\in\kk^{r+1}$. Then
\[
[\mathbf e(a),\mathbf u]=\det\Bigl(\begin{pmatrix}C\\ \mathbf e(a)\\ \mathbf u\end{pmatrix}N^{(r+2)}_a\Bigr)=\det\begin{pmatrix}v_0&C'\\ 1&0\ \cdots\ 0\\ u^\circ&\hat{\mathbf u}\end{pmatrix}=(-1)^r\det\begin{pmatrix}C'\\ \hat{\mathbf u}\end{pmatrix},
\]
by expansion along the row with index $r$, whose only non-zero entry is the entry $1$ in column $0$. For $\mathbf u=\mathbf e_\infty$ we have $\hat{\mathbf u}=(0,\dots,0,1)\in\kk^{r+1}$, and expansion along the last row gives $\det\binom{C'}{\hat{\mathbf u}}=\detE(v'_0,\dots,v'_{r-1})=\Cas_0(\partial_aF)$; this is the first formula of \eqref{eq:eval}. For $\mathbf u=\mathbf e(b)$ we have $\hat{\mathbf u}=(b-a)(1,b,\dots,b^r)$. Multiplying the $(r+1)\times(r+1)$ matrix with the rows of $C'$ and the last row $(1,b,\dots,b^r)$ from the right by $N^{(r+1)}_b$ turns the columns of $C'$ into $v'_0,\ v'_1-bv'_0,\ \dots,\ v'_r-bv'_{r-1}$ and the last row into $(1,0,\dots,0)$; expanding along the last row (index $r$, non-zero entry in column $0$) gives
\[
\det\begin{pmatrix}C'\\ (1,b,\dots,b^r)\end{pmatrix}=(-1)^r\detE\bigl(v'_1-bv'_0,\dots,v'_r-bv'_{r-1}\bigr)=(-1)^r\,\Cas_0(\partial_b\partial_aF).
\]
Hence $[\mathbf e(a),\mathbf e(b)]=(-1)^r(b-a)(-1)^r\Cas_0(\partial_a\partial_bF)$, which is the second formula of \eqref{eq:eval}.

\emph{Step 3.} Apply \eqref{eq:pl3} to $\mathbf u_1=\mathbf e(a_1)$, $\mathbf u_2=\mathbf e(a_2)$, $\mathbf u_3=\mathbf e(a_3)$, $\mathbf u_4=\mathbf e_\infty$ and insert \eqref{eq:eval}:
\begin{align*}
0={}&(a_2-a_1)\Cas_0(\partial_{a_1}\partial_{a_2}F)\cdot(-1)^r\Cas_0(\partial_{a_3}F)-(a_3-a_1)\Cas_0(\partial_{a_1}\partial_{a_3}F)\cdot(-1)^r\Cas_0(\partial_{a_2}F)\\
&+(-1)^r\Cas_0(\partial_{a_1}F)\cdot(a_3-a_2)\Cas_0(\partial_{a_2}\partial_{a_3}F).
\end{align*}
Multiplying by $(-1)^{r+1}$ gives \eqref{eq:HM} for $m=0$.
\end{proof}

\subsection{The quotient ring and the operator \texorpdfstring{$\Lambda$}{Lambda}}\label{ss:Lambda}
In this subsection $\kk$ is a field, $x$ is an element of $\kk$, $k\ge1$, $P_m=U_m-xU_{m-1}\in\kk[X]$ for $m\in\Z$, and $R:=\kk[X]/(P_k)$, a $\kk$-vector space of dimension $k$. For $j\in\Z$ let $F_j\in R$ be the residue class of $U_{j-1}$. The classes $F_1,\dots,F_k$ of $U_0,\dots,U_{k-1}$ form a monic basis of $R$; we write $\detR$ for $\detx{P_k}$ and put
\[
w(m):=\Cas_m(F)=\detR(F_m,F_{m+1},\dots,F_{m+k-1})\qquad(m\in\Z).
\]
By the recurrence of the $U_m$ and by Lemma~\ref{lem:cheb}(ii),
\begin{equation}\label{eq:Frec}
F_{j+1}=XF_j-F_{j-1},\qquad F_0=0,\qquad F_{-j}=-F_j\qquad(j\in\Z).
\end{equation}
In the situation of Theorem~\ref{thm:C} we have $\kk=K:=\Q(x)$, the field of rational functions in the indeterminate $x$, and $w(m)$ is the quantity defined there. Parts (i) and (iii) of the following lemma are \cite[Lemma~3.6(ii),(iii)]{Yu}, and the description of $\Lambda$ is the first part of \cite[Remark~4.3]{Yu}. Part (ii) is not stated there, but it follows from the second part of that remark. In \cite{Yu} the relations ($\Lambda1$), ($\Lambda2$) are obtained from a model of $R$ by antisymmetric Laurent polynomials. We give a direct proof in $R$, valid over every field.

\begin{lemma}\label{lem:Lambda}
Let $\Lambda\colon R\to R$ be the multiplication by the residue class of $-P_{k-1}$.
\begin{enumerate}
\item[(i)] For all $j\in\Z$,
\[
(\Lambda1)\quad \Lambda F_j=F_{j+k+1}-xF_{j+k},\qquad\qquad (\Lambda2)\quad \Lambda F_j=F_{j-k-1}-xF_{j-k}.
\]
\item[(ii)] $\Lambda^2$ is the multiplication by $1+x^2-xX$.
\item[(iii)] $\det\Lambda=\eps_k:=(-1)^{k(k+1)/2}$; in particular $\Lambda$ is invertible, and the class of $P_{k-1}$ is a unit of $R$.
\end{enumerate}
\end{lemma}

\begin{proof}
(i) Each of the three sequences $(\Lambda F_j)_j$, $(F_{j+k+1}-xF_{j+k})_j$ and $(F_{j-k-1}-xF_{j-k})_j$ satisfies the recurrence $f_{j+1}=Xf_j-f_{j-1}$ in $R$, because $(F_j)$ does and $\Lambda$ commutes with the multiplication by $X$. So it suffices to compare them for $j=0$ and $j=1$. For $j=0$: $\Lambda F_0=0$; $F_{k+1}-xF_k$ is the class of $U_k-xU_{k-1}=P_k$, which is $0$; and $F_{-k-1}-xF_{-k}=-(F_{k+1}-xF_k)=0$ by \eqref{eq:Frec}. For $j=1$: $\Lambda F_1$ is the class of $-P_{k-1}$; $F_{k+2}-xF_{k+1}$ is the class of $P_{k+1}=XP_k-P_{k-1}$ (Lemma~\ref{lem:cheb}(iii)), which is the class of $-P_{k-1}$; and $F_{-k}-xF_{-k+1}=-(F_k-xF_{k-1})$ is the class of $-P_{k-1}$.

(ii) By (i), first in the form ($\Lambda1$) and then in the form ($\Lambda2$),
\begin{align*}
\Lambda^2F_j&=\Lambda F_{j+k+1}-x\Lambda F_{j+k}=(F_j-xF_{j+1})-x(F_{j-1}-xF_j)\\
&=(1+x^2)F_j-x(F_{j+1}+F_{j-1})=(1+x^2-xX)F_j
\end{align*}
for all $j$, and $F_1,\dots,F_k$ is a basis.

(iii) For $m\in\Z$ put $G_m:=F_{m+1}-xF_m=(\partial_xF)_m$, the class of $P_m$. By ($\Lambda1$) with $j=m-k$,
\begin{equation}\label{eq:GLambda}
G_m=\Lambda F_{m-k}\qquad(m\in\Z).
\end{equation}
The classes $G_0,\dots,G_{k-1}$ of $P_0,\dots,P_{k-1}$ form a monic basis, so
\begin{align*}
1&=\detR(G_0,\dots,G_{k-1})=\detR(\Lambda F_{-k},\dots,\Lambda F_{-1})=\det(\Lambda)\cdot\detR(F_{-k},\dots,F_{-1})\\
&=\det(\Lambda)\cdot\detR(-F_k,-F_{k-1},\dots,-F_1)=\det(\Lambda)\cdot(-1)^k(-1)^{k(k-1)/2}\detR(F_1,\dots,F_k)\\
&=\det(\Lambda)\cdot(-1)^{k(k+1)/2};
\end{align*}
here we used \eqref{eq:Frec}, and reversing the order of $k$ rows has the sign $(-1)^{k(k-1)/2}$.
\end{proof}

In the language of \S\ref{ss:HM}, \eqref{eq:GLambda} and ($\Lambda2$) say that, up to the map $\Lambda$ and a scalar, the operation $\partial_x$ is the translation of the index by $-k$, and, if $x\ne0$, the operation $\partial_{x^{-1}}$ is the translation by $k+1$:
\begin{equation}\label{eq:shifts}
(\partial_xF)_j=\Lambda F_{j-k},\qquad(\partial_{x^{-1}}F)_j=-x^{-1}\bigl(F_j-xF_{j+1}\bigr)=-x^{-1}\Lambda F_{j+k+1}\qquad(j\in\Z).
\end{equation}

\begin{lemma}[{\cite[Lemma~3.4]{Yu}}]\label{lem:reduction}
Let $\kk=K=\Q(x)$. For all $n\ge0$,
\[
D_{k,n}(x,c)=(-1)^{nk}\detR(P_n,P_{n+1},\dots,P_{n+k-1})=(-1)^{nk}\,\Cas_n(\partial_xF).
\]
\end{lemma}

\begin{proof}
By Corollary~\ref{cor:cindep}, $D_{k,n}(x,c)=\det\bigl(L_x(X^{i+j}P_k)\bigr)_{i,j<n}$. Now apply \eqref{eq:christoffel} to the weights $s=(x,0,0,\dots)$, $\lambda=(1,1,\dots)$, for which $p_m=P_m$ (Lemma~\ref{lem:cheb}(iii)) and $h_m=1$, and to $f=P_k$.
\end{proof}

\subsection{The bilinear identity and the proofs of Theorems \texorpdfstring{\ref{thm:B} and~\ref{thm:C}}{1.1 and 1.2}}\label{ss:proofs}
For $\kk=\Q(x)$ the following proposition is the identity \eqref{eq:C2}, that is, \cite[Proposition~3.7]{Yu}.

\begin{proposition}\label{prop:bilinear}
Let $\kk$ be a field, $x\in\kk$, $k\ge1$, and let $w(m)$ be defined as in \S$\ref{ss:Lambda}$. Then for every $\ell\in\Z$,
\[
w(\ell+k+1)\,w(\ell-k-1)-x^2\,w(\ell+k)\,w(\ell-k)=(1-x^2)\,w(\ell)^2.
\]
\end{proposition}

\begin{proof}
Recall that $\det\Lambda\ne0$ (Lemma~\ref{lem:Lambda}(iii)). If $x=0$, then ($\Lambda1$) and ($\Lambda2$) read $F_{j+k+1}=\Lambda F_j=F_{j-k-1}$, hence $w(\ell+k+1)=\det(\Lambda)\,w(\ell)=w(\ell-k-1)$. Moreover $\Lambda^2$ is the identity by Lemma~\ref{lem:Lambda}(ii), so $(\det\Lambda)^2=1$, and the identity holds. (We use (ii) rather than the value of $\det\Lambda$ in (iii) so that the argument also applies in Corollary~\ref{thm:module}.) Let $x\ne0$. By \eqref{eq:shifts}, and since $\partial_x$ commutes with translations of the index and with the termwise application of $\Lambda$,
\[
(\partial_x\partial_{x^{-1}}F)_j=-x^{-1}\Lambda\,(\partial_xF)_{j+k+1}=-x^{-1}\Lambda^2F_{j+1}.
\]
Therefore, for $m=\ell-1$,
\begin{align*}
\Cas_m(\partial_0F)&=w(\ell), &\Cas_m(\partial_x\partial_{x^{-1}}F)&=(-x)^{-k}(\det\Lambda)^2\,w(\ell),\\
\Cas_m(\partial_xF)&=\det(\Lambda)\,w(\ell-k-1), &\Cas_m(\partial_0\partial_{x^{-1}}F)&=(-x)^{-k}\det(\Lambda)\,w(\ell+k+1),\\
\Cas_m(\partial_{x^{-1}}F)&=(-x)^{-k}\det(\Lambda)\,w(\ell+k), &\Cas_m(\partial_0\partial_xF)&=\det(\Lambda)\,w(\ell-k).
\end{align*}
Lemma~\ref{lem:HM} with $(a_1,a_2,a_3)=(0,x,x^{-1})$ gives
\[
(-x)^{-k}(\det\Lambda)^2\Bigl[(x-x^{-1})\,w(\ell)^2+x^{-1}\,w(\ell-k-1)\,w(\ell+k+1)-x\,w(\ell+k)\,w(\ell-k)\Bigr]=0.
\]
Dividing by $(-x)^{-k}(\det\Lambda)^2$ and multiplying by $x$ gives the assertion.
\end{proof}

\begin{proof}[Proof of Theorem~\ref{thm:C}]
Let $\kk=K$. Since $P_k=U_k-xU_{k-1}$ is monic in $X$ with coefficients in $\Z[x]$, division with remainder by $P_k$ can be carried out in $\Z[x][X]$, so the residue class of each $U_{j-1}\in\Z[X]$ $(j\in\Z)$ is represented by a polynomial in $\Z[x][X]$ of degree less than $k$. As $U_0,\dots,U_{k-1}$ are monic of degrees $0,\dots,k-1$ with integer coefficients, its coordinates with respect to $F_1,\dots,F_k$ lie in $\Z[x]$. Hence $w(m)\in\Z[x]$ for all $m\in\Z$. By Lemma~\ref{lem:reduction}, \eqref{eq:GLambda} and Lemma~\ref{lem:Lambda}(iii),
\[
D_{k,n}=(-1)^{nk}\detR(\Lambda F_{n-k},\dots,\Lambda F_{n-1})=(-1)^{nk}\det(\Lambda)\,w(n-k)=(-1)^{nk}\eps_k\,w(n-k)
\]
for all $n\ge0$,
which is \eqref{eq:C1}. The identity \eqref{eq:C2} is Proposition~\ref{prop:bilinear}.
\end{proof}

\begin{proof}[Proof of Theorem~\ref{thm:B}]
For $k=0$ we have $D_{k,n}=1$ (Corollary~\ref{cor:cindep}), and \eqref{eq:main} reads $1-x^2=1-x^2$. Let $k\ge1$, $n\ge2k+2$ and $\ell:=n-2k-1\ge1$. All of $n$, $n-1$, $n-k-1$, $n-2k-1$, $n-2k-2$ are $\ge0$, so \eqref{eq:C1} applies to each of them:
\begin{align*}
D_{k,n}&=(-1)^{nk}\eps_k\,w(\ell+k+1), & D_{k,n-2k-2}&=(-1)^{(n-2k-2)k}\eps_k\,w(\ell-k-1),\\
D_{k,n-1}&=(-1)^{(n-1)k}\eps_k\,w(\ell+k), & D_{k,n-2k-1}&=(-1)^{(n-2k-1)k}\eps_k\,w(\ell-k),\\
D_{k,n-k-1}&=(-1)^{(n-k-1)k}\eps_k\,w(\ell).
\end{align*}
The sign of $D_{k,n}D_{k,n-2k-2}$ and the sign of $D_{k,n-1}D_{k,n-2k-1}$ are both $(-1)^{2nk-2k^2-2k}\eps_k^2=1$, and $(D_{k,n-k-1})^2=w(\ell)^2$. Hence, by \eqref{eq:C2},
\begin{align*}
D_{k,n}D_{k,n-2k-2}-x^2D_{k,n-1}D_{k,n-2k-1}&=w(\ell+k+1)\,w(\ell-k-1)-x^2\,w(\ell+k)\,w(\ell-k)\\
&=(1-x^2)\,w(\ell)^2=(1-x^2)(D_{k,n-k-1})^2.
\end{align*}
This is an identity in $K$ between elements of $\Z[x]$, hence an identity in $\Z[x]$, and $D_{k,n}(x,c)$ does not depend on $c$.
\end{proof}

\begin{remark}[Comparison with the proof in \cite{Yu}]\label{rem:compare}\label{rem:firstproof}
Both proofs of \eqref{eq:C2} rest on a three-term Pl\"ucker relation and on ($\Lambda1$), ($\Lambda2$). In \cite[\S3.5]{Yu} the Pl\"ucker relation is applied to a basis of the two-dimensional kernel of the window $F_{\ell-1},\dots,F_{\ell+k}$. The resulting factor, a $2\times2$ determinant formed from the values at $x$ and $1/x$ of the two polynomials whose coefficient vectors form this basis, is evaluated by a divisibility argument in a model of $R$ by Laurent polynomials, and the windows of rank less than $k$, which do occur \cite[Remark~4.2]{Yu}, are treated separately. In the proof above the corresponding factor is $[\mathbf e(x),\mathbf e(x^{-1})]=(x^{-1}-x)\Cas(\partial_x\partial_{x^{-1}}F)$. Since $\partial_x\partial_{x^{-1}}F$ is $-x^{-1}\Lambda^2$ applied to the translate of $F$ by one step, this factor is again a Casoratian of $F$. Also, the divisibility argument of Step~4 of the proof of \cite[Proposition~3.7]{Yu} gives a second proof of the identity \eqref{eq:twovar} of \S\ref{ss:spectral} (with $\cW_m$ and $M_i(m)$ as in \S\ref{ss:spectral}). By Cramer's rule, $\sum_{i=0}^k(-1)^iM_i(m)F_{m+i}=0$ for all $m$. By \cite[Lemma~3.5(v)]{Yu}, for $m=\ell-1$ this says that $Q(t)$ divides $\Psi_\ell(t)-\Psi_\ell(t^{-1})$ in $K[t,t^{-1}]$, where $\Psi_\ell(t):=t^{\ell-1}\cW_{\ell-1}(t)$. Hence $Q$ divides $\Psi_\ell(t)\Psi_{\ell+1}(t^{-1})-\Psi_\ell(t^{-1})\Psi_{\ell+1}(t)=t^{-1}\cW_{\ell-1}(t)\cW_\ell(t^{-1})-t\,\cW_{\ell-1}(t^{-1})\cW_\ell(t)$. This Laurent polynomial has exponents in $[-k-1,k+1]$, so it is a constant multiple of $t^{-k-1}Q(t)$. The constant is its coefficient of $t^{k+1}$, which is $(-1)^{k+1}M_0(\ell-1)M_k(\ell)=(-1)^{k+1}w(\ell)^2$ by \eqref{eq:Wzero}. Multiplying by $t^{k+1}$ gives \eqref{eq:twovar} with $t$ in place of $y$.
\end{remark}

\section{Generalizations of Theorem \texorpdfstring{\ref{thm:B}}{1.1}}\label{sec:general}

In this section we state the principle behind the proof of Proposition~\ref{prop:bilinear} for $x\ne0$ (\S\ref{ss:reduction}), apply it to arbitrary solutions of the three-term recurrence and to general moment functionals (\S\ref{ss:general}), and derive an identity with a free parameter that holds for every monic polynomial in place of $P_k$ (\S\ref{ss:spectral}).

\subsection{The reduction principle}\label{ss:reduction}
We use the notation of \S\ref{ss:HM}: $\kk$ is a field and $E$ a $\kk$-vector space of dimension $r\ge1$ with a fixed basis.

\begin{proposition}[Reduction principle]\label{thm:reduction}
Let $F$ be a sequence in $E$ and $w(m):=\Cas_m(F)$ for $m\in\Z$. Let $a_1,a_2,a_3\in\kk$, and suppose that for each $i\in\{1,2,3\}$ there are an integer $\pi_i$ and an invertible $\kk$-linear map $\Lambda_i\colon E\to E$ such that
\begin{equation}\label{eq:contig}
F_{j+1}-a_iF_j=\Lambda_i\,F_{j+\pi_i}\qquad\text{for all }j\in\Z.
\end{equation}
Then for all $m\in\Z$
\begin{equation}\label{eq:GRgeneral}
\begin{aligned}
&(a_2-a_3)\,w(m+\pi_1)\,w(m+\pi_2+\pi_3)+(a_3-a_1)\,w(m+\pi_2)\,w(m+\pi_1+\pi_3)\\
&\qquad+(a_1-a_2)\,w(m+\pi_3)\,w(m+\pi_1+\pi_2)=0.
\end{aligned}
\end{equation}
\end{proposition}

\begin{proof}
By \eqref{eq:contig}, $(\partial_{a_i}F)_j=\Lambda_iF_{j+\pi_i}$, hence $\Cas_m(\partial_{a_i}F)=\det(\Lambda_i)\,w(m+\pi_i)$. For $i\ne j$ and $l\in\Z$,
\[
(\partial_{a_i}\partial_{a_j}F)_l=(\partial_{a_j}F)_{l+1}-a_i(\partial_{a_j}F)_l=\Lambda_j\bigl(F_{l+1+\pi_j}-a_iF_{l+\pi_j}\bigr)=\Lambda_j\Lambda_i\,F_{l+\pi_i+\pi_j},
\]
hence $\Cas_m(\partial_{a_i}\partial_{a_j}F)=\det(\Lambda_i)\det(\Lambda_j)\,w(m+\pi_i+\pi_j)$. Inserting this into \eqref{eq:HM} shows that the left-hand side of \eqref{eq:GRgeneral}, multiplied by $\det(\Lambda_1)\det(\Lambda_2)\det(\Lambda_3)\ne0$, is zero.
\end{proof}

For $a_1=0$ the hypothesis \eqref{eq:contig} holds trivially with $\pi_1=1$ and $\Lambda_1=\mathrm{id}$. This gives the following form of the principle.

\begin{corollary}\label{cor:twoparam}
Let $F$ be a sequence in $E$ and $w(m)=\Cas_m(F)$. Let $a,b\in\kk$, $\pi,\pi'\in\Z$ and let $\Lambda_a,\Lambda_b$ be invertible $\kk$-linear maps $E\to E$ such that $F_{j+1}-aF_j=\Lambda_aF_{j+\pi}$ and $F_{j+1}-bF_j=\Lambda_bF_{j+\pi'}$ for all $j\in\Z$. Then for all $m\in\Z$
\begin{equation}\label{eq:GRtwo}
(a-b)\,w(m+1)\,w(m+\pi+\pi')+b\,w(m+\pi)\,w(m+\pi'+1)-a\,w(m+\pi')\,w(m+\pi+1)=0.
\end{equation}
\end{corollary}

For $x\ne0$, $a=x$, $\pi=-k$, $\Lambda_a=\Lambda$ and $b=x^{-1}$, $\pi'=k+1$, $\Lambda_b=-x^{-1}\Lambda$ the hypotheses are \eqref{eq:shifts}, and \eqref{eq:GRtwo}, multiplied by $x$, is Proposition~\ref{prop:bilinear}.

\begin{remark}\label{rem:GR}
The equation \eqref{eq:GRgeneral} is a three-term bilinear recurrence of Gale--Robinson type, since the three products have the same index sum $2m+\pi_1+\pi_2+\pi_3$. Proposition~\ref{thm:reduction} gives a sufficient condition for a Casorati determinant solution to be compatible with a reduction of the kind considered in \cite{HKW}, namely that each operation $\partial_{a_i}$ is a translation followed by a fixed invertible linear map. For solutions built from exponentials see Remark~\ref{rem:soliton}.
\end{remark}

\subsection{Arbitrary solutions of the three-term recurrence; general moment functionals}\label{ss:general}

\begin{corollary}\label{thm:module}
Let $k\ge1$, let $\kk$ be a field, $x\in\kk$, and $R=\kk[X]/(U_k-xU_{k-1})$. Let $\Mod$ be an $R$-module of finite dimension $r\ge1$ over $\kk$, with a fixed $\kk$-basis, and let $f=(f_j)_{j\in\Z}$ be a sequence in $\Mod$ with $f_{j+1}=Xf_j-f_{j-1}$ for all $j\in\Z$. Then the Casoratians $\Cas_m(f)=\detM(f_m,f_{m+1},\dots,f_{m+r-1})$ satisfy
\begin{equation}\label{eq:Wmodule}
\Cas_{\ell+k+1}(f)\,\Cas_{\ell-k-1}(f)-x^2\,\Cas_{\ell+k}(f)\,\Cas_{\ell-k}(f)=(1-x^2)\,\Cas_\ell(f)^2\qquad\text{for all }\ell\in\Z.
\end{equation}
\end{corollary}

\begin{proof}
For all $j\in\Z$ we have $f_j=U_{j-1}(X)f_1-U_{j-2}(X)f_0$: both sides satisfy the recurrence in $j$, and they agree for $j=0$ and $j=1$ (recall $U_{-1}=0$, $U_{-2}=-1$, $U_0=1$). Let $\Lambda_{\Mod}\colon\Mod\to\Mod$ be the action of the class of $-P_{k-1}$, where $P_m=U_m-xU_{m-1}$. By Lemma~\ref{lem:Lambda}(i), in $R$ the class of $-P_{k-1}U_{j-1}$ equals the class of $U_{j+k}-xU_{j+k-1}$ and also the class of $U_{j-k-2}-xU_{j-k-1}$, for every $j\in\Z$. Applying this to $f_1$ and, with $j-1$ in place of $j$, to $f_0$, we get
\[
\Lambda_{\Mod}f_j=f_{j+k+1}-xf_{j+k}=f_{j-k-1}-xf_{j-k}\qquad(j\in\Z).
\]
By Lemma~\ref{lem:Lambda}(iii) the class of $P_{k-1}$ is a unit of $R$, so $\Lambda_{\Mod}$ is invertible. By Lemma~\ref{lem:Lambda}(ii), applied to $F_1=1$, the class of $P_{k-1}^2$ is that of $1+x^2-xX$. Hence $\Lambda_{\Mod}^2$ is the action of $1+x^2-xX$, and $\Lambda_{\Mod}^2=\mathrm{id}$ if $x=0$. Now the proof of Proposition~\ref{prop:bilinear} in \S\ref{ss:proofs} applies verbatim, with $E=\Mod$, with $f$ and $\Lambda_{\Mod}$ in place of $F$ and $\Lambda$, and with Casoratians of size $r$.
\end{proof}

For $\Mod=R$, $f_0=0$ and $f_1=1$ this is Proposition~\ref{prop:bilinear}. The size $r$ of the determinants need not be $k$, and $f_0,f_1\in\Mod$ are arbitrary. For example, let $\rho_1,\dots,\rho_k$ be the roots of $P_k$ in an algebraic closure $\Omega$ of $\kk$, choose $\omega_j\in\Omega$ with $\omega_j+\omega_j^{-1}=\rho_j$, let $j_1,\dots,j_r\in\{1,\dots,k\}$, and let $\Mod=\Omega^r$, where $X$ acts by the diagonal matrix with entries $\rho_{j_1},\dots,\rho_{j_r}$. This is a module over $\Omega[X]/(P_k)$, and for arbitrary $\alpha_i,\beta_i\in\Omega$ the vectors $f_m:=\bigl(\alpha_i\omega_{j_i}^{m}+\beta_i\omega_{j_i}^{-m}\bigr)_{1\le i\le r}$ clearly satisfy the recurrence. So the identity \eqref{eq:C2} holds for the determinants
\begin{equation}\label{eq:solitondet}
\det\bigl(\alpha_i\,\omega_{j_i}^{\,m+l}+\beta_i\,\omega_{j_i}^{-m-l}\bigr)_{0\le l\le r-1,\ 1\le i\le r}
\end{equation}
formed with any of the pairs $(\omega_j,\omega_j^{-1})$ and arbitrary coefficients. When the $2k$ elements $\omega_j^{\pm1}$ are pairwise distinct, $w(m)$ is, up to a constant factor, the case of all $k$ pairs and $\alpha_i=-\beta_i$ (see Remark~\ref{rem:symplectic}).

Corollary~\ref{thm:module} has the following consequence for Hankel determinants. We use the notation of \S\ref{ss:matrices} for weights $(s,\lambda)$ with values in a commutative ring $A$.

\begin{theorem}\label{thm:hankelgen}
Let $k\ge1$ and $N\ge0$ be integers, let $x\in A$, and let $(s,\lambda)$ be weights with $s_n=0$ and $\lambda_n=1$ for all $n>N$. Let $L=L^{(s,\lambda)}$ and put
\[
b_n:=L\bigl(X^n\,(U_k-xU_{k-1})\bigr)\quad(n\ge0),\qquad \cH_n:=\det\bigl(b_{i+j}\bigr)_{0\le i,j\le n-1},\quad \cH_0:=1.
\]
Then for all $n\ge 2k+2+\max\bigl(0,\,N-\lfloor k/2\rfloor\bigr)$
\begin{equation}\label{eq:Hgen}
\cH_n\,\cH_{n-2k-2}-x^2\,\cH_{n-1}\,\cH_{n-2k-1}=(1-x^2)\,\cH_{n-k-1}^{\,2}\qquad\text{in }A.
\end{equation}
\end{theorem}

For $s=(x,0,0,\dots)$ and $\lambda=(1,1,\dots)$ we have $N=0$ and $b_n=a_{n,k}(x,0)$ by Lemma~\ref{lem:cheb}(iii), so that Theorem~\ref{thm:hankelgen} contains Theorem~\ref{thm:B}. In general the functional $L$ and the polynomial $U_k-xU_{k-1}$ are unrelated. The recurrence comes from the factor $U_k-xU_{k-1}$ alone, provided that above height $N$ the weights are those of the polynomials $U_m$.

\begin{proof}
The elements $b_n$ and $\cH_n$ are the values, at the given elements of $A$, of polynomials with integer coefficients in $x,s_0,\dots,s_N,\lambda_1,\dots,\lambda_N$. So it suffices to prove \eqref{eq:Hgen} when these $2N+2$ elements are indeterminates and $A$ is the polynomial ring over $\Z$ generated by them; let $\kk$ be its field of fractions, and let $P_k=U_k-xU_{k-1}$, $R=\kk[X]/(P_k)$. By \eqref{eq:christoffel},
\begin{equation}\label{eq:Hchr}
\cH_n=(-1)^{nk}\,h_0h_1\cdots h_{n-1}\,\detR(p_n,\dots,p_{n+k-1})\qquad(n\ge0).
\end{equation}
Since $s_m=0$ and $\lambda_m=1$ for $m>N$, we have $p_{m+1}=Xp_m-p_{m-1}$ for all $m>N$. Let $f=(f_j)_{j\in\Z}$ be the sequence in $R$ with $f_{j+1}=Xf_j-f_{j-1}$ for all $j\in\Z$ whose terms $f_N$ and $f_{N+1}$ are the residue classes of $p_N$ and $p_{N+1}$; then $f_j$ is the residue class of $p_j$ for all $j\ge N$. By Corollary~\ref{thm:module}, the Casoratians $\Cas_m(f)=\detR(f_m,\dots,f_{m+k-1})$ satisfy \eqref{eq:Wmodule}. Put $\Theta:=h_N$ and $\theta:=\bigl(\prod_{i=0}^{N-1}h_i\bigr)\Theta^{-N}$. Since $h_i=\Theta$ for $i\ge N$, we have $h_0\cdots h_{m-1}=\theta\,\Theta^m$ for $m\ge N$. We claim that
\begin{equation}\label{eq:HW}
\cH_m=(-1)^{mk}\,\theta\,\Theta^m\,\Cas_m(f)\qquad\text{for all }m\ge\max\bigl(0,\,N-\lfloor k/2\rfloor\bigr).
\end{equation}
For $m\ge N$ this is \eqref{eq:Hchr}. Suppose the claim is proved. Let $n\ge2k+2+\max(0,N-\lfloor k/2\rfloor)$ and $\ell:=n-k-1$. Then \eqref{eq:HW} applies to the five indices $n$, $n-1$, $n-k-1$, $n-2k-1$, $n-2k-2$, and each of the products $\cH_n\cH_{n-2k-2}$, $\cH_{n-1}\cH_{n-2k-1}$, $\cH_{n-k-1}^{\,2}$ is $\theta^2\Theta^{2n-2k-2}$ times the corresponding product of Casoratians of $f$ (the signs $(-1)^{(2n-2k-2)k}$ are $1$). So \eqref{eq:Hgen} is \eqref{eq:Wmodule} multiplied by $\theta^2\Theta^{2n-2k-2}$.

It remains to prove \eqref{eq:HW} for $N-\lfloor k/2\rfloor\le m<N$, $m\ge0$. For $0\le i\le N$ put $\lambda^{(i)}:=\lambda_N\lambda_{N-1}\cdots\lambda_{N-i+1}$ (so $\lambda^{(0)}=1$). We write $p_j$ also for the residue class of $p_j$ in $R$, and show
\begin{equation}\label{eq:backwards}
f_{N-i}\in\lambda^{(i)}\,p_{N-i}+\operatorname{span}\bigl(p_{N-i+1},\dots,p_{N+i-1}\bigr)\qquad(0\le i\le N).
\end{equation}
For $i=0$ this is $f_N=p_N$; for $i=1$ it is $f_{N-1}=Xf_N-f_{N+1}=Xp_N-p_{N+1}=\lambda_Np_{N-1}+s_Np_N$; and if it holds for $i-1$ and $i$ (with $1\le i<N$), it holds for $i+1$, because $f_{N-i-1}=Xf_{N-i}-f_{N-i+1}$ and $Xp_j=p_{j+1}+s_jp_j+\lambda_jp_{j-1}$ for $j\ge1$. Now let $m=N-i_0$ with $1\le i_0\le\min(N,\lfloor k/2\rfloor)$. The family $f_m,\dots,f_{m+k-1}$ consists of $f_{N-i_0},\dots,f_{N-1}$ followed by $p_N,\dots,p_{m+k-1}$, and $m+k-1\ge N+i_0-1$ because $k\ge2i_0$. Replace successively, for $i=1,2,\dots,i_0$, the member $f_{N-i}$ by $\lambda^{(i)}p_{N-i}$. By \eqref{eq:backwards}, each replacement changes $f_{N-i}$ by a linear combination of $p_{N-i+1},\dots,p_{N+i-1}$, which at that stage are other members of the family (they belong to it because $N+i-1\le m+k-1$) up to the factors $\lambda^{(i')}$, $i'<i$; these factors are non-zero elements of $\kk$. So the determinant does not change, and
\[
\Cas_m(f)=\lambda^{(1)}\lambda^{(2)}\cdots\lambda^{(i_0)}\,\detR(p_m,\dots,p_{m+k-1}).
\]
Since $\lambda^{(i)}h_{N-i}=h_N=\Theta$, we have $\lambda^{(1)}\cdots\lambda^{(i_0)}=\Theta^{N-m}/(h_mh_{m+1}\cdots h_{N-1})$, and \eqref{eq:Hchr} gives
\begin{align*}
\cH_m&=(-1)^{mk}\,h_0\cdots h_{m-1}\,\detR(p_m,\dots,p_{m+k-1})\\
&=(-1)^{mk}\,h_0\cdots h_{N-1}\,\Theta^{-(N-m)}\,\Cas_m(f)=(-1)^{mk}\,\theta\,\Theta^m\,\Cas_m(f).\qedhere
\end{align*}
\end{proof}

\begin{remark}\label{rem:range}
(a) \emph{Constant diagonal.} If $s_n=c$ and $\lambda_n=1$ for all $n>N$, with some $c\in A$, the statement of Theorem~\ref{thm:hankelgen} holds for $b_n=L\bigl(X^n\,P_k(X-c)\bigr)$, where $P_k=U_k-xU_{k-1}$. Indeed, let $L^\circ$ be the functional of the weights $(s_j-c)_{j\ge0}$, $\lambda$, to which the theorem applies; by Lemma~\ref{lem:jacobi}(iv) the Hankel determinants of the sequences $\bigl(L(X^nP_k(X-c))\bigr)_{n\ge0}$ and $\bigl(L^\circ(X^nP_k)\bigr)_{n\ge0}$ coincide.

(b) \emph{The range of $n$.} If moreover $N\ge1$ and $s_N=0$, then the span in \eqref{eq:backwards} can be replaced by $\operatorname{span}(p_{N-i+1},\dots,p_{N+i-2})$ for $i\ge1$, and the same argument gives \eqref{eq:Hgen} for all $n\ge2k+2+\max(0,N-\lceil k/2\rceil)$. In $56$ random examples in exact arithmetic ($1\le k\le7$, $1\le N\le4$, $\lambda_N\ne1$, so that $N$ is minimal) the smallest $n$ from which \eqref{eq:Hgen} holds was exactly the bound with $\lfloor k/2\rfloor$ if $s_N\ne0$ and the bound with $\lceil k/2\rceil$ if $s_N=0$; this is a numerical observation. We do not know whether these bounds are best possible for generic weights.
\end{remark}

\subsection{An identity with a spectral parameter; symplectic characters}\label{ss:spectral}
For an arbitrary monic polynomial $P$ in place of $P_k$ there is no map $\Lambda$, but Lemma~\ref{lem:HM} still gives an identity with a free parameter.

Let $\kk$ be a field, let $P\in\kk[X]$ be monic of degree $k\ge1$, and for $j\in\Z$ let $F_j$ be the residue class of $U_{j-1}$ in $\kk[X]/(P)$; the classes $F_1,\dots,F_k$ form a monic basis. For $m\in\Z$ and $0\le i\le k$ put
\[
w_P(m):=\detx{P}(F_m,\dots,F_{m+k-1}),\qquad M^P_i(m):=\detx{P}\bigl(F_m,\dots,\widehat{F_{m+i}},\dots,F_{m+k}\bigr),
\]
and define polynomials in an indeterminate $y$ by
\[
\cW^P_m(y):=\sum_{i=0}^k(-y)^i\,M^P_i(m),\qquad\cW^{P*}_m(y):=\sum_{i=0}^k(-1)^iM^P_i(m)\,y^{k-i}=y^k\,\cW^P_m(1/y)\qquad\text{in }\kk[y].
\]
By Lemma~\ref{lem:expansion}, for every element $y_0$ of an extension field of $\kk$,
\begin{equation}\label{eq:Wdet}
\begin{aligned}
\cW^P_m(y_0)&=\detx{P}\bigl(F_{m+1+i}-y_0F_{m+i}\bigr)_{0\le i\le k-1}=\Cas_m(\partial_{y_0}F),\\
\cW^{P*}_m(y_0)&=\detx{P}\bigl(y_0F_{m+1+i}-F_{m+i}\bigr)_{0\le i\le k-1},
\end{aligned}
\end{equation}
where $\detx{P}(\mathbf y_i)_{0\le i\le k-1}$ stands for $\detx{P}(\mathbf y_0,\dots,\mathbf y_{k-1})$. In particular
\begin{equation}\label{eq:Wzero}
\cW^P_m(0)=M^P_0(m)=w_P(m+1),\qquad\cW^{P*}_m(0)=(-1)^kM^P_k(m)=(-1)^kw_P(m).
\end{equation}
Finally let
\[
Q_P(y):=y^k\,(y^2-1)\,P\bigl(y+y^{-1}\bigr)\in\kk[y],
\]
a polynomial of degree $2k+2$. For $P=P_k=U_k-xU_{k-1}$ we omit the index $P$ and write $w(m)$, $M_i(m)$, $\cW_m(y)$, $\cW^*_m(y)$; in this case $Q_P$ is the polynomial
\[
Q(y)=y^{2k+2}-xy^{2k+1}+xy-1
\]
of \cite[\S3.4]{Yu}, by \eqref{eq:cheb}.

\begin{proposition}[Identity with a spectral parameter]\label{thm:spectral}
Let $P\in\kk[X]$ be monic of degree $k\ge1$. Then for every $\ell\in\Z$ we have, in $\kk[y]$,
\begin{equation}\label{eq:spectral}
y^2\,\cW^P_\ell(y)\,\cW^{P*}_{\ell-1}(y)-\cW^{P*}_\ell(y)\,\cW^P_{\ell-1}(y)=(-1)^k\,Q_P(y)\,w_P(\ell)^2.
\end{equation}
\end{proposition}

\begin{proof}
We work over the field $\kk(y)$ and apply Lemma~\ref{lem:HM} to the sequence $F$ in $E=\kk(y)[X]/(P)$, with $(a_1,a_2,a_3)=(0,y,y^{-1})$ and $m=\ell-1$. By \eqref{eq:Wdet} we have $\Cas_{\ell-1}(\partial_0F)=w_P(\ell)$, $\Cas_{\ell-1}(\partial_yF)=\cW^P_{\ell-1}(y)$, $\Cas_{\ell-1}(\partial_0\partial_yF)=\Cas_\ell(\partial_yF)=\cW^P_\ell(y)$, and the same with $y^{-1}$ in place of $y$. Moreover
\[
(\partial_y\partial_{y^{-1}}F)_j=F_{j+2}-(y+y^{-1})F_{j+1}+F_j=\bigl(X-(y+y^{-1})\bigr)F_{j+1}
\]
by the recurrence of the $U_m$. The characteristic polynomial of the multiplication by $X$ on $E$ is $P$, so the multiplication by $X-\xi$ has determinant $(-1)^kP(\xi)$ for every $\xi\in\kk(y)$. Hence
\[
\Cas_{\ell-1}(\partial_y\partial_{y^{-1}}F)=(-1)^kP(y+y^{-1})\,w_P(\ell),
\]
and \eqref{eq:HM} becomes
\[
(y-y^{-1})(-1)^kP(y+y^{-1})\,w_P(\ell)^2+y^{-1}\,\cW^P_{\ell-1}(y)\,\cW^P_\ell(y^{-1})-y\,\cW^P_{\ell-1}(y^{-1})\,\cW^P_\ell(y)=0.
\]
Multiplying by $y^{k+1}$ gives \eqref{eq:spectral}.
\end{proof}

\begin{corollary}[The case $P=P_k$]\label{cor:twovar}
Let $k\ge1$, $\kk=K=\Q(x)$ and $P=P_k$. Then $M_i(m)\in\Z[x]$, and for every $\ell\in\Z$
\begin{equation}\label{eq:twovar}
y^2\,\cW_\ell(y)\,\cW^*_{\ell-1}(y)-\cW^*_\ell(y)\,\cW_{\ell-1}(y)=(-1)^k\,Q(y)\,w(\ell)^2\qquad\text{in }\Z[x,y].
\end{equation}
Moreover, for all $m\in\Z$,
\begin{equation}\label{eq:Wspecial}
\begin{aligned}
\cW_m(0)&=w(m+1), &\cW^*_m(0)&=(-1)^kw(m),\\
\cW_m(x)&=\eps_k\,w(m-k), &\cW^*_m(x)&=(-1)^k\eps_k\,w(m+k+1),
\end{aligned}
\end{equation}
and for $y=x$ the identity \eqref{eq:twovar} is the identity of Proposition~\ref{prop:bilinear}.
\end{corollary}

\begin{proof}
The elements $M_i(m)$ lie in $\Z[x]$ for the reason given in the proof of Theorem~\ref{thm:C}. The first two formulas of \eqref{eq:Wspecial} are \eqref{eq:Wzero}. By \eqref{eq:shifts}, $F_{m+1+i}-xF_{m+i}=\Lambda F_{m+i-k}$, so \eqref{eq:Wdet} gives $\cW_m(x)=\det(\Lambda)\,w(m-k)$. By ($\Lambda2$), $F_{m+i}-xF_{m+1+i}=\Lambda F_{m+i+k+1}$, so \eqref{eq:Wdet} gives $\cW^*_m(x)=(-1)^k\detR\bigl(F_{m+i}-xF_{m+1+i}\bigr)_{0\le i\le k-1}=(-1)^k\det(\Lambda)\,w(m+k+1)$. Finally $\det\Lambda=\eps_k$. Inserting \eqref{eq:Wspecial} into \eqref{eq:twovar} at $y=x$, where $Q(x)=x^2-1$, gives $(-1)^kx^2w(\ell-k)w(\ell+k)-(-1)^kw(\ell+k+1)w(\ell-k-1)=(-1)^k(x^2-1)w(\ell)^2$.
\end{proof}

By \eqref{eq:christoffel}, the polynomials $\cW_m(y)$ with $m\ge0$ are, up to sign, Hankel determinants:

\begin{corollary}[Hankel form]\label{cor:hankel2}
Let $k\ge1$, and for $n\ge0$ let
\[
\cH_n(x;y):=\det\bigl(L_y(X^{i+j}P_k)\bigr)_{0\le i,j\le n-1}\in\Z[x,y],\qquad\cH^*_n(x;y):=y^k\,\cH_n(x;1/y).
\]
\begin{enumerate}
\item[(i)] For every $c$,
\[
\cH_n(x;y)=\det\Bigl(a_{i+j,k}(y,c)+(y-x)\sum_{l=0}^{k-1}y^l\,a_{i+j,k-1-l}(y,c)\Bigr)_{0\le i,j\le n-1},
\]
and $\cH_n(x;y)=(-1)^{nk}\,\cW_n(y)$; in particular $\cH_n(x;y)$ is a polynomial in $y$ of degree at most $k$, and $\cH^*_n(x;y)\in\Z[x,y]$.
\item[(ii)] For all $n\ge1$,
\begin{equation}\label{eq:hankel2}
y^2\,\cH_n(x;y)\,\cH^*_{n-1}(x;y)-\cH^*_n(x;y)\,\cH_{n-1}(x;y)=Q(y)\,D_{k,n+k}(x)^2,
\end{equation}
and for all $n\ge0$
\begin{equation}\label{eq:special}
\begin{aligned}
\cH_n(x;x)&=D_{k,n}, &\cH^*_n(x;x)&=D_{k,n+2k+1},\\
\cH_n(x;0)&=\eps_k\,D_{k,n+k+1}, &\cH^*_n(x;0)&=\eps_k\,D_{k,n+k}.
\end{aligned}
\end{equation}
In particular, for every $c$,
\begin{equation}\label{eq:companion}
\det\bigl(a_{i+j,k}(0,c)-x\,a_{i+j,k-1}(0,c)\bigr)_{0\le i,j\le n-1}=\eps_k\,D_{k,n+k+1}(x)\qquad(n\ge0).
\end{equation}
\item[(iii)] For all $n\ge2k+2$,
\[
\cH_n(x;y)\,\cH_{n-2k-2}(x;y)-x^2\,\cH_{n-1}(x;y)\,\cH_{n-2k-1}(x;y)=(1-x^2)\,\cH_{n-k-1}(x;y)^2\qquad\text{in }\Z[x,y].
\]
\end{enumerate}
For $y=x$ the identity \eqref{eq:hankel2} is \eqref{eq:main} with $n$ replaced by $n+2k+1$, and {\rm(iii)} is Theorem~$\ref{thm:B}$.
\end{corollary}

\begin{proof}
(i) Since $P_{j,y}=U_j-yU_{j-1}$, the sum $\sum_{l=0}^{k-1}y^lP_{k-1-l,y}$ telescopes to $U_{k-1}$, so $P_k=P_{k,y}+(y-x)\sum_{l=0}^{k-1}y^lP_{k-1-l,y}$. By \eqref{eq:shiftm}, over the ring $\Z[x,y,c]$, the $(i,j)$ entry of the matrix in (i) is $L_y\bigl((X+c)^{i+j}P_k\bigr)$, and its determinant is $\cH_n(x;y)$ by Lemma~\ref{lem:jacobi}(iv). By \eqref{eq:christoffel} over the field $\Q(x,y)$, for the weights $(y,0,0,\dots)$, $(1,1,\dots)$ and $f=P_k$, and by \eqref{eq:Wdet},
\begin{align*}
\cH_n(x;y)&=(-1)^{nk}\detR\bigl(P_{n,y},\dots,P_{n+k-1,y}\bigr)\\
&=(-1)^{nk}\detR\bigl(F_{n+1+i}-yF_{n+i}\bigr)_{0\le i\le k-1}=(-1)^{nk}\,\cW_n(y).
\end{align*}

(ii) By (i), \eqref{eq:hankel2} follows from \eqref{eq:twovar} with $\ell=n$ (the sign is $(-1)^{nk+(n-1)k}(-1)^k=1$) and from $w(n)^2=D_{k,n+k}^{\,2}$, which holds by \eqref{eq:C1}. The identities \eqref{eq:special} follow from \eqref{eq:Wspecial} and \eqref{eq:C1}:
\begin{align*}
\cH_n(x;x)&=(-1)^{nk}\eps_kw(n-k)=D_{k,n}, & \cH^*_n(x;x)&=(-1)^{nk+k}\eps_kw(n+k+1)=D_{k,n+2k+1},\\
\cH_n(x;0)&=(-1)^{nk}w(n+1)=\eps_kD_{k,n+k+1}, & \cH^*_n(x;0)&=(-1)^{nk+k}w(n)=\eps_kD_{k,n+k},
\end{align*}
where we used $(-1)^{k(k+1)}=1$. For $y=0$ only the term $l=0$ of the sum in (i) survives, which gives \eqref{eq:companion}.

(iii) This is Theorem~\ref{thm:hankelgen} for $A=\Z[x,y]$, the weights $s=(y,0,0,\dots)$, $\lambda=(1,1,\dots)$ and $N=0$.

Finally, by \eqref{eq:special}, \eqref{eq:hankel2} at $y=x$ reads $x^2D_{k,n}D_{k,n+2k}-D_{k,n+2k+1}D_{k,n-1}=(x^2-1)D_{k,n+k}^{\,2}$.
\end{proof}

By (i) for $y=0$, the quantity $w(m)$ of Theorem~\ref{thm:C} is itself a Hankel determinant: for $m\ge1$, $w(m)=(-1)^{(m-1)k}\,\cH_{m-1}(x;0)=(-1)^{(m-1)k}\det\bigl(a_{i+j,k}(0,0)-x\,a_{i+j,k-1}(0,0)\bigr)_{0\le i,j\le m-2}$, where the $a_{n,j}(0,0)$ are the ballot numbers.

\begin{remark}[Symplectic characters]\label{rem:symplectic}
Let $P=\prod_{j=1}^k(X-\omega_j-\omega_j^{-1})$ with $\omega_1,\dots,\omega_k$ in an extension field of $\kk$, such that the $2k$ elements $\omega_j^{\pm1}$ are pairwise distinct. Then the roots $\rho_j=\omega_j+\omega_j^{-1}$ of $P$ are pairwise distinct, and $\omega_j\ne\omega_j^{-1}$. The computation of \cite[Remark~4.4]{Yu} applies to $P$. If $C_m$ is the matrix whose rows are the coordinate vectors of $F_m,\dots,F_{m+k-1}$, so that $w_P(m)=\det C_m$, then evaluation at $\rho_1,\dots,\rho_k$ gives $\bigl(U_{m-1+i}(\rho_j)\bigr)_{i,j}=C_m\bigl(U_l(\rho_j)\bigr)_{l,j}$, and $\det\bigl(U_l(\rho_j)\bigr)_{0\le l\le k-1,\,1\le j\le k}=\prod_{i<j}(\rho_j-\rho_i)\ne0$. By \eqref{eq:cheb} with $t=\omega_j$, $(\omega_j-\omega_j^{-1})U_n(\rho_j)=\omega_j^{n+1}-\omega_j^{-n-1}$. Hence, by the Weyl character formula for $\Sp(2k)$ \cite[\S24.2, (24.18)]{FH},
\[
w_P(m)=\frac{\det\bigl(\omega_j^{m+i}-\omega_j^{-m-i}\bigr)_{0\le i\le k-1,\,1\le j\le k}}{\det\bigl(\omega_j^{1+i}-\omega_j^{-1-i}\bigr)_{0\le i\le k-1,\,1\le j\le k}}=\symp_{((m-1)^k)}(\omega_1,\dots,\omega_k)\qquad(m\ge1),
\]
the character of the irreducible representation of $\Sp(2k)$ indexed by the rectangular partition $((m-1)^k)$, evaluated at $\omega=(\omega_1,\dots,\omega_k)$. In the same way $M^P_i(m)=\symp_{(m^{k-i},(m-1)^i)}(\omega)$ for $m\ge1$, so that
\[
\cW^P_m(y)=\sum_{i=0}^k(-y)^i\,\symp_{(m^{k-i},(m-1)^i)}(\omega)\qquad(m\ge1).
\]
Thus, for $\ell\ge2$, \eqref{eq:spectral} is an identity between symplectic characters of rectangular and nearly rectangular shape at an arbitrary point $\omega$, with a free parameter $y$, and it can be viewed as a generating function of three-term Pl\"ucker relations. For Schur functions, Kleber \cite{Kle} proves Pl\"ucker relations for arbitrary shapes, which for rectangular shapes are the discrete Hirota relations, and he adds spectral parameters for quantum analogues of the characters (\cite[Theorem~4.4]{Kle}, numbering of arXiv:math/9907177v1). Fulmek and Kleber \cite{FK} give bijective proofs. Krattenthaler \cite{K98} gives identities of a different, linear kind (branching and tensor product rules) for classical group characters of nearly rectangular shape.

In the situation of Theorem~\ref{thm:C} the hypothesis is satisfied, $\omega_1^{\pm1},\dots,\omega_k^{\pm1}$ are the roots of $Q(t)/(t^2-1)=t^kP_k(t+t^{-1})$, and $w(m)=\symp_{((m-1)^k)}(\omega)$ for $m\ge1$, $x=\symp_{(1)}(\omega)$ \cite[Remark~4.4]{Yu}. So for $\ell\ge k+2$ the identity \eqref{eq:C2} is a bilinear relation between rectangular symplectic characters at the point $\omega$. As observed in \cite{Yu}, it does not hold at a general point. For example, let $k=2$, $x=3$ and $P=U_2-3U_1-\frac12U_0$, so that $\symp_{(1)}(\omega)=3$. A direct computation gives $w_P(1)=1$, $w_P(2)=-\frac12$, $w_P(4)=\frac{139}8$, $w_P(6)=-\frac{8277}{32}$ and $w_P(7)=-\frac{18605}{64}$. So for $\ell=4$ the left-hand side of \eqref{eq:C2}, with $w_P$ in place of $w$, is $-\frac{46549}{32}$, and the right-hand side is $-\frac{19321}8$.

This gives a partial answer to the question of \cite[Remark~4.4]{Yu} whether \eqref{eq:C2} is a specialization of a known identity for symplectic characters. The identity \eqref{eq:C2} is obtained from \eqref{eq:spectral}, which holds at every point, by putting $y=x$ and restricting $\omega$ to the points with $P=U_k-xU_{k-1}$, where \eqref{eq:Wspecial} holds. These points are described by the vanishing of fundamental characters. Since
\[
\prod_{j=1}^k\bigl(X-\omega_j-\omega_j^{-1}\bigr)=\sum_{i=0}^k(-1)^i\,\symp_{(1^i)}(\omega)\,U_{k-i}(X)
\]
(this follows from $\symp_{(1^i)}=\sigma_i-\sigma_{i-2}$, where $\sigma_i$ is the $i$-th elementary symmetric polynomial of the $2k$ elements $\omega_j^{\pm1}$, and from \eqref{eq:cheb}), we have $P=U_k-xU_{k-1}$ if and only if
\[
\symp_{(1)}(\omega)=x\qquad\text{and}\qquad\symp_{(1^i)}(\omega)=0\quad(2\le i\le k).
\]
Thus, for $\ell\ge k+2$, the identity \eqref{eq:C2} says that on the line $\{\symp_{(1^2)}=\dots=\symp_{(1^k)}=0\}$ the rectangular characters $\symp_{((\ell-1)^k)}$ satisfy a bilinear recurrence in $\ell$ whose coefficients are $\symp_{(1)}^2$ and $1-\symp_{(1)}^2$. We do not know whether \eqref{eq:spectral} can be obtained by specialization from known determinantal identities for symplectic characters.
\end{remark}
\section{Shifted Hankel determinants}\label{sec:shifted}

As in \cite{C23} we write, for $m\in\Z$ and $k,n\ge0$,
\[
D_{m,k,n}(x,c):=\det\bigl(a_{i+j+m,k}(x,c)\bigr)_{0\le i,j\le n-1}\in\Z[x,c],\qquad a_{r,k}(x,c):=0\ \text{ for } r<0,
\]
so that $D_{0,k,n}=D_{k,n}$. Throughout this section $k$ is fixed, $K=\Q(x)$, and
\[
d_n:=D_{k,n}(x)\in\Z[x].
\]
For $k\ge1$ we use the notation of \S\ref{ss:Lambda} and \S\ref{ss:spectral} for $\kk=K$ and $P=P_k$; recall \eqref{eq:C1},
\begin{equation}\label{eq:dw}
d_n=(-1)^{nk}\eps_k\,w(n-k)\qquad(n\ge0).
\end{equation}
Moreover $d_0=1$, $d_n=0$ for $1\le n\le k$ (Corollary~\ref{cor:cindep}) and $d_{k+1}=\eps_k$: for $n=k+1$ the matrix defining $d_n$ has the entries $1$ on the antidiagonal and $0$ above it (Lemma~\ref{lem:jacobi}(i)). For $k\ge1$ and $m,n\ge0$ let
\[
D^{\circ}_{m,k,n}(x,c):=\det\bigl(a_{i+j+m,k}(0,c)-x\,a_{i+j+m,k-1}(0,c)\bigr)_{0\le i,j\le n-1}\in\Z[x,c].
\]
By \eqref{eq:Dshift}, and by \eqref{eq:shiftm} for $y=0$ and Lemma~\ref{lem:jacobi}(iv), we have for all $m,n\ge0$
\begin{equation}\label{eq:DL}
\begin{aligned}
D_{m,k,n}(x,c)&=\det\bigl(L_x(X^{i+j}(X+c)^mP_k)\bigr)_{0\le i,j\le n-1},\\
D^{\circ}_{m,k,n}(x,c)&=\det\bigl(L_0(X^{i+j}(X+c)^mP_k)\bigr)_{0\le i,j\le n-1};
\end{aligned}
\end{equation}
in particular $D^{\circ}_{0,k,n}(x,c)=\cH_n(x;0)=\eps_k\,d_{n+k+1}$ by \eqref{eq:special}.

\subsection{Auxiliary lemmas}\label{ss:lemmas}

\begin{lemma}\label{lem:geronimus}
Let $A$ be a commutative ring and $y\in A$. Then for all $f\in A[X]$
\[
L_y\bigl((1-yX+y^2)f\bigr)=L_0(f).
\]
\end{lemma}

\begin{proof}
By Lemma~\ref{lem:jacobi}(iii) with $n=0$, every functional $L=L^{(s,\lambda)}$ satisfies $L(p_m)=h_m\,a^{(s,\lambda)}_{0,m}=\delta_{m,0}$. The polynomials $U_j$ $(j\ge0)$ form a basis of $A[X]$, and $XU_j=U_{j+1}+U_{j-1}$, hence
\[
(1-yX+y^2)U_j=(U_j-yU_{j-1})-y(U_{j+1}-yU_j)=P_{j,y}-yP_{j+1,y},
\]
so $L_y\bigl((1-yX+y^2)U_j\bigr)=\delta_{j,0}=L_0(U_j)$ for all $j\ge0$, by Lemma~\ref{lem:cheb}(iii).
\end{proof}

In the language of orthogonal polynomials (see e.g.\ \cite{Z}), Lemma~\ref{lem:geronimus} says that, for invertible $y$, $L_0$ is, up to the factor $-y$, the Christoffel transform of $L_y$ at the point $y+y^{-1}$, and that $-y\,L_y$ and $-y^{-1}L_{1/y}$ are Geronimus transforms of $L_0$ at the same point.

\begin{lemma}\label{lem:extra}
Let $\kk$ be a field, $P\in\kk[X]$ monic of degree $k\ge1$, $\xi\in\kk$, and $b_0,\dots,b_k\in\kk[X]$. Then
\[
P(\xi)\cdot\detx{(X-\xi)P}(b_0,\dots,b_k)=\sum_{j=0}^k(-1)^{k+j}\,b_j(\xi)\,\detx{P}\bigl(b_0,\dots,\widehat{b_j},\dots,b_k\bigr).
\]
\end{lemma}

\begin{proof}
Let $\beta_0,\dots,\beta_{k-1}$ be monic polynomials with $\deg\beta_i=i$. Then $\beta_0,\dots,\beta_{k-1},P$ give a monic basis of $\kk[X]/((X-\xi)P)$. For $f\in\kk[X]$ write $f=\sum_{i<k}\alpha_i(f)\beta_i+Pg_f$ with $\alpha_i(f)\in\kk$ and $g_f\in\kk[X]$ (division with remainder by $P$). Then $(\alpha_i(f))_{i<k}$ is the coordinate vector of $f$ in $\kk[X]/(P)$. Since $P\cdot(g_f-g_f(\xi))$ is divisible by $(X-\xi)P$, the coordinate vector of $f$ in $\kk[X]/((X-\xi)P)$ is $(\alpha_0(f),\dots,\alpha_{k-1}(f),g_f(\xi))$, and $P(\xi)g_f(\xi)=f(\xi)-\sum_{i<k}\alpha_i(f)\beta_i(\xi)$. Hence the left-hand side of the lemma is the determinant of the $(k+1)\times(k+1)$ matrix with the rows $\bigl(\alpha_0(b_j),\dots,\alpha_{k-1}(b_j),\,b_j(\xi)-\sum_i\alpha_i(b_j)\beta_i(\xi)\bigr)$, $0\le j\le k$. Adding $\beta_i(\xi)$ times the column with index $i$ to the last column, for all $i<k$, turns the last column into $(b_j(\xi))_{0\le j\le k}$, and expanding along the last column gives the right-hand side.
\end{proof}

The backward shifts are handled by the following classical duality between the Hankel determinants of a power series and of its reciprocal. It is stated in \cite[p.~5, (13)]{C24}, where it is attributed to \cite[Theorem~2]{AW} and \cite[Proposition~2.5]{C13}. The case $s=0$ is credited to Andrews and Wimp in \cite[Proposition~2.4]{C13}. We include a short proof.

\begin{lemma}\label{lem:reciprocal}
Let $A$ be a commutative ring, let $\mathbf u(z)=\sum_{r\ge0}u_rz^r\in A[[z]]$ with $u_0=1$, and let $\mathbf v(z)=\mathbf u(z)^{-1}=\sum_{r\ge0}v_rz^r$. Put $u_r=v_r=0$ for $r<0$. Then for all integers $s,t\ge0$
\[
\det\bigl(u_{i+j-s}\bigr)_{0\le i,j\le t+s}=(-1)^{t+\binom{s+1}2}\det\bigl(v_{i+j+s+2}\bigr)_{0\le i,j\le t-1}.
\]
\end{lemma}

\begin{proof}
Let $\varrho:=t+s+1$ and $\varrho':=t+\varrho$, and let $T_{\mathbf u}=(u_{i-j})_{0\le i,j\le\varrho'-1}$ and $T_{\mathbf v}=(v_{i-j})_{0\le i,j\le\varrho'-1}$. Both matrices are lower unitriangular, and $T_{\mathbf u}T_{\mathbf v}=I_{\varrho'}$, because $\sum_lu_{i-l}v_{l-j}$ is the coefficient of $z^{i-j}$ in $\mathbf u\mathbf v=1$. Consider the square matrices
\[
\mathcal U:=\bigl(u_{t+i-j}\bigr)_{0\le i,j\le\varrho-1},\qquad\mathcal V:=\bigl(v_{\varrho+i-j}\bigr)_{0\le i,j\le t-1},
\]
formed by the last $\varrho$ rows and the first $\varrho$ columns of $T_{\mathbf u}$, and by the last $t$ rows and the first $t$ columns of $T_{\mathbf v}$. Let $Z$ be the $\varrho'\times\varrho'$ matrix whose first $\varrho$ columns are the standard basis vectors $e_0,\dots,e_{\varrho-1}$ and whose last $t$ columns are the first $t$ columns of $T_{\mathbf v}$. Then $Z$ is block upper triangular with the diagonal blocks $I_\varrho$ and $\mathcal V$, so $\det Z=\det\mathcal V$. The first $\varrho$ columns of $T_{\mathbf u}Z$ are the first $\varrho$ columns of $T_{\mathbf u}$, and its last $t$ columns are $e_0,\dots,e_{t-1}$. Moving these $t$ columns to the front, which changes the determinant by the sign $(-1)^{t\varrho}$, produces a block upper triangular matrix with the diagonal blocks $I_t$ and $\mathcal U$. As $\det T_{\mathbf u}=1$, we get $\det\mathcal V=(-1)^{t\varrho}\det\mathcal U$.

Reversing the order of the columns turns $\mathcal U$ into the matrix $(u_{i+j-s})_{0\le i,j\le\varrho-1}$, with the sign $(-1)^{\binom\varrho2}$, and $\mathcal V$ into the matrix $(v_{i+j+s+2})_{0\le i,j\le t-1}$, with the sign $(-1)^{\binom t2}$. Now
\[
\tbinom\varrho2=\tbinom t2+t(s+1)+\tbinom{s+1}2\qquad\text{and}\qquad t\varrho=t^2+t(s+1),
\]
so the total sign is $(-1)^{t^2+\binom{s+1}2}=(-1)^{t+\binom{s+1}2}$.
\end{proof}

\begin{lemma}\label{lem:gf}
Let $\Gamma=\Gamma(z)\in\Z[c][[z]]$ be the unique formal power series with $\Gamma=1+cz\Gamma+z^2\Gamma^2$. Then for all $j\ge0$
\[
\sum_{n\ge0}a_{n,j}(0,c)\,z^n=z^j\,\Gamma^{j+1}\qquad\text{and}\qquad\sum_{n\ge0}a_{n,j}(x,c)\,z^n=\frac{z^j\,\Gamma^{j+1}}{1-xz\Gamma}.
\]
\end{lemma}

\begin{proof}
The first formula is \cite[Lemma~A.2]{Yu}. For the second one let $J$ and $J_0$ be the matrices of the weights $(c+x,c,c,\dots)$ and $(c,c,c,\dots)$, with $\lambda=(1,1,\dots)$. Then $J-J_0=x\,e_0e_0^{\T}$. From $J^n-J_0^n=\sum_{i=0}^{n-1}J_0^i(J-J_0)J^{n-1-i}$ we get $a_{n,j}(x,c)=a_{n,j}(0,c)+x\sum_{i=0}^{n-1}a_{i,0}(0,c)\,a_{n-1-i,j}(x,c)$, that is, $\mathcal G_j=z^j\Gamma^{j+1}+xz\,\Gamma\,\mathcal G_j$ for $\mathcal G_j:=\sum_na_{n,j}(x,c)z^n$. (This computation is part of the proof of \cite[Lemma~A.3]{Yu}.)
\end{proof}

\subsection{Backward shifts for arbitrary \texorpdfstring{$(x,c)$}{(x,c)}}\label{ss:backward}
Cigler remarks in \cite{C23} that for the determinants $D_{-m,k,n}(x,c)$ with $x\ne0$ there seems to be no analogue of the duality between backward and forward shifts which holds for $x=0$ \cite[Theorem~2]{C24}, apart from the point $(x,c)=(-1,2)$ and from the identities (32), (33) of \cite{C23} for $m=1$, $x=1$ (see Corollary~\ref{rem:C32}). The next proposition gives such an analogue for all $(x,c)$, in which the forward shifted determinants are those of the sequence $n\mapsto a_{n,k}(0,c)-x\,a_{n,k-1}(0,c)$.

\begin{proposition}[Backward shifts]\label{prop:backward}
Let $k\ge1$. For all $m\ge0$ and $n\ge0$, in $\Z[x,c]$,
\[
D_{-m,k,n+m+k+1}(x,c)=(-1)^{\binom{m+k+1}2}\,D^{\circ}_{m,k,n}(x,c),
\]
and $D_{-m,k,n}(x,c)=0$ for $1\le n\le m+k$.
\end{proposition}

For $x=0$ this is \cite[Theorem~2]{C24}, and for $m=0$ it is \eqref{eq:companion}. On the curve $xc+x^2+1=0$ the proposition is equivalent to \cite[Theorem~4]{C24a}, by \eqref{eq:Dcirc} and Corollary~\ref{cor:narayana} below. The proof follows those of \cite[Theorem~2]{C24} and \cite[Theorem~4]{C24a}, which also rest on Lemma~\ref{lem:reciprocal}. For $x=0$ the proposition also follows from \cite[Theorem~2.4]{LWZZ}, applied to the series $\Gamma$ of Lemma~\ref{lem:gf}, which satisfies $1+(cz-1)\Gamma+z^2\Gamma^2=0$, and to its power $\Gamma^{k+1}$. For $j\ge1$ the polynomials $\gamma_j$ in the proof below are the polynomials $L_j$ of \cite{LWZZ} for this series.

\begin{proof}
If $1\le n\le m+k$, the first row of the matrix $(a_{i+j-m,k}(x,c))_{0\le i,j\le n-1}$ consists of the elements $a_{j-m,k}$ with $j-m\le k-1$, which vanish by Lemma~\ref{lem:jacobi}(i) and the convention $a_{r,k}=0$ for $r<0$.

By Lemma~\ref{lem:gf},
\[
\mathbf u(z):=\sum_{r\ge0}a_{r+k,k}(x,c)\,z^r=\frac{\Gamma^{k+1}}{1-xz\Gamma}\in\Z[x,c][[z]],\qquad u_0=1.
\]
Dividing $\Gamma=1+cz\Gamma+z^2\Gamma^2$ by $\Gamma$ gives $\Gamma^{-1}+z^2\Gamma=1-cz$. Hence the series $\gamma_j:=\Gamma^{-j}+z^{2j}\Gamma^j$ $(j\ge0)$ satisfy $\gamma_0=2$, $\gamma_1=1-cz$ and $\gamma_{j+1}=\gamma_1\gamma_j-z^2\gamma_{j-1}$ for $j\ge1$, so that $\gamma_j$ is a polynomial in $z$ of degree at most $j$. Therefore
\[
\mathbf u(z)^{-1}=\Gamma^{-k-1}-xz\,\Gamma^{-k}=\bigl(\gamma_{k+1}-xz\,\gamma_k\bigr)-z^{k+2}\bigl(z^k\Gamma^{k+1}-x\,z^{k-1}\Gamma^k\bigr),
\]
and since $\gamma_{k+1}-xz\gamma_k$ is a polynomial in $z$ of degree at most $k+1$, Lemma~\ref{lem:gf} shows that the coefficients $v_r$ of $\mathbf u^{-1}$ satisfy
\[
v_r=-\bigl(a_{r-k-2,k}(0,c)-x\,a_{r-k-2,k-1}(0,c)\bigr)\qquad(r\ge k+2).
\]
Now apply Lemma~\ref{lem:reciprocal} with $s=m+k$ and $t=n$. We have $u_{i+j-s}=a_{i+j-m,k}(x,c)$ for all $i,j\ge0$ (both sides vanish if $i+j-m<k$), so the determinant on the left of Lemma~\ref{lem:reciprocal} is $D_{-m,k,n+m+k+1}(x,c)$. On the right we have
\[
v_{i+j+s+2}=-\bigl(a_{i+j+m,k}(0,c)-x\,a_{i+j+m,k-1}(0,c)\bigr),
\]
so the right-hand side of Lemma~\ref{lem:reciprocal} equals $(-1)^{n+\binom{m+k+1}2}(-1)^n\,D^{\circ}_{m,k,n}(x,c)$.
\end{proof}

\subsection{The curve \texorpdfstring{$xc+x^2+1=0$}{xc+x\textasciicircum2+1=0}}\label{ss:curve}
For $k=0$ the determinants $D_{1,0,n}(x,c)$ and $D_{2,0,n}(x,c)$ are given in \cite[(16)]{C23}, and Cigler adds: ``For $k>0$ and general $b$ we did not find nice formulas'' \cite[p.~5]{C23}. The next theorem expresses the determinants with the shifts $1$, $2$ and $-1$, for all $x$ and all $k\ge1$, by the determinants $d_n=D_{k,n}(x)$ when $c=-(x+x^{-1})$.

\begin{theorem}\label{thm:curve}
Let $k\ge1$ and $d_n=D_{k,n}(x)$. The following identities hold in $\Z[x,x^{-1}]$:
\begin{align}
D_{1,k,n}\bigl(x,-x-x^{-1}\bigr)&=(-1)^n\eps_k\,x^{-n}\,d_{n+k+1}&&(n\ge0),\label{eq:D1}\\
(x^2-1)\,D_{2,k,n}\bigl(x,-x-x^{-1}\bigr)&=x^{k+2}\,d_{n+1}-x^{-2n}\,d_{n+2k+2}&&(n\ge0),\label{eq:D2}\\
(x^2-1)\,D_{-1,k,n}\bigl(x,-x-x^{-1}\bigr)&=(-1)^{n-1}\eps_k\bigl(x^{n}\,d_{n-k-1}-x^{k+2-n}\,d_{n+k}\bigr)&&(n\ge k+1).\label{eq:D3}
\end{align}
Moreover, for all $m\ge1$ and $n\ge0$,
\begin{equation}\label{eq:D4}
D_{1-m,k,n+m+k}\bigl(x,-x-x^{-1}\bigr)=(-1)^{\binom{m+k}2}\,(-x)^n\,D_{m,k,n}\bigl(x,-x-x^{-1}\bigr).
\end{equation}
\end{theorem}

By Corollary~\ref{cor:narayana} below, the identity \eqref{eq:D4} is Cigler's duality \cite[Theorem~4]{C24a} for the Hankel determinants of the odd convolution powers of the Narayana polynomials, written in a different notation (see \S\ref{ss:related} for the case $x=-1$). Its case $m=1$ is \eqref{eq:D1}, since $D_{0,k,n}$ does not depend on $c$. We include the proof of \eqref{eq:D4} because it is a special case of Proposition~\ref{prop:backward}. For an arbitrary power series, Han and Yang \cite[Theorems~2.1--2.3]{HY} express the Hankel determinants with the shifts $1$, $2$ and $-1$ through its Hankel continued fraction. The identities \eqref{eq:D1}--\eqref{eq:D3} express them, for the sequence $n\mapsto a_{n,k}(x,-x-x^{-1})$, through the unshifted determinants $d_n$.

\begin{proof}
Put $\zeta:=x+x^{-1}\in K$. Since both sides of the identities \eqref{eq:D1}--\eqref{eq:D4} lie in $\Z[x,x^{-1}]$, we may compute all determinants in $K$. Note that
\begin{equation}\label{eq:zeta}
X-\zeta=-x^{-1}(1-xX+x^2),\qquad U_m(\zeta)=\frac{x^{m+1}-x^{-m-1}}{x-x^{-1}}\quad(m\in\Z),\qquad P_k(\zeta)=x^{-k},
\end{equation}
by \eqref{eq:cheb} at $t=x$ and Lemma~\ref{lem:cheb}(ii). By Lemma~\ref{lem:geronimus} with $y=x$ we have $L_x\bigl((X-\zeta)f\bigr)=-x^{-1}L_0(f)$ for all $f\in K[X]$. Hence, by \eqref{eq:DL}, for all $m,n\ge0$
\begin{equation}\label{eq:Dcirc}
D_{m+1,k,n}(x,-\zeta)=\det\bigl(L_x(X^{i+j}(X-\zeta)^{m+1}P_k)\bigr)_{i,j<n}=(-x)^{-n}\,D^{\circ}_{m,k,n}(x,-\zeta).
\end{equation}

\emph{Proof of \eqref{eq:D1}.} This is \eqref{eq:Dcirc} for $m=0$, since $D^\circ_{0,k,n}=\eps_k\,d_{n+k+1}$.

\emph{Proof of \eqref{eq:D2}.} By \eqref{eq:Dcirc} for $m=1$, \eqref{eq:DL}, and \eqref{eq:christoffel} for the weights $s=(0,0,\dots)$, $\lambda=(1,1,\dots)$ (whose polynomials $p_m$ are the $U_m$) and $f=(X-\zeta)P_k$,
\begin{align*}
D_{2,k,n}(x,-\zeta)&=(-x)^{-n}\det\bigl(L_0(X^{i+j}(X-\zeta)P_k)\bigr)_{i,j<n}\\
&=(-x)^{-n}(-1)^{n(k+1)}\detx{(X-\zeta)P_k}(U_n,\dots,U_{n+k}).
\end{align*}
We evaluate the last determinant by Lemma~\ref{lem:extra} with $\xi=\zeta$ and $b_j=U_{n+j}$. Since the class of $U_{n+l}$ in $R$ is $F_{n+1+l}$, we have $\detR(U_n,\dots,\widehat{U_{n+j}},\dots,U_{n+k})=M_j(n+1)$, and by \eqref{eq:zeta} and the definition of $\cW_{n+1}$
\begin{align*}
\sum_{j=0}^k(-1)^jU_{n+j}(\zeta)M_j(n+1)&=\frac{x^{n+1}\cW_{n+1}(x)-x^{-n-1}\cW_{n+1}(x^{-1})}{x-x^{-1}}\\
&=\frac{x^{n+1}\cW_{n+1}(x)-x^{-n-1-k}\cW^*_{n+1}(x)}{x-x^{-1}}.
\end{align*}
By \eqref{eq:Wspecial} and \eqref{eq:dw}, $\cW_{n+1}(x)=\eps_kw(n+1-k)=(-1)^{(n+1)k}d_{n+1}$ and $\cW^*_{n+1}(x)=(-1)^k\eps_kw(n+k+2)=(-1)^{(n+1)k}d_{n+2k+2}$. With $P_k(\zeta)=x^{-k}$, Lemma~\ref{lem:extra} gives
\[
\detx{(X-\zeta)P_k}(U_n,\dots,U_{n+k})=(-1)^{nk}\,\frac{x^{n+k+1}d_{n+1}-x^{-n-1}d_{n+2k+2}}{x-x^{-1}},
\]
and therefore
\[
D_{2,k,n}(x,-\zeta)=x^{-n}\,\frac{x^{n+k+1}d_{n+1}-x^{-n-1}d_{n+2k+2}}{x-x^{-1}}=\frac{x^{k+2}d_{n+1}-x^{-2n}d_{n+2k+2}}{x^2-1}.
\]

\emph{Proof of \eqref{eq:D4}.} By \eqref{eq:Dcirc} and Proposition~\ref{prop:backward}, for all $m\ge0$,
\[
D_{-m,k,n+m+k+1}(x,-\zeta)=(-1)^{\binom{m+k+1}2}\,(-x)^n\,D_{m+1,k,n}(x,-\zeta),
\]
which is \eqref{eq:D4} with $m+1$ in place of $m$.

\emph{Proof of \eqref{eq:D3}.} For $n=k+1$ the left-hand side vanishes by Proposition~\ref{prop:backward} (with $m=1$), and the right-hand side is $(-1)^k\eps_k(x^{k+1}d_0-x\,d_{2k+1})=0$, because $d_0=1$ and $d_{2k+1}=x^k$; the latter follows from \eqref{eq:special} for $n=0$, where $\cH_0=1$ and hence $\cH^*_0(x;y)=y^k$. Let $n\ge k+2$. By \eqref{eq:D4} with $m=2$ and with $n-k-2$ in place of $n$, and by \eqref{eq:D2},
\begin{align*}
(x^2-1)\,D_{-1,k,n}(x,-\zeta)&=(-1)^{\binom{k+2}2}(-x)^{n-k-2}\,(x^2-1)\,D_{2,k,n-k-2}(x,-\zeta)\\
&=(-1)^{\binom{k+2}2+n-k}\,x^{n-k-2}\bigl(x^{k+2}d_{n-k-1}-x^{-2(n-k-2)}d_{n+k}\bigr)\\
&=(-1)^{\binom{k+2}2+n-k}\bigl(x^nd_{n-k-1}-x^{k+2-n}d_{n+k}\bigr),
\end{align*}
and $(-1)^{\binom{k+2}2+n-k}=(-1)^{n-1}\eps_k$, because $\binom{k+2}2=\binom{k+1}2+k+1$.
\end{proof}

\begin{remark}\label{rem:curve}
(1) By Proposition~\ref{prop:backward}, $D_{-1,k,n}(x,c)=0$ for $1\le n\le k+1$ and all $c$, so \eqref{eq:D3} covers all non-trivial cases.

(2) For $k=0$ the identities \eqref{eq:D1}--\eqref{eq:D4} hold with $d_n=1$ and $\eps_0=1$. The proofs above apply with $P_0=1$ and $D^{\circ}_{m,0,n}(x,c):=D_{m,0,n}(0,c)$, using $\det\bigl(L_0(X^{i+j})\bigr)_{i,j<n}=1$ in the proof of \eqref{eq:D1} and $\det\bigl(L_0(X^{i+j}(X-\zeta))\bigr)_{i,j<n}=(-1)^nU_n(\zeta)$ (Lemma~\ref{lem:christoffel}) in the proof of \eqref{eq:D2}. In the proof of Proposition~\ref{prop:backward} one has $\mathbf u^{-1}=1-(c+x)z-z^2\Gamma$, so that $v_r=-a_{r-2,0}(0,c)$ for $r\ge2$, which gives $D_{-m,0,n+m+1}(x,c)=(-1)^{\binom{m+1}2}\,D_{m,0,n}(0,c)$ $(m,n\ge0)$, that is, \cite[Theorem~1]{C24} for the weights $(c+x,c,c,\dots)$. For $n=1$ the right-hand side of \eqref{eq:D3} is $x\,d_0-x\,d_1=0$. The resulting formulas are $D_{1,0,n}(x,-\zeta)=(-x)^{-n}$ and $D_{2,0,n}(x,-\zeta)=\sum_{j=0}^nx^{-2j}$, in agreement with \cite[(16)]{C23}.

(3) For $x=1$, $c=-2$ the identity \eqref{eq:D1} agrees with \cite[Theorems~1.6 and~1.7]{CS}. At $x=\pm1$ the identities \eqref{eq:D2} and \eqref{eq:D3} have to be divided by $x^2-1$ in $\Z[x,x^{-1}]$ before $x$ is specialized.
\end{remark}

\begin{corollary}[Theorem~\ref{thm:B} and the Desnanot--Jacobi identity]\label{cor:DJB}
Let $k\ge1$. For all $m\in\Z$ and $n\ge1$ the Desnanot--Jacobi identity \cite[Proposition~10]{K99} gives
\begin{equation}\label{eq:DJ}
D_{m,k,n+1}\,D_{m+2,k,n-1}=D_{m,k,n}\,D_{m+2,k,n}-D_{m+1,k,n}^{\,2}\qquad\text{in }\Z[x,c].
\end{equation}
For $m=0$ and $c=-(x+x^{-1})$, after inserting \eqref{eq:D1} and \eqref{eq:D2} and multiplying by $(x^2-1)x^{2n}$, the identity \eqref{eq:DJ} becomes
\[
d_{n+2k+2}\,d_n-x^2\,d_{n+2k+1}\,d_{n+1}=(1-x^2)\,d_{n+k+1}^{\,2},
\]
that is, \eqref{eq:main} with $n$ replaced by $n+2k+2$. Moreover, on the curve $xc+x^2+1=0$ the determinants $D_{1,k,n}$ satisfy the recurrence of Theorem~\ref{thm:B}:
\[
D_{1,k,n}D_{1,k,n-2k-2}-x^2D_{1,k,n-1}D_{1,k,n-2k-1}=(1-x^2)D_{1,k,n-k-1}^{\,2}\qquad(n\ge2k+2,\ c=-x-x^{-1}).
\]
\end{corollary}

\begin{proof}
The matrix $(a_{i+j+m,k})_{0\le i,j\le n}$ has the minors $D_{m,k,n+1}$ (the full determinant), $D_{m,k,n}$ (last row and column deleted), $D_{m+2,k,n}$ (first row and column deleted), $D_{m+1,k,n}$ (first row and last column deleted, or last row and first column deleted) and $D_{m+2,k,n-1}$ (first and last rows and columns deleted); so \eqref{eq:DJ} is the Desnanot--Jacobi identity for this matrix. For $m=0$ and $c=-(x+x^{-1})$, multiplying \eqref{eq:DJ} by $(x^2-1)x^{2n}$ and inserting \eqref{eq:D1}, \eqref{eq:D2} gives
\[
d_{n+1}\bigl(x^{2n+k+2}d_n-x^2d_{n+2k+1}\bigr)=d_n\bigl(x^{2n+k+2}d_{n+1}-d_{n+2k+2}\bigr)-(x^2-1)d_{n+k+1}^{\,2},
\]
which is the asserted identity. For $n=0$, \eqref{eq:D2} reads $x^2-1=x^{k+2}d_1-d_{2k+2}$ (as $D_{2,k,0}=1$), and $d_1=0$. The last assertion follows from \eqref{eq:D1} and \eqref{eq:main} for $n+k+1$ in place of $n$: the signs $(-1)^n\eps_k$ cancel in each of the three products, and the powers of $x$ are $x^{-(2n-2k-2)}$ in all three terms.
\end{proof}

The proofs of \eqref{eq:D1} and \eqref{eq:D2} rest on Lemmas~\ref{lem:christoffel}, \ref{lem:expansion}, \ref{lem:Lambda}, \ref{lem:geronimus} and~\ref{lem:extra}, through \eqref{eq:C1}, \eqref{eq:Wspecial} and \eqref{eq:special}, and they are independent of Lemma~\ref{lem:HM} and Proposition~\ref{prop:bilinear}. So Corollary~\ref{cor:DJB} gives a third proof of Theorem~\ref{thm:B}. The remaining case $n=2k+2$ of \eqref{eq:main} reads $d_{2k+2}=1-x^2$, because $d_0=1$, $d_1=0$ and $d_{k+1}^{\,2}=1$, and it is \eqref{eq:D2} for $n=0$. Similarly, \eqref{eq:DJ} for $m=-1$, together with \eqref{eq:D1} and \eqref{eq:D3}, is equivalent to \eqref{eq:main}. In the lattice of the determinants $D_{m,k,n}$, the identity \eqref{eq:DJ} is the discrete-time Toda equation in the form \cite[(2.8)]{SNvdK}.

\begin{remark}[The Desnanot--Jacobi identity and \eqref{eq:hankel2}]\label{rem:DJ}
For $n\ge1$ the identity \eqref{eq:hankel2} can also be derived from the Desnanot--Jacobi identity. Let $\zeta_y:=y+y^{-1}$. By Lemma~\ref{lem:geronimus}, applied to $y$ and to $y^{-1}$, we have $L_y\bigl((X-\zeta_y)f\bigr)=-y^{-1}L_0(f)$ and $y^{-2}L_{1/y}\bigl((X-\zeta_y)f\bigr)=-y^{-1}L_0(f)$ for all $f$. Hence the two sequences $g'_r:=L_y((X-\zeta_y)^rP_k)$ and $g_r:=y^{-2}L_{1/y}((X-\zeta_y)^rP_k)$ $(r\ge0)$ agree for $r\ge1$, where both are equal to $-y^{-1}L_0((X-\zeta_y)^{r-1}P_k)$, and $g'_0=g_0+\varsigma$ with $\varsigma=(1-y^{-2})P_k(\zeta_y)$. For any sequence $(g_r)_{r\ge0}$ and any scalar $\varsigma$, the linearity of the determinant in the entry $(0,0)$ and the Desnanot--Jacobi identity give
\begin{multline*}
\det\bigl(g_{i+j}+\varsigma\delta_{i+j,0}\bigr)_{i,j<n}\,\det(g_{i+j})_{i,j<n-1}-\det(g_{i+j})_{i,j<n}\,\det\bigl(g_{i+j}+\varsigma\delta_{i+j,0}\bigr)_{i,j<n-1}\\
=\varsigma\,\det(g_{i+j+1})_{i,j<n-1}^{\,2}.
\end{multline*}
By Lemma~\ref{lem:jacobi}(iv), the Hankel determinants of a functional may be computed in the basis $(X-\zeta_y)^r$. So the four determinants on the left are $\cH_n(x;y)$, $y^{-2(n-1)-k}\cH^*_{n-1}(x;y)$, $y^{-2n-k}\cH^*_n(x;y)$ and $\cH_{n-1}(x;y)$, and the determinant on the right is $(-y)^{-(n-1)}\cH_{n-1}(x;0)$. Multiplying by $y^{2n+k}$ and using $\cH_{n-1}(x;0)^2=D_{k,n+k}^{\,2}$ and $y^k(y^2-1)P_k(\zeta_y)=Q(y)$ gives \eqref{eq:hankel2}. The same computation applies to an arbitrary monic polynomial $P$ in place of $P_k$ and gives \eqref{eq:spectral} for $\ell\ge1$. Thus an identity of the shape \eqref{eq:hankel2} holds for every monic polynomial. What is special about $P_k$ is \eqref{eq:special}, by which, for $y=x$, all three sequences of Hankel determinants involved are shifts of the single sequence $(d_n)$. This comes from the operator $\Lambda$.
\end{remark}
\subsection{Narayana weights}\label{ss:narayana}
The curve $xc+x^2+1=0$ is the curve of the \emph{Narayana weights}. Let $q$ be an indeterminate and let $\Nar_{n,j}(q):=a^{(s,\lambda)}_{n,j}\in\Z[q]$ for $s=(1,1+q,1+q,\dots)$ and $\lambda=(q,q,\dots)$, as in \S\ref{ss:results}. For fixed $k$ let
\[
\cN_{m,n}(q):=\det\bigl(\Nar_{i+j+m,k}(q)\bigr)_{0\le i,j\le n-1}\quad(\Nar_{r,k}:=0\text{ for }r<0),\qquad\fn_n:=\cN_{0,n}(q).
\]

\begin{corollary}[Narayana form]\label{cor:narayana}
Let $k\ge0$.
\begin{enumerate}
\item[(i)] Let $\nu_0:=\sum_{n\ge0}\Nar_{n,0}(q)\,z^n$ and $\nu_1:=1-q+q\,\nu_0$ in $\Z[q][[z]]$. Then $\nu_0=1+z\,\nu_0\nu_1$ and
\[
\sum_{n\ge0}\Nar_{n,k}(q)\,z^n=z^k\,\nu_0^{\,k+1}\nu_1^{\,k}.
\]
Thus $\nu_0$ is the generating function of the Narayana polynomials, and the $k$-th column of the triangle $(\Nar_{n,j}(q))$ is the sequence of the $(2k+1)$-th convolution powers of the Narayana polynomials in the sense of \cite[Section~3]{C24a}, preceded by $k$ zeros.
\item[(ii)] In $\Z[q,q^{-1}]$ we have
\begin{align*}
\fn_n\,\fn_{n-2k-2}-q^{2k+2}\,\fn_{n-1}\,\fn_{n-2k-1}&=(1-q)\,q^{(k+1)^2}\,\fn_{n-k-1}^{\,2},\\
\cN_{1,n}&=\eps_k\,q^{-n(k+1)}\,\fn_{n+k+1},\\
(q-1)\,\cN_{2,n}&=(-1)^kq^{k+1}\,\fn_{n+1}-q^{-(k+1)(2n+k+1)}\,\fn_{n+2k+2},\\
(q-1)\,\cN_{-1,n}&=\eps_k\bigl((-1)^kq^{-(k+1)(n-1)}\,\fn_{n+k}-q^{(k+1)(n-k-1)}\,\fn_{n-k-1}\bigr),\\
\cN_{1-m,n+m+k}&=(-1)^{\binom{m+k}2}\,q^{n(k+1)}\,\cN_{m,n},
\end{align*}
the first identity for $n\ge2k+2$, the second and the third for $n\ge0$, the fourth for $n\ge k+1$, and the last one for all $m\ge1$ and $n\ge0$.
\end{enumerate}
\end{corollary}

In the notation of \cite{C24a}, where our $z$ and $q$ are called $x$ and $t$, $\Nar_{n,k}(q)$ is $C_{2k+1,n-k}(q)$ and $\cN_{m,n}(q)$ is the determinant denoted there by $D_{2k+1,m-k,q}(n)$ \cite[(28), (42)]{C24a}. By (i), the last identity of (ii) is \cite[Theorem~4]{C24a} with $k+1$ in place of the $k$ there, and the second identity is its case $m=1$. The first identity is Theorem~\ref{thm:B}. Thus the conjecture~(18) of \cite{C23} is a quadratic recurrence for the Hankel determinants $\fn_n=D_{2k+1,-k,q}(n)$ of the odd convolution powers of the Narayana polynomials, shifted backward by $k$. By the second identity the determinants $\cN_{1,n}$, that is, the shift $1-k$, satisfy the same recurrence, since the powers of $q$ cancel. For $k\ge2$ the unshifted Hankel determinants $\cN_{k,n}(q)=D_{2k+1,0,q}(n)$ do not satisfy the recurrence in general. For instance, a direct computation shows that it fails for $k=2$, $q=3/5$ and $n=7$. For $q=1$ (the Catalan numbers) the right-hand side of the first identity vanishes. This agrees with \cite[(21)]{C24a}, according to which the determinants $\fn_n$ are then periodic up to sign. The even convolution powers correspond in the same way to the triangle of constant type $s=(1+q,1+q,\dots)$, $\lambda=(q,q,\dots)$, that is, to $a_{n,k}(0,c)$. Their Hankel determinants are studied in \cite{CK11,C24a,H25}. For $k=0$ we have $\fn_n=q^{\binom n2}$, and the second and the third identity reduce to the evaluations $\cN_{1,n}=q^{\binom n2}$ and $\cN_{2,n}=q^{\binom n2}(1+q+\dots+q^n)$ of the shifted Hankel determinants of the Narayana polynomials. For $k=1$ the determinants $\cN_{1,n}(q)$ are evaluated in closed form in \cite[Theorem~5.2]{C18}, and by the second identity this also evaluates $\fn_n$ for $k=1$.

\begin{proof}
A Motzkin path from $(0,0)$ to $(n,j)$ with $\beta$ down steps has $j+\beta$ up steps. If it has $\ell_0$ level steps at height $0$ and $\ell_+$ level steps at positive height, then $(j+\beta)+\beta+\ell_0+\ell_+=n$, so $\ell_0+\ell_++2\beta=n-j$. Its weight in $a_{n,j}(x,-x-x^{-1})$ is $(-x^{-1})^{\ell_0}(-x-x^{-1})^{\ell_+}$, and its weight in $\Nar_{n,j}(x^2)$ is $(1+x^2)^{\ell_+}x^{2\beta}$, which is $(-x)^{\ell_0+\ell_++2\beta}=(-x)^{n-j}$ times the former. Hence
\begin{equation}\label{eq:Nar-a}
\Nar_{n,j}(x^2)=(-x)^{n-j}\,a_{n,j}(x,-x-x^{-1})\qquad(n,j\ge0).
\end{equation}
Both parts of the corollary follow from this relation.

(i) Let $c=-(x+x^{-1})$ and $\widetilde\Gamma(z):=\Gamma(-xz)\in\Z[x,x^{-1}][[z]]$, with the series $\Gamma$ of Lemma~\ref{lem:gf}. Then $\widetilde\Gamma=1+(1+x^2)z\widetilde\Gamma+x^2z^2\widetilde\Gamma^2$, that is,
\[
(1+x^2z\widetilde\Gamma)(1+z\widetilde\Gamma)=\widetilde\Gamma,
\]
and by \eqref{eq:Nar-a} and Lemma~\ref{lem:gf}, $\sum_{n\ge0}\Nar_{n,k}(x^2)z^n=z^k\,\widetilde\Gamma^{k+1}/(1+x^2z\widetilde\Gamma)$. For $k=0$ this reads $\nu_0\big|_{q=x^2}=\widetilde\Gamma/(1+x^2z\widetilde\Gamma)=1+z\widetilde\Gamma$. Hence $1+x^2z\widetilde\Gamma=1-x^2+x^2\nu_0=\nu_1$ and $\widetilde\Gamma=\nu_0\nu_1$ at $q=x^2$, so that $\nu_0=1+z\,\nu_0\nu_1$ and $\sum_n\Nar_{n,k}(x^2)z^n=z^k\,\widetilde\Gamma^k\,\nu_0=z^k\,\nu_0^{\,k+1}\nu_1^{\,k}$ at $q=x^2$. Since $\Z[q]\to\Z[x,x^{-1}]$, $q\mapsto x^2$, is injective, both identities hold in $\Z[q][[z]]$. The relations $\nu_1=1-q+q\,\nu_0$ and $\nu_0=1+z\,\nu_0\nu_1$ are the relations (25) and (26) of \cite{C24a}, with $(z,q)$ in place of $(x,t)$, for the generating function $c_0(x,t)$ of the Narayana polynomials and the series $c_1(x,t)$ considered there, and they determine a power series $\nu_0$ uniquely. Hence $\nu_0$ and $\nu_1$ are these two series, and $\nu_0^{\,k+1}\nu_1^{\,k}$ is the series $c(x,t)^{(2k+1)}$ of \cite[(27)]{C24a}.

(ii) By \eqref{eq:Nar-a},
\[
\cN_{m,n}(x^2)=(-x)^{n(n-1)+n(m-k)}\,D_{m,k,n}(x,-x-x^{-1})\qquad(m\in\Z,\ n\ge0).
\]
In particular $d_N=(-x)^{-N(N-1-k)}\fn_N(x^2)$. Inserting this into \eqref{eq:main} and \eqref{eq:D1}--\eqref{eq:D4} $($for $k=0$ see Remark~\ref{rem:curve}(2)$)$ and collecting the powers of $-x$, all of which have even exponents, gives the five identities with $q=x^2$. Since $\Z[q,q^{-1}]\to\Z[x,x^{-1}]$, $q\mapsto x^2$, is injective, they hold in $\Z[q,q^{-1}]$. For example, with $\theta(N):=N(N-1-k)$ we have $\theta(n)+\theta(n-2k-2)=2\theta(n-k-1)+2(k+1)^2$ and $\theta(n-1)+\theta(n-2k-1)=2\theta(n-k-1)+2k^2$, which turns \eqref{eq:main} into the first identity. Similarly, $n(n-k)-\theta(n+k+1)=-n(2k+1)$ turns \eqref{eq:D1} into the second one, and for the last one the exponent of $-x$ is $(n+m+k)n+n-n(n-1)-n(m-k)=2n(k+1)$.
\end{proof}

\subsection{The shifts \texorpdfstring{$1$ and $-1$ for arbitrary $(x,c)$}{1 and -1 for arbitrary (x,c)}}\label{ss:shift1}
For general $c$ the first shifted determinant is governed by the polynomials $\cW_m(y)$ of \S\ref{ss:spectral}. Let $\phi_m(c):=U_{m-1}(c)$ $(m\in\Z)$ be the Fibonacci polynomials (denoted by $F_m(c)$ in \cite{C23}).

\begin{proposition}\label{prop:shift1}
Let $k\ge1$. For all $n\ge0$, in $\Z[x,c]$,
\begin{equation}\label{eq:shift1}
\bigl(\phi_{k+1}(c)+x\phi_k(c)\bigr)\,D_{1,k,n}(x,c)=(-1)^{nk}\eps_k\sum_{i=0}^kM_i(n-k)\,\bigl(\phi_{n+i+1}(c)+x\phi_{n+i}(c)\bigr),
\end{equation}
where $M_i(m)\in\Z[x]$ are the coefficients of $\cW_m(y)=\sum_{i=0}^k(-y)^iM_i(m)$. Equivalently, for $c=-(\tau+\tau^{-1})$ with a new variable $\tau$,
\begin{equation}\label{eq:shift1tau}
\tau^{-k-1}Q(\tau)\,D_{1,k,n}\bigl(x,-\tau-\tau^{-1}\bigr)=(-1)^{nk+n+k}\eps_k\Bigl(\tau^n(\tau-x)\,\cW_{n-k}(\tau)-\tau^{-n}(\tau^{-1}-x)\,\cW_{n-k}(\tau^{-1})\Bigr)
\end{equation}
in $\Z[x][\tau,\tau^{-1}]$. For $k=1$, \eqref{eq:shift1} gives
\[
(x+c)\,D_{1,1,n}(x,c)=(-1)^{n+1}\Bigl(\phi_n(x)\bigl(\phi_{n+1}(c)+x\phi_n(c)\bigr)+\phi_{n-1}(x)\bigl(\phi_{n+2}(c)+x\phi_{n+1}(c)\bigr)\Bigr).
\]
\end{proposition}

\begin{proof}
We work over the field $\Q(x,c)$. By \eqref{eq:DL}, by \eqref{eq:christoffel} for the weights $(x,0,0,\dots)$, $(1,1,\dots)$ and $f=(X+c)P_k$, and by Lemma~\ref{lem:extra} for $\xi=-c$, $b_j=P_{n+j}$,
\[
P_k(-c)\,D_{1,k,n}(x,c)=(-1)^{n(k+1)}\sum_{j=0}^k(-1)^{k+j}P_{n+j}(-c)\,\detR\bigl(P_n,\dots,\widehat{P_{n+j}},\dots,P_{n+k}\bigr).
\]
By \eqref{eq:GLambda} the class of $P_{n+l}$ is $\Lambda F_{n-k+l}$, so the last determinant equals $\det(\Lambda)M_j(n-k)=\eps_kM_j(n-k)$. Finally $U_m(-c)=(-1)^mU_m(c)=(-1)^m\phi_{m+1}(c)$, so
\begin{equation}\label{eq:Pphi}
P_m(-c)=(-1)^m\bigl(\phi_{m+1}(c)+x\phi_m(c)\bigr)\qquad(m\ge0),
\end{equation}
and the signs combine to $(-1)^{n(k+1)+k+n+k}=(-1)^{nk}$. This proves \eqref{eq:shift1}.

For \eqref{eq:shift1tau} put $c=-(\tau+\tau^{-1})$ in \eqref{eq:shift1}. By \eqref{eq:Pphi} and \eqref{eq:cheb},
\[
(\tau-\tau^{-1})\bigl(\phi_{m+1}(c)+x\phi_m(c)\bigr)=(-1)^m\bigl(\tau^m(\tau-x)-\tau^{-m}(\tau^{-1}-x)\bigr)\qquad(m\ge0),
\]
and for $m=k$ the right-hand side equals $(-1)^k\tau^{-k-1}Q(\tau)$. Multiplying \eqref{eq:shift1} by $\tau-\tau^{-1}$ and using $\sum_i(-\tau^{\pm1})^iM_i(n-k)=\cW_{n-k}(\tau^{\pm1})$ gives \eqref{eq:shift1tau}.

For $k=1$ we have $R=K[X]/(X-x)$, so $F_j$ is $U_{j-1}(x)$ times the class of $1$; hence $M_0(m)=U_m(x)=\phi_{m+1}(x)$ and $M_1(m)=U_{m-1}(x)=\phi_m(x)$. Moreover $\eps_1=-1$ and $\phi_2(c)+x\phi_1(c)=c+x$.
\end{proof}

By \eqref{eq:Wspecial}, four values of the polynomial $\cW_m(y)$, which has degree at most $k$, are known in terms of the sequence $(d_n)$, namely its values at $y=x$ and (after reversal) at $y=x^{-1}$, its constant term and its leading coefficient. For $k\le3$ these determine $\cW_m$ over $K$ (for $k=3$ the coefficients $M_1,M_2$ solve a linear system with determinant $x-x^{-1}$), so that \eqref{eq:shift1} then expresses $D_{1,k,n}(x,c)$, as an element of $\Z[x,c]$, through at most four terms of $(d_n)$. On the curve $xc+x^2+1=0$, that is, for $\tau=x$, the first term on the right of \eqref{eq:shift1tau} vanishes and only the value $\cW_{n-k}(x^{-1})$ survives, which gives \eqref{eq:D1} again.

The same method, combined with Proposition~\ref{prop:backward}, gives the determinants with the backward shift $1$. For $x=1$ the sums collapse, and one obtains closed formulas which complement the results of \cite{CS} on $D_{1,k,n}(1,c)$.

\begin{corollary}[Backward shift $1$]\label{cor:back1}
Let $k\ge1$ and write $\phi_m$ for $\phi_m(c)$.
\begin{enumerate}
\item[(i)] For all $n\ge0$, in $\Z[x,c]$,
\[
\bigl(\phi_{k+1}+x\phi_k\bigr)\,D_{-1,k,n+k+2}(x,c)=(-1)^{\binom{k+2}2+nk}\sum_{j=0}^kM_j(n+1)\,\phi_{n+j+1}.
\]
\item[(ii)] Let $x=1$, and put $\Phi_\ell:=\phi_{(2k+1)\ell+k+1}+\phi_{(2k+1)\ell+k}$ for $\ell\in\Z$ and $Y:=\Phi_0=\phi_{k+1}+\phi_k$. Then for all $\ell\ge0$, in $\Z[c]$,
\begin{align*}
Y\cdot D_{-1,k,(2k+1)\ell+k+1}(1,c)&=(-1)^{\binom{k+2}2+k}\,\phi_{(2k+1)\ell},\\
Y\cdot D_{-1,k,(2k+1)\ell+k+2}(1,c)&=(-1)^{\binom{k+2}2}\,\Phi_\ell,\\
Y\cdot D_{-1,k,(2k+1)(\ell+1)}(1,c)&=-\Phi_\ell,\\
Y\cdot D_{-1,k,(2k+1)(\ell+1)+1}(1,c)&=(-1)^{k+1}\,\phi_{(2k+1)(\ell+1)},
\end{align*}
and $D_{-1,k,n}(1,c)=0$ for all $n\ge1$ which are not congruent to $0$, $1$, $k+1$ or $k+2$ modulo $2k+1$. $($For $k=1$ the second and the third formula concern the same sizes $3\ell+3$, and they agree, because $\binom32$ is odd.$)$
\end{enumerate}
\end{corollary}

\begin{proof}[Proof of {\rm(i)}]
We work over $\Q(x,c)$. By \eqref{eq:DL}, by \eqref{eq:christoffel} for the weights $s=(0,0,\dots)$, $\lambda=(1,1,\dots)$ and the polynomial $(X+c)P_k$, and by Lemma~\ref{lem:extra} for $\xi=-c$ and $b_j=U_{n+j}$,
\[
P_k(-c)\,D^{\circ}_{1,k,n}(x,c)=(-1)^{n(k+1)}\sum_{j=0}^k(-1)^{k+j}\,U_{n+j}(-c)\,\detR\bigl(U_n,\dots,\widehat{U_{n+j}},\dots,U_{n+k}\bigr),
\]
and the last determinant is $M_j(n+1)$. Inserting $U_m(-c)=(-1)^m\phi_{m+1}(c)$ and \eqref{eq:Pphi} for $m=k$ we get
\[
\bigl(\phi_{k+1}+x\phi_k\bigr)\,D^{\circ}_{1,k,n}(x,c)=(-1)^{nk}\sum_{j=0}^kM_j(n+1)\,\phi_{n+j+1},
\]
and (i) follows from Proposition~\ref{prop:backward} with $m=1$.
\end{proof}

The proof of (ii) is given in Appendix~\ref{app:C32}.

\begin{corollary}[The identities (32) and (33) of \cite{C23}]\label{rem:C32}
Let $k\ge1$. For all $n\ge0$, in $\Z[c]$,
\begin{align*}
D_{-1,k,(2k+1)(n+1)}(1,c)&=-D_{1,k,(2k+1)n}(1,c),\\
D_{-1,k,(2k+1)n+k+2}(1,c)&=-D_{1,k,(2k+1)n+k}(1,c),\\
(c+2)\,D_{-1,k,(2k+1)(n+1)+1}(1,c)&=-D_{1,k,(2k+1)n+2k}(1,c),\\
(c+2)\,D_{-1,k,(2k+1)(n+1)+k+1}(1,c)&=-D_{1,k,(2k+1)n+k+1}(1,c),
\end{align*}
and $D_{-1,k,n}(1,c)=0$ for all $n\ge1$ that are not congruent to $0$, $1$, $k+1$ or $k+2$ modulo $2k+1$.
\end{corollary}

The first two identities are (32) of \cite{C23}, and the last two are (33) of \cite{C23} with $n$ replaced by $n+1$ on the left-hand sides. In \cite[(33)]{C23} the sizes on the left-hand sides are printed as $(2k+1)n+1$ and $(2k+1)n+k+1$, but the numerical example given there for $k=1$, $c=3$ agrees with the indexing used here. In \cite{C23} these identities and the vanishing of the other determinants are stated as observations. The proof, in Appendix~\ref{app:C32}, also gives another proof of the evaluation of $D_{1,k,n}(1,c)$ in \cite[Theorem~1.7]{CS}.

\section{An analogue of odd orthogonal type and paths with a weighted return step}\label{sec:typeB}

In this section $k\ge1$ and $K=\Q(x)$. Recall that $W_m=U_m+U_{m-1}$ $(m\in\Z)$ are the Chebyshev polynomials of the fourth kind, and put
\[
r_m:=W_m-xW_{m-1}=P_m+P_{m-1}\qquad(m\in\Z).
\]
Thus $W_0=1$, $W_{m+1}=XW_m-W_{m-1}$ and $W_{-m-1}=-W_m$ for all $m\in\Z$ (Lemma~\ref{lem:cheb}(ii)), each $W_m$ with $m\ge0$ is monic of degree $m$, $r_0=1+x$, and $r_m$ is monic of degree $m$ for $m\ge1$. We consider two sequences of Hankel determinants. The first one belongs to the weights $s=(c+x-1,c,c,\dots)$, $\lambda=(1+x,1,1,\dots)$, and the second to the sums of two adjacent columns of the triangle $(a_{n,j}(x,c))$,
\[
\cE_n:=\det\bigl(a^{\flat}_{i+j,k}\bigr)_{0\le i,j\le n-1},\qquad\cO_n:=\det\bigl(a_{i+j,k}(x,c)+a_{i+j,k-1}(x,c)\bigr)_{0\le i,j\le n-1},
\]
where $a^{\flat}_{n,j}:=a^{(s,\lambda)}_{n,j}\in\Z[x,c]$ for $s=(c+x-1,c,c,\dots)$, $\lambda=(1+x,1,1,\dots)$. Thus $a^\flat_{n,j}$ is the weighted number of Motzkin paths from $(0,0)$ to $(n,j)$ in which level steps have weight $c+x-1$ at height $0$ and $c$ at positive height, down steps ending at height $0$ have weight $1+x$, and all other steps have weight $1$. When the dependence on $k$ and $x$ matters we write $\cE^{(k)}_n(x)$ and $\cO^{(k)}_n(x)$.

\begin{theorem}\label{thm:typeB}
Let $k\ge1$. The determinants $\cE_n$ and $\cO_n$ do not depend on $c$, and for all $n\ge2k+1$
\begin{align}
(-1)^k\bigl(\cE_n\,\cE_{n-2k-1}-x^2\,\cE_{n-1}\,\cE_{n-2k}\bigr)&=(1-x)\,\cO_{n-k-1}^{\,2},\label{eq:E1}\\
(-1)^k\bigl(\cO_n\,\cO_{n-2k-1}-x^2\,\cO_{n-1}\,\cO_{n-2k}\bigr)&=(1-x)(1+x)^2\,\cE_{n-k}^{\,2}\label{eq:E2}
\end{align}
in $\Z[x]$.
\end{theorem}

The two recurrences have the shape of \eqref{eq:main}, with span $2k+1$ instead of $2k+2$, and in each of them the right-hand side is the square of a term of the other sequence. The proof follows that of Theorem~\ref{thm:B} in Section~\ref{sec:HM}, with the quotient ring $K[X]/(r_k)$ in place of $R$. The role of $\Lambda$ is played by two multiplication maps, which exchange the sequence of the classes of the polynomials $U_{j-1}$ and the sequence of the classes of the polynomials $W_j$.

\subsection{The quotient ring \texorpdfstring{$K[X]/(r_k)$}{K[X]/(r\_k)}}\label{ss:modelB}
Let $\Rt:=K[X]/(r_k)$, a $K$-vector space of dimension $k$, and for $j\in\Z$ let
\[
F_j\in\Rt\ \text{ be the class of }U_{j-1},\qquad\Ft_j:=F_{j+1}+F_j\ \text{ the class of }W_j.
\]
Both sequences satisfy the recurrence $f_{j+1}=Xf_j-f_{j-1}$, and $F_0=0$, $F_{-j}=-F_j$, $\Ft_{-j-1}=-\Ft_j$. The families $(F_1,\dots,F_k)$ and $(\Ft_0,\dots,\Ft_{k-1})$ are monic bases of $\Rt$. With $\detx{r_k}$ as in \S\ref{ss:det}, put
\[
e(m):=\Cas_m(F)=\detx{r_k}(F_m,\dots,F_{m+k-1}),\qquad o(m):=\Cas_m(\Ft)=\detx{r_k}(\Ft_m,\dots,\Ft_{m+k-1})
\]
for $m\in\Z$.
As in the proof of Theorem~\ref{thm:C}, $e(m),o(m)\in\Z[x]$, because $r_k$ is monic with coefficients in $\Z[x]$.

\begin{lemma}\label{lem:modelB}
Let $\Lambda'$ and $\Lambda''$ be the multiplications by the classes of $-r_{k-1}$ and of $-P_{k-1}$ on $\Rt$. Then for all $j\in\Z$:
\begin{enumerate}
\item[(i)] $\Lambda'F_j=\Ft_{j+k}-x\,\Ft_{j+k-1}=\Ft_{j-k-1}-x\,\Ft_{j-k}$;
\item[(ii)] $\Lambda''\Ft_j=F_{j+k+1}-x\,F_{j+k}=F_{j-k}-x\,F_{j-k+1}$;
\item[(iii)] $\det\Lambda'=\eps_k\,(1+x)$ and $\det\Lambda''=\eps_k$.
\end{enumerate}
\end{lemma}

\begin{proof}
All sequences occurring in (i) and (ii) satisfy the recurrence $f_{j+1}=Xf_j-f_{j-1}$, so it suffices to verify each identity for two consecutive values of $j$. Note that $\Ft_m-x\Ft_{m-1}$ is the class of $r_m$ and $F_{m+1}-xF_m$ is the class of $P_m$, for all $m\in\Z$, and that $r_{m+1}=Xr_m-r_{m-1}$.

(i) For $j=0$: $\Lambda'F_0=0$; $\Ft_k-x\Ft_{k-1}$ is the class of $r_k$, which is $0$; and $\Ft_{-k-1}-x\Ft_{-k}=-(\Ft_k-x\Ft_{k-1})=0$. For $j=1$: $\Lambda'F_1$ is the class of $-r_{k-1}$; $\Ft_{k+1}-x\Ft_k$ is the class of $r_{k+1}=Xr_k-r_{k-1}$; and $\Ft_{-k}-x\Ft_{-k+1}=-(\Ft_{k-1}-x\Ft_{k-2})$ is the class of $-r_{k-1}$.

(ii) We use $j=0$ and $j=-1$, where $\Ft_0=1$ and $\Ft_{-1}=-1$. For $j=0$: $\Lambda''\Ft_0$ is the class of $-P_{k-1}$; $F_{k+1}-xF_k$ is the class of $P_k=r_k-P_{k-1}$; and $F_{-k}-xF_{-k+1}=-(F_k-xF_{k-1})$ is the class of $-P_{k-1}$. For $j=-1$: $\Lambda''\Ft_{-1}$ is the class of $P_{k-1}$; $F_k-xF_{k-1}$ is the class of $P_{k-1}$; and $F_{-k-1}-xF_{-k}=-(F_{k+1}-xF_k)$ is the class of $-P_k$, which equals the class of $P_{k-1}$.

(iii) By (i), $\Lambda'F_{m-k}$ is the class of $r_m$ for all $m\in\Z$. Since $r_0=1+x$ and $r_1,\dots,r_{k-1}$ are monic of degrees $1,\dots,k-1$, we get, as in the proof of Lemma~\ref{lem:Lambda}(iii),
\[
1+x=\detx{r_k}(r_0,\dots,r_{k-1})=\det(\Lambda')\,\detx{r_k}(F_{-k},\dots,F_{-1})=\det(\Lambda')\,\eps_k.
\]
By (ii), $\Lambda''\Ft_{m-k}$ is the class of $P_m$ for all $m\in\Z$, and $\Ft_{-k},\dots,\Ft_{-1}$ are $-\Ft_{k-1},\dots,-\Ft_0$; hence
\begin{align*}
1=\detx{r_k}(P_0,\dots,P_{k-1})&=\det(\Lambda'')\,\detx{r_k}(\Ft_{-k},\dots,\Ft_{-1})\\
&=\det(\Lambda'')\,(-1)^k(-1)^{k(k-1)/2}=\det(\Lambda'')\,\eps_k.\qedhere
\end{align*}
\end{proof}

\begin{proposition}\label{prop:bilinearB}
For every $\ell\in\Z$,
\begin{align}
o(\ell+k)\,o(\ell-k-1)-x^2\,o(\ell+k-1)\,o(\ell-k)&=(1-x^2)(1+x)\,e(\ell)^2,\label{eq:oo}\\
e(\ell+k+1)\,e(\ell-k)-x^2\,e(\ell+k)\,e(\ell-k+1)&=(1-x)\,o(\ell)^2.\label{eq:ee}
\end{align}
\end{proposition}

\begin{proof}
By Lemma~\ref{lem:modelB}(ii), applied to $j-k$ and to $j+k$, and by Lemma~\ref{lem:modelB}(i), applied to $j+1-k$ and to $j+k+1$, we have for all $j\in\Z$
\begin{align*}
(\partial_xF)_j&=\Lambda''\Ft_{j-k}, & (\partial_{x^{-1}}F)_j&=-x^{-1}\bigl(F_j-xF_{j+1}\bigr)=-x^{-1}\Lambda''\Ft_{j+k},\\
(\partial_x\Ft)_j&=\Lambda'F_{j+1-k}, & (\partial_{x^{-1}}\Ft)_j&=-x^{-1}\bigl(\Ft_j-x\Ft_{j+1}\bigr)=-x^{-1}\Lambda'F_{j+k+1}.
\end{align*}
Moreover $(\partial_x\partial_{x^{-1}}f)_j=\bigl(X-(x+x^{-1})\bigr)f_{j+1}$ for $f=F$ and for $f=\Ft$, and the determinant of the multiplication by $X-(x+x^{-1})$ on $\Rt$ is $(-1)^kr_k(x+x^{-1})=(-x)^{-k}(1+x)$, by Lemma~\ref{lem:cheb}(ii). Let $m=\ell-1$. For the sequence $F$ we obtain
\begin{align*}
\Cas_m(\partial_0F)&=e(\ell), &\Cas_m(\partial_x\partial_{x^{-1}}F)&=(-x)^{-k}(1+x)\,e(\ell),\\
\Cas_m(\partial_xF)&=\det(\Lambda'')\,o(\ell-k-1), &\Cas_m(\partial_0\partial_{x^{-1}}F)&=(-x)^{-k}\det(\Lambda'')\,o(\ell+k),\\
\Cas_m(\partial_{x^{-1}}F)&=(-x)^{-k}\det(\Lambda'')\,o(\ell+k-1), &\Cas_m(\partial_0\partial_xF)&=\det(\Lambda'')\,o(\ell-k),
\end{align*}
and Lemma~\ref{lem:HM} with $(a_1,a_2,a_3)=(0,x,x^{-1})$, multiplied by $x(-x)^k$, gives
\[
(x^2-1)(1+x)\,e(\ell)^2+(\det\Lambda'')^2\bigl(o(\ell-k-1)\,o(\ell+k)-x^2\,o(\ell+k-1)\,o(\ell-k)\bigr)=0,
\]
which is \eqref{eq:oo}, because $(\det\Lambda'')^2=1$. For the sequence $\Ft$ we obtain
\begin{align*}
\Cas_m(\partial_0\Ft)&=o(\ell), &\Cas_m(\partial_x\partial_{x^{-1}}\Ft)&=(-x)^{-k}(1+x)\,o(\ell),\\
\Cas_m(\partial_x\Ft)&=\det(\Lambda')\,e(\ell-k), &\Cas_m(\partial_0\partial_{x^{-1}}\Ft)&=(-x)^{-k}\det(\Lambda')\,e(\ell+k+1),\\
\Cas_m(\partial_{x^{-1}}\Ft)&=(-x)^{-k}\det(\Lambda')\,e(\ell+k), &\Cas_m(\partial_0\partial_x\Ft)&=\det(\Lambda')\,e(\ell-k+1),
\end{align*}
and Lemma~\ref{lem:HM} gives in the same way
\[
(x^2-1)(1+x)\,o(\ell)^2+(\det\Lambda')^2\bigl(e(\ell-k)\,e(\ell+k+1)-x^2\,e(\ell+k)\,e(\ell-k+1)\bigr)=0,
\]
which is \eqref{eq:ee}, because $(\det\Lambda')^2=(1+x)^2$.
\end{proof}

In terms of \S\ref{ss:reduction}, $\partial_x$ does not act on the sequence $F$ as a translation followed by a linear map in general. For instance, for $k=1$ we have $e(m)=U_{m-1}(x-1)$ and $o(m)=W_m(x-1)$, and $\Cas_m(\partial_xF)=\det(\Lambda'')\,o(m-k)$ is not proportional to $e(m+p)$ for any fixed $p$. Instead, $\partial_x$ maps the pair $(F,\Ft)$ to translates of the pair $(\Ft,F)$, so that $\partial_x^2$ is a translation by $1-2k$ followed by $\Lambda'\Lambda''$. This is why the two families of Casoratians are coupled.

\begin{proof}[Proof of Theorem~\ref{thm:typeB}]
We show that for all $n\ge0$
\begin{equation}\label{eq:EOeo}
\cE_n=(-1)^{nk}\eps_k\,e(n-k),\qquad\cO_n=(-1)^{nk}\eps_k\,o(n-k),
\end{equation}
which includes the independence of $c$.

\emph{The determinants $\cO_n$.} By \eqref{eq:shiftm}, $a_{i,k}(x,c)+a_{i,k-1}(x,c)=L_x\bigl((X+c)^i\,r_k\bigr)$ for $i\ge0$. Hence, by Lemma~\ref{lem:jacobi}(iv), by \eqref{eq:christoffel} for the weights $(x,0,0,\dots)$, $(1,1,\dots)$ and $f=r_k$, and by Lemma~\ref{lem:modelB}(ii),(iii),
\begin{align*}
\cO_n=\det\bigl(L_x(X^{i+j}r_k)\bigr)_{i,j<n}&=(-1)^{nk}\detx{r_k}(P_n,\dots,P_{n+k-1})\\
&=(-1)^{nk}\detx{r_k}\bigl(\Lambda''\Ft_{n-k+i}\bigr)_{0\le i\le k-1}=(-1)^{nk}\eps_k\,o(n-k).
\end{align*}

\emph{The determinants $\cE_n$.} Let $L^\flat$ and $p^\flat_m$ be the functional and the polynomials of \S\ref{ss:matrices} for the weights $s=(x-1,0,0,\dots)$, $\lambda=(1+x,1,1,\dots)$, so that $h_0=1$ and $h_m=1+x$ for $m\ge1$. Then $p^\flat_0=1$ and $p^\flat_m=r_m$ for $m\ge1$: indeed $p^\flat_1=X-x+1=r_1$, $p^\flat_2=Xp^\flat_1-(1+x)=X^2+(1-x)X-1-x=r_2$, and for $m\ge2$ both sequences satisfy the same recurrence. By Lemma~\ref{lem:jacobi}(iii),(iv), $(1+x)\,a^\flat_{i,k}=L^\flat\bigl((X+c)^i\,r_k\bigr)$ for $i\ge0$. Hence, by Lemma~\ref{lem:jacobi}(iv), \eqref{eq:christoffel} and Lemma~\ref{lem:modelB}(i),(iii), for $n\ge1$,
\begin{align*}
(1+x)^n\,\cE_n=\det\bigl(L^\flat(X^{i+j}r_k)\bigr)_{i,j<n}&=(1+x)^{n-1}(-1)^{nk}\detx{r_k}(r_n,\dots,r_{n+k-1})\\
&=(1+x)^{n}(-1)^{nk}\eps_k\,e(n-k),
\end{align*}
because the class of $r_m$ is $\Lambda'F_{m-k}$. For $n=0$ we have $\cE_0=1$, and $e(-k)$ is the determinant of $F_{-k},\dots,F_{-1}$, which is $\eps_k$.

\emph{The recurrences.} Let $n\ge2k+1$ and $\ell:=n-2k-1\ge0$. By \eqref{eq:EOeo}, $e(\ell+k+1)e(\ell-k)=(-1)^{(2\ell+2k+1)k}\cE_n\cE_{n-2k-1}=(-1)^k\cE_n\cE_{n-2k-1}$, $e(\ell+k)e(\ell-k+1)=(-1)^k\cE_{n-1}\cE_{n-2k}$ and $o(\ell)^2=\cO_{n-k-1}^{\,2}$, so \eqref{eq:ee} gives \eqref{eq:E1}. Similarly, with $\ell:=n-2k\ge1$, $o(\ell+k)o(\ell-k-1)=(-1)^k\cO_n\cO_{n-2k-1}$, $o(\ell+k-1)o(\ell-k)=(-1)^k\cO_{n-1}\cO_{n-2k}$ and $e(\ell)^2=\cE_{n-k}^{\,2}$, so \eqref{eq:oo} gives \eqref{eq:E2}, because $(1-x^2)(1+x)=(1-x)(1+x)^2$.
\end{proof}

\begin{remark}\label{rem:typeB}
(1) For $k=1$ we have $r_1=X+1-x$, hence $e(m)=U_{m-1}(x-1)$, $o(m)=W_m(x-1)$ and
\[
\cE_n=(-1)^{n+1}U_{n-2}(x-1),\qquad\cO_n=(-1)^{n+1}\bigl(U_{n-1}(x-1)+U_{n-2}(x-1)\bigr).
\]
For $k=2$ the first terms are $(\cE_n)_{n\ge0}=(1,0,0,-1,x,1-x,x^2-x^3,\dots)$ and $(\cO_n)_{n\ge0}=(1,0,-1,0,x^2,1+x-x^2-x^3,\dots)$.

(2) The same two sequences appear as Hankel determinants of linear combinations of columns of the triangles $(a_{n,j}(0,c))$ and $(a_{n,j}(-1,c))$: for every $c$ and all $n\ge0$,
\begin{align*}
\det\bigl(a_{i+j,k}(0,c)+(1-x)\,a_{i+j,k-1}(0,c)-x\,a_{i+j,k-2}(0,c)\bigr)_{0\le i,j\le n-1}&=\eps_k\,\cE_{n+k+1},\\
\det\bigl(a_{i+j,k}(-1,c)-x\,a_{i+j,k-1}(-1,c)\bigr)_{0\le i,j\le n-1}&=(-1)^k\eps_k\,\cO_{n+k}
\end{align*}
$(a_{r,-1}:=0)$. Indeed, $r_k=U_k+(1-x)U_{k-1}-xU_{k-2}=P_{k,-1}-xP_{k-1,-1}$, so by \eqref{eq:shiftm}, Lemma~\ref{lem:jacobi}(iv) and \eqref{eq:christoffel} the two determinants are
\begin{align*}
\det\bigl(L_0(X^{i+j}r_k)\bigr)_{i,j<n}&=(-1)^{nk}\detx{r_k}(U_n,\dots,U_{n+k-1})=(-1)^{nk}e(n+1),\\
\det\bigl(L_{-1}(X^{i+j}r_k)\bigr)_{i,j<n}&=(-1)^{nk}\detx{r_k}(W_n,\dots,W_{n+k-1})=(-1)^{nk}o(n),
\end{align*}
and \eqref{eq:EOeo} applies. Moreover, since $a_{n,j}(-x,-c)=(-1)^{n-j}a_{n,j}(x,c)$, the Hankel determinants of the differences are $\det\bigl(a_{i+j,k}(x,c)-a_{i+j,k-1}(x,c)\bigr)_{i,j<n}=(-1)^{nk}\cO^{(k)}_n(-x)$.

(3) For $x=1$ the numbers $a^\flat_{n,j}$ count all paths from $(0,0)$ to $(n,j)$ with steps $(1,1)$, $(1,-1)$ and $(1,0)$, the latter with weight $c$, without the restriction to non-negative heights (reflect the excursions below the axis). The determinants $\cE_n$ are then evaluated in \cite[Theorem~3]{CK11}, and \eqref{eq:E1}, \eqref{eq:E2} degenerate to $\cE_n\cE_{n-2k-1}=\cE_{n-1}\cE_{n-2k}$ and $\cO_n\cO_{n-2k-1}=\cO_{n-1}\cO_{n-2k}$. For $x=0$, $\cE_n=D_{k,n}(-1)$.

(4) Let $\rho_1,\dots,\rho_k$ be the roots of $r_k$ in an algebraic closure of $K$ and $\omega_j+\omega_j^{-1}=\rho_j$. Then $\omega_1^{\pm1},\dots,\omega_k^{\pm1}$ are the roots of $(t^{2k+1}-xt^{2k}+xt-1)/(t-1)=t^k\,r_k(t+t^{-1})$, and they are pairwise distinct, by the argument used for $P_k$ in \cite[Remark~4.4]{Yu} (for $x=0$ the polynomial $r_k=W_k$ has $k$ distinct roots, and $r_k(2)=2k+1-(2k-1)x$ and $r_k(-2)=(-1)^k(1+x)$ are non-zero). As in Remark~\ref{rem:symplectic}, evaluation at the roots gives $e(m)=\symp_{((m-1)^k)}(\omega)$ for $m\ge1$. Since $(\omega^{1/2}-\omega^{-1/2})\,W_n(\omega+\omega^{-1})=\omega^{n+1/2}-\omega^{-n-1/2}$ by \eqref{eq:cheb}, it also gives, for $m\ge0$,
\[
o(m)=\frac{\det\bigl(\omega_j^{m+i+1/2}-\omega_j^{-m-i-1/2}\bigr)_{0\le i\le k-1,\,1\le j\le k}}{\det\bigl(\omega_j^{i+1/2}-\omega_j^{-i-1/2}\bigr)_{0\le i\le k-1,\,1\le j\le k}}=\sorth_{(m^k)}(\omega_1,\dots,\omega_k),
\]
the character of $\SO(2k+1)$ indexed by the rectangular partition $(m^k)$ \cite[\S24.2, (24.28)]{FH}, and the quotient does not depend on the choice of the square roots $\omega_j^{1/2}$. Moreover $x=1+\rho_1+\dots+\rho_k=\sorth_{(1)}(\omega)$. Thus Proposition~\ref{prop:bilinearB} is a coupled system of bilinear relations between rectangular symplectic and odd orthogonal characters at the point $\omega$.
\end{remark}

\subsection{Down steps ending at height \texorpdfstring{$0$ with weight $u$}{0 with weight u}}\label{ss:u}
Let $u$ be an indeterminate and let $a^{[u]}_{n,j}(c):=a^{(s,\lambda)}_{n,j}\in\Z[u,c]$ for $s=(c,c,c,\dots)$ and $\lambda=(u,1,1,\dots)$. This is the generating polynomial of the Motzkin paths from $(0,0)$ to $(n,j)$ that never go below height $0$, where a path with $l$ level steps and with $\beta_0$ down steps ending at height $0$ has the weight $c^{\,l}u^{\beta_0}$. For $u=1$ we have $a^{[1]}_{n,j}(c)=a_{n,j}(0,c)$, and for $u=2$ the numbers $a^{[2]}_{n,j}(c)$ count the paths without the restriction to non-negative heights, as in Remark~\ref{rem:typeB}(3). Let
\[
D^{[u]}_{k,n}:=\det\bigl(a^{[u]}_{i+j,k}(c)\bigr)_{0\le i,j\le n-1}\in\Z[u,c].
\]

\begin{theorem}\label{thm:return}
Let $k\ge1$. The determinants $D^{[u]}_{k,n}$ do not depend on $c$, and:
\begin{enumerate}
\item[(i)] If $k=2k'+1$ is odd, then for all $n\ge0$
\[
D^{[u]}_{k,2n+1}=0,\qquad D^{[u]}_{k,2n}=(-1)^n\bigl(D_{k',n}(u-1)\bigr)^2,
\]
where $D_{k',n}(x)$ is the determinant of Theorem~\ref{thm:B} for $k'$ in place of $k$ $($so $D_{k',n}=1$ if $k'=0)$.
\item[(ii)] If $k=2k'$ is even, then for all $n\ge0$
\[
D^{[u]}_{k,n}=\cE^{(k')}_{\lceil n/2\rceil}(u-1)\cdot\cO^{(k')}_{\lfloor n/2\rfloor}(u-1),
\]
where $\cE^{(k')}_n(x)$ and $\cO^{(k')}_n(x)$ are the determinants of Theorem~\ref{thm:typeB} for $k'$ in place of $k$.
\end{enumerate}
\end{theorem}

Thus, for odd $k$, the sequence $\bigl((-1)^nD^{[u]}_{k,2n}\bigr)_{n\ge0}$ consists of the squares of a sequence satisfying the quadratic recurrence of Theorem~\ref{thm:B} (with $(k-1)/2$ and $u-1$ in place of $k$ and $x$), and for even $k$ the sequence $(D^{[u]}_{k,n})_{n\ge0}$ is obtained by interlacing the two sequences of Theorem~\ref{thm:typeB}, which satisfy \eqref{eq:E1} and \eqref{eq:E2}. In the proof below, Step~3 is the classical factorization of a determinant whose entries vanish when the sum of the row and column indices is odd \cite[Lemma~5]{CK11}, and Step~4 is its analogue for entries that vanish when this sum is even. What is new is the identification of the factors, in Step~2.

\begin{proof}
\emph{Step 1: independence of $c$.} Let $L$ and $p_m$ be the functional and the polynomials of the weights $(0,0,\dots)$, $(u,1,1,\dots)$ over $\Z[u]$, so that $h_k=u$. By Lemma~\ref{lem:jacobi}(iii),(iv), $u\,a^{[u]}_{i,k}(c)=L\bigl((X+c)^ip_k\bigr)$ for all $i\ge0$, and therefore $u^n\,D^{[u]}_{k,n}=\det\bigl(L(X^{i+j}p_k)\bigr)_{i,j<n}$, which does not contain $c$. So we may take $c=0$, and we write $b_i:=a^{[u]}_{i,k}(0)$.

\emph{Step 2: the squares of the transfer matrix.} Let $J$ be the matrix of the weights $s=(0,0,\dots)$, $\lambda=(u,1,1,\dots)$. Its non-zero entries are $J_{m,m+1}=1$ and $J_{m+1,m}=\lambda_{m+1}$, so $(J^2)_{i,j}=0$ if $i-j$ is odd, and for $i\ge0$
\[
(J^2)_{i,i+2}=1,\qquad(J^2)_{i,i}=\lambda_{i+1}+\lambda_i,\qquad(J^2)_{i+2,i}=\lambda_{i+2}\lambda_{i+1}\qquad(\lambda_0:=0).
\]
Consequently $J^2$ maps the span of the $e_m$ with $m$ odd to itself, and so for the $e_m$ with $m$ even. In the basis $e_1,e_3,e_5,\dots$ its matrix is tridiagonal with the diagonal entries $\lambda_2+\lambda_1=1+u$ and $\lambda_{2i+2}+\lambda_{2i+1}=2$ $(i\ge1)$, and with all entries next to the diagonal equal to $1$: this is the matrix of the weights $(c+x,c,c,\dots)$, $(1,1,\dots)$ for $c=2$ and $x=u-1$. In the basis $e_0,e_2,e_4,\dots$ its matrix is tridiagonal with the diagonal entries $\lambda_1=u$ and $\lambda_{2i+1}+\lambda_{2i}=2$ $(i\ge1)$, the entries $1$ above the diagonal, and the entries $\lambda_2\lambda_1=u$ and $\lambda_{2i+2}\lambda_{2i+1}=1$ $(i\ge1)$ below it: this is the matrix of the weights $(c+x-1,c,c,\dots)$, $(1+x,1,1,\dots)$ for $c=2$ and $x=u-1$. Since $(J^{2n+1})_{0,2j+1}=(J^{2n})_{1,2j+1}$ and $(J^{2n+2})_{0,2j}=(J^{2n+1})_{0,2j-1}+\lambda_{2j+1}(J^{2n+1})_{0,2j+1}$, we obtain, for all $n,j\ge0$,
\begin{equation}\label{eq:aerated}
\begin{aligned}
a^{[u]}_{2n+1,2j+1}(0)&=a_{n,j}(u-1,2),\qquad a^{[u]}_{2n,2j}(0)=a^\flat_{n,j}\big|_{x=u-1,\,c=2},\\
a^{[u]}_{2n+2,2j}(0)&=a_{n,j}(u-1,2)+a_{n,j-1}(u-1,2)\quad(j\ge1),
\end{aligned}
\end{equation}
and $a^{[u]}_{n,j}(0)=0$ if $n-j$ is odd.

\emph{Step 3: $k=2k'$ even.} Then $b_i=0$ for odd $i$, so the entry $b_{i+j}$ of the Hankel matrix vanishes unless $i\equiv j\pmod2$. Permuting rows and columns simultaneously so that the even indices come first produces a block diagonal matrix with the blocks $(b_{2i+2j})_{0\le i,j<\lceil n/2\rceil}$ and $(b_{2i+2j+2})_{0\le i,j<\lfloor n/2\rfloor}$. By \eqref{eq:aerated} the first block is $\bigl(a^\flat_{i+j,k'}\bigr)$ and the second one is $\bigl(a_{i+j,k'}(u-1,2)+a_{i+j,k'-1}(u-1,2)\bigr)$; their determinants are $\cE^{(k')}_{\lceil n/2\rceil}(u-1)$ and $\cO^{(k')}_{\lfloor n/2\rfloor}(u-1)$, which do not depend on $c$ by Theorem~\ref{thm:typeB}. This is (ii).

\emph{Step 4: $k=2k'+1$ odd.} Then $b_i=0$ for even $i$, so the entry $b_{i+j}$ vanishes unless $i\not\equiv j\pmod2$. If $n$ is odd, the $\lceil n/2\rceil$ rows with even index have non-zero entries only in the $\lfloor n/2\rfloor$ columns with odd index, so they are linearly dependent, and $D^{[u]}_{k,n}=0$. If $n=2N$, the same permutation as in Step~3 produces the matrix $\bigl(\begin{smallmatrix}0&C\\ C^{\T}&0\end{smallmatrix}\bigr)$ with $C=(b_{2i+2j+1})_{0\le i,j<N}=\bigl(a_{i+j,k'}(u-1,2)\bigr)_{0\le i,j<N}$ by \eqref{eq:aerated}, whose determinant is $(-1)^{N^2}\det(C)^2=(-1)^N\bigl(D_{k',N}(u-1)\bigr)^2$. This is (i).
\end{proof}

\begin{remark}[Paths without level steps]\label{cor:dyck}
For $c=0$ the paths have no level steps. Let $a^{\mathrm D}_{n,j}:=a^{[u]}_{n,j}(0)\in\Z[u]$ be the generating polynomial of the paths with steps $(1,1)$ and $(1,-1)$ from $(0,0)$ to $(n,j)$ that never go below height $0$, where a down step ending at height $0$ has weight $u$ and all other steps have weight $1$. By \eqref{eq:aerated}, $\det\bigl(a^{\mathrm D}_{2i+2j+1,2k'+1}\bigr)_{0\le i,j\le n-1}=D_{k',n}(u-1)$ for $k'\ge0$. For $k'\ge1$ the determinants $\mathcal A_n:=\det\bigl(a^{\mathrm D}_{2i+2j,\,2k'}\bigr)_{0\le i,j\le n-1}$ and $\mathcal B_n:=\det\bigl(a^{\mathrm D}_{2i+2j+2,\,2k'}\bigr)_{0\le i,j\le n-1}$ are $\cE^{(k')}_n(u-1)$ and $\cO^{(k')}_n(u-1)$, again by \eqref{eq:aerated}. So \eqref{eq:E2} and \eqref{eq:E1}, for $k'$, $u-1$ and $n+k'+1$ in place of $k$, $x$ and $n$, give for all $n\ge k'$, in $\Z[u]$,
\begin{align}
(u-1)^2\,\mathcal B_{n+k'}\,\mathcal B_{n+1-k'}-\mathcal B_{n-k'}\,\mathcal B_{n+1+k'}&=(-1)^{k'}\,(u-2)\,u^2\,\mathcal A_{n+1}^2,\label{eq:dyck1}\\
(u-1)^2\,\mathcal A_{n+k'}\,\mathcal A_{n+1-k'}-\mathcal A_{n-k'}\,\mathcal A_{n+1+k'}&=(-1)^{k'}\,(u-2)\,\mathcal B_n^2.\label{eq:dyck2}
\end{align}
\end{remark}

\begin{remark}\label{rem:u}
For $u=1$ (paths restricted to non-negative heights) and $u=2$ (unrestricted paths) the determinants $D^{[u]}_{k,n}$ are evaluated in \cite[Theorems~1 and~3]{CK11}. Theorem~\ref{thm:return} interpolates between these two cases. A computer check for $k\le7$ shows that for $u=1,2$ it agrees with these evaluations. The Hankel determinants $D^{[u]}_{k,n}$ of the full sequence are products of two factors or squares, and the factors satisfy recurrences of the shape of \eqref{eq:main}. For the sequences $(D^{[u]}_{k,n})_n$ themselves we found no non-trivial three-term recurrence (see \S\ref{ss:numerical}). For even $k$ and $n>k$ the two factors in Theorem~\ref{thm:return}(ii) are, up to sign, a symplectic and an odd orthogonal character of rectangular shape at the same point (Remark~\ref{rem:typeB}(4)). Factorizations of rectangular Schur functions into such products are given in \cite[Corollary~8]{AB} (see also \cite{CiK} for analogous factorizations of $s_{(M^n)}(x_1^{\pm1},\dots,x_n^{\pm1})$).
\end{remark}
\section{Further remarks}\label{sec:remarks}

\subsection{The recurrence of Han and Pedon}\label{ss:HP}
\begin{remark}\label{rem:HP}
The argument of Corollary~\ref{cor:DJB} works for arbitrary sequences. Let $(f_i)_{i\ge0}$ be a sequence in a commutative ring and $\Delta^{(\ell)}_j:=\det(f_{a+b+\ell})_{0\le a,b\le j-1}$ $(\ell,j\ge0)$. The Desnanot--Jacobi identity for the Hankel matrix $(f_{a+b+\ell-1})_{0\le a,b\le j}$ is the discrete-time Toda equation \cite[(2.8)]{SNvdK}
\begin{equation}\label{eq:DJgen}
\bigl(\Delta^{(\ell)}_j\bigr)^2=\Delta^{(\ell+1)}_j\,\Delta^{(\ell-1)}_j-\Delta^{(\ell-1)}_{j+1}\,\Delta^{(\ell+1)}_{j-1}\qquad(\ell\ge1,\ j\ge1).
\end{equation}
If $\Delta^{(\ell+1)}_j=\vartheta_j\Delta^{(\ell)}_{j+p}$ and $\Delta^{(\ell)}_j=\vartheta'_j\Delta^{(\ell-1)}_{j+p}$ for all $j\ge0$, with signs $\vartheta_j,\vartheta'_j$ and an integer $p\ge1$, then \eqref{eq:DJgen} gives $(\Delta^{(\ell)}_j)^2=\vartheta_j\vartheta'_{j-p}\Delta^{(\ell)}_{j+p}\Delta^{(\ell)}_{j-p}-\vartheta'_{j+1-p}\vartheta_{j-1}\Delta^{(\ell)}_{j+1-p}\Delta^{(\ell)}_{j-1+p}$ for $j\ge p$, a recurrence of Gale--Robinson type with span $2p$. For the $q$-deformed $m$-th metallic number and $0\le\ell\le m$, Han and Pedon prove $\Delta^{(\ell+1)}_j=(-1)^{j+m(m+2\ell-1)/2}\Delta^{(\ell)}_{j+m+1}$ \cite[Theorem~D]{HP} (they write $n$ for $m$). For $1\le\ell\le m$ this applies with $p=m+1$, both products of signs equal $-1$ (their exponents are $2j+m^2+2m\ell-3m-1$, which is odd), and we obtain
\[
\Delta^{(\ell)}_{j+m+1}\,\Delta^{(\ell)}_{j-m-1}=\Delta^{(\ell)}_{j+m}\,\Delta^{(\ell)}_{j-m}-\bigl(\Delta^{(\ell)}_j\bigr)^2\qquad(j\ge m+1).
\]
This is the recurrence of \cite[Theorem~C]{HP}, conjectured in \cite{OP}. In \cite{HP} it is proved for $\ell=0$ by an explicit evaluation of the determinants and a case analysis (for $m\ge3$, the cases $m=1,2$ being in \cite{OP}), and then transferred to the other values of $\ell$ by \cite[Theorem~D]{HP}.
\end{remark}

\subsection{Solitons and a nodal curve}\label{ss:soliton}

\begin{remark}\label{rem:soliton}
If the maps $\Lambda_1,\Lambda_2,\Lambda_3$ of Proposition~\ref{thm:reduction} are simultaneously diagonalizable, a coordinate of $F_j$ with respect to a common eigenbasis is a scalar sequence $\varphi$ with $\varphi_{j+1}-a_i\varphi_j=\mu_i\varphi_{j+\pi_i}$ $(i=1,2,3)$, where $\mu_i$ is the eigenvalue of $\Lambda_i$. The sequences $\varphi_j=\omega^j$ with $\omega^{-\pi_i}(\omega-a_i)=\mu_i$ have this property, and so has $\alpha\,\omega^j+\beta\,\omega'^j$ whenever
\begin{equation}\label{eq:pair}
\omega^{-\pi_i}(\omega-a_i)=\omega'^{-\pi_i}(\omega'-a_i)\qquad(i=1,2,3).
\end{equation}
In the case of \S\ref{ss:proofs}, $(a_1,a_2,a_3)=(0,x,x^{-1})$ and $(\pi_1,\pi_2,\pi_3)=(1,-k,k+1)$, the condition \eqref{eq:pair} is empty for $i=1$, and for $\omega'=\omega^{-1}$ the conditions for $i=2$ and $i=3$ both reduce to $\omega^{k}(\omega-x)=\omega^{-k}(\omega^{-1}-x)$, that is, to $Q(\omega)=0$. Accordingly, the Casoratian \eqref{eq:solitondet} is, after division by $\prod_i\omega_{j_i}^{-m}$, a sum of $2^r$ exponentials in $m$, of the form of an $r$-soliton $\tau$-function, and $D_{k,n}$ involves all $k$ pairs $(\omega_j,\omega_j^{-1})$. The condition $Q(\omega_j)=0$ says that $\mu(t):=t^k(t-x)$ takes the same value at $t=\omega_j^{\pm1}$. Together with $X=t+t^{-1}$, the function $\mu$ defines a map from the projective line to the curve
\[
\mu^2-\bigl(P_{k+1}(X)-P_{k-1}(X)\bigr)\mu+(1+x^2-xX)=0,
\]
because $\mu(t)+\mu(t^{-1})=(P_{k+1}-P_{k-1})(t+t^{-1})$ and $\mu(t)\mu(t^{-1})=1+x^2-x(t+t^{-1})$. Since $\mu(t)-\mu(t^{-1})=(t-t^{-1})P_k(t+t^{-1})$ by \eqref{eq:cheb}, the discriminant of this quadratic equation is $(X^2-4)P_k(X)^2$, so the curve is rational, with $k$ nodes at the roots of $P_k$. On the eigenline belonging to the root $\rho_j$ of $P_k$ the map $\Lambda$ acts by $-P_{k-1}(\rho_j)=\mu(\omega_j)$, and $\Lambda^2=1+x^2-xX$ (Lemma~\ref{lem:Lambda}(ii)) is the product of the two values of $\mu$ over a point $X$. The Hankel determinant solutions of Somos-type recurrences obtained from hyperelliptic curves of genus one and two \cite{Hon} are of a different, non-degenerate kind.
\end{remark}

\subsection{Numerical observations}\label{ss:numerical}
We searched by computer, in exact arithmetic, for three-term relations $D_nD_{n-N}=\alpha D_{n-i}D_{n-N+i}+\beta D_{n-j}D_{n-N+j}$ $(n\ge N)$, with constants $\alpha,\beta$, $1\le i<j\le N/2$ and span $N\le4k+8$, among the terms with $n\le8k+20$ of a sequence $(D_n)_{n\ge0}$. Indices at which all three products vanish were discarded, and a relation was counted only if at least eight indices remained and they determined $\alpha$ and $\beta$. This excludes relations that hold trivially, such as $D^{[u]}_{k,n}D^{[u]}_{k,n-N}=0$ for odd $k$ and odd $N$ (Theorem~\ref{thm:return}(i)). We found no non-trivial three-term relation for $(D^{[u]}_{k,n})_n$ $(k=3,4$, $u\in\{4,-3/7\})$, for $(\cE_n)_n$ or $(\cO_n)_n$ $(k=2,3$, $x=3)$, and for the Hankel determinants $D_n$ of the $k$-th column for the weights $s=(c+x,c,c,\dots)$, $\lambda=(u,1,1,\dots)$ with $(x,u)\in\{(2,5),(3,5)\}$ and $k=2,3$. For $k=2$ and $u=x+1$, however, these $D_n$ satisfy $D_nD_{n-12}=\varkappa^2D_{n-3}D_{n-9}+(1-\varkappa^2)D_{n-6}^{\,2}$ $(n\ge12)$ with $\varkappa=\rho^3-3\rho$ and $\rho:=1+x$. The reason is that $p_2=X^2-xX-u$ then factors as $(X+1)(X-\rho)$. By \eqref{eq:christoffel} and Lemma~\ref{lem:jacobi}, $D_n=u^{-1}\detx{p_2}(p_n,p_{n+1})$ for $n\ge1$ (with $c=0$, which suffices by Lemma~\ref{lem:jacobi}(iv)), and evaluating at the two roots gives, for all $j\ge0$,
\[
D_{3j}=-\frac{U_{3j-2}(\rho)+U_{3j-3}(\rho)}{\rho+1},\qquad D_{3j+1}=D_{3j+2}=\frac{U_{3j-1}(\rho)}{\rho+1}
\]
(for $j=0$ these are $D_0=1$ and $D_1=D_2=0$). Since $U_{m+3}+U_{m-3}=(X^3-3X)\,U_m$ for all $m\in\Z$, each of the sequences $y_j=D_{3j+r}$ $(r=0,1,2)$ satisfies $y_{j+1}=\varkappa\,y_j-y_{j-1}$ for $j\ge1$, and every sequence with this property satisfies $y_{j+2}\,y_{j-2}=\varkappa^2\,y_{j+1}\,y_{j-1}+(1-\varkappa^2)\,y_j^{\,2}$.

\appendix
\section{Proof of the identities (32) and (33) of Cigler}\label{app:C32}

\begin{proof}[Proof of Corollary~\ref{cor:back1}(ii)]
We specialize $x$ to $1$ in Corollary~\ref{cor:back1}(i). Since $P_k$ is monic in $X$, division with remainder by $P_k$ commutes with this specialization. Hence $M_j(m)\big|_{x=1}$ is the determinant, with respect to a monic basis of $\Q[X]/(U_k-U_{k-1})$, of the classes $\upsilon_m,\dots,\allowbreak\widehat{\upsilon_{m+j}},\dots,\allowbreak\upsilon_{m+k}$, where $\upsilon_i$ denotes the class of $U_{i-1}$ $(i\in\Z)$. These are the elements $F_i$ of \S\ref{ss:Lambda} for $\kk=\Q$ and $x=1$. We claim that
\begin{equation}\label{eq:upsilon}
\upsilon_0=0,\qquad\upsilon_{-i}=-\upsilon_i,\qquad\upsilon_{i+2k+1}=-\upsilon_i,\qquad\upsilon_{2k+1-i}=\upsilon_i\qquad(i\in\Z).
\end{equation}
The first two relations are \eqref{eq:Frec}, and the fourth one follows from the second and the third. For the third one, the equality of ($\Lambda1$) and ($\Lambda2$) for $x=1$ reads $\upsilon_{j+k+1}+\upsilon_{j-k}=\upsilon_{j+k}+\upsilon_{j-k-1}$. Hence $\upsilon_{j+k+1}+\upsilon_{j-k}$ does not depend on $j$, and for $j=0$ it is $\upsilon_{k+1}-\upsilon_k=0$, the class of $U_k-U_{k-1}$.

Now let $m\in\Z$ and write $m=(2k+1)\ell+r$ with $\ell\in\Z$ and $0\le r\le2k$. By \eqref{eq:upsilon}, $(\upsilon_m,\dots,\upsilon_{m+k})$ is $(-1)^\ell$ times
\[
(0,\upsilon_1,\dots,\upsilon_k),\qquad(\upsilon_1,\dots,\upsilon_k,\upsilon_k),\qquad(\upsilon_k,\upsilon_k,\upsilon_{k-1},\dots,\upsilon_1),\qquad(\upsilon_k,\upsilon_{k-1},\dots,\upsilon_1,0)
\]
for $r=0$, $r=1$, $r=k$ and $r=k+1$, respectively (for $k=1$ the cases $r=1$ and $r=k$ coincide). Since $\det(\upsilon_1,\dots,\upsilon_k)=1$ and $\det(\upsilon_k,\dots,\upsilon_1)=(-1)^{\binom k2}$, the non-zero values of $M_j(m)\big|_{x=1}$ are
\begin{align*}
M_0&=(-1)^{\ell k}\quad(r=0), & M_{k-1}=M_k&=(-1)^{\ell k}\quad(r=1),\\
M_0=M_1&=(-1)^{\ell k+\binom k2}\quad(r=k), & M_k&=(-1)^{\ell k+\binom k2}\quad(r=k+1).
\end{align*}
For all other $r$ (they exist only if $k\ge2$) all $M_j(m)\big|_{x=1}$ vanish, because the vectors $\upsilon_m,\dots,\upsilon_{m+k}$ lie in a subspace of dimension at most $k-1$. For $2\le r\le k-1$ they are, up to signs, among $\upsilon_2,\dots,\upsilon_k$, and for $k+2\le r\le2k$ they are, up to signs, among $0,\upsilon_1,\dots,\upsilon_{k-1}$.

Insert this into (i) with $m=n+1$, that is, $n=(2k+1)\ell+r-1$, where now $\ell\ge0$. Note that $Y$ is a non-zero element of $\Z[c]$ and that the size of the determinant on the left of (i) is $n+k+2=(2k+1)\ell+r+k+1$. For $r=0$ (then $\ell\ge1$) the sum in (i) is $(-1)^{\ell k}\phi_{n+1}$ and $(-1)^{nk+\ell k}=(-1)^k$, as $n\equiv\ell-1\pmod2$, which gives the first formula. (For $\ell=0$ both sides of the first formula vanish, by Proposition~\ref{prop:backward} and $\phi_0=0$.) For $r=1$ the sum is $(-1)^{\ell k}(\phi_{n+k}+\phi_{n+k+1})=(-1)^{\ell k}\Phi_\ell$ and $(-1)^{nk+\ell k}=1$, as $n\equiv\ell$, which gives the second formula. For $r=k$ the sum is $(-1)^{\ell k+\binom k2}(\phi_{n+1}+\phi_{n+2})=(-1)^{\ell k+\binom k2}\Phi_\ell$ and $(-1)^{nk+\ell k}=(-1)^{(k-1)k}=1$, so the sign is $(-1)^{\binom{k+2}2+\binom k2}=-1$, and $n+k+2=(2k+1)(\ell+1)$, which gives the third formula. For $r=k+1$ the sum is $(-1)^{\ell k+\binom k2}\phi_{n+k+1}$ and $(-1)^{nk+\ell k}=(-1)^{k^2}=(-1)^k$, so the sign is $(-1)^{\binom{k+2}2+\binom k2+k}=(-1)^{k+1}$, and $n+k+2=(2k+1)(\ell+1)+1$, which gives the fourth formula. The sizes $n+k+2$ with $r\in\{0,1,k,k+1\}$ are those congruent to $k+1$, $k+2$, $0$, $1$ modulo $2k+1$. For all other sizes $n+k+2\ge k+2$ the sum in (i) vanishes, and for the sizes $1,\dots,k+1$ the determinant vanishes by Proposition~\ref{prop:backward}.
\end{proof}

\begin{proof}[Proof of Corollary~\ref{rem:C32}]
Let $x=1$, and let $Y$ and $\Phi_\ell$ be as in Corollary~\ref{cor:back1}(ii). For $x=1$ the identity \eqref{eq:shift1} reads $Y\,D_{1,k,n}(1,c)=(-1)^{nk}\eps_k\sum_{i=0}^kM_i(n-k)\big|_{x=1}\,(\phi_{n+i+1}+\phi_{n+i})$. We insert the values of $M_i(m)\big|_{x=1}$ for $m=n-k$ which were determined in the proof of Corollary~\ref{cor:back1}(ii), and use $\phi_{a+1}+\phi_{a-1}=c\,\phi_a$. The cases $r=k+1$, $r=0$, $r=1$ and $r=k$ of that proof give, in this order, for all $\ell\ge0$,
\begin{align*}
Y\,D_{1,k,(2k+1)\ell}(1,c)&=\Phi_\ell, & Y\,D_{1,k,(2k+1)\ell+k}(1,c)&=(-1)^{\binom k2}\,\Phi_\ell,
\end{align*}
\begin{align*}
Y\,D_{1,k,(2k+1)\ell+k+1}(1,c)&=(-1)^{\binom{k+1}2}(c+2)\,\phi_{(2k+1)(\ell+1)},\\
Y\,D_{1,k,(2k+1)\ell+2k}(1,c)&=(-1)^k(c+2)\,\phi_{(2k+1)(\ell+1)},
\end{align*}
and $D_{1,k,n}(1,c)=0$ for all other $n\ge0$. $($In the first case $n-k=(2k+1)(\ell-1)+k+1$. The signs are computed as in the proof of Corollary~\ref{cor:back1}(ii), with $\eps_k(-1)^{\binom k2}=(-1)^k$.$)$ This is \cite[Theorem~1.7]{CS}, formerly Conjecture~8 of \cite{C23}, so we obtain another proof of it. In \cite[Theorem~1.7]{CS} the values are written as $(\phi_{(2k+1)\ell}+\phi_{(2k+1)(\ell+1)})/\phi_{2k+1}$ and $(\psi_k+\psi_{k+1})\,\phi_{(2k+1)(\ell+1)}/\phi_{2k+1}$ in place of $\Phi_\ell/Y$ and $(c+2)\,\phi_{(2k+1)(\ell+1)}/Y$, where $\psi_m:=\phi_{m+1}-\phi_{m-1}$ are the Lucas polynomials. The two forms agree, because for all $a\in\Z$
\[
\phi_{a+2k+1}+\phi_a=(\phi_{a+k+1}+\phi_{a+k})(\phi_{k+1}-\phi_k),\qquad\psi_{k+1}+\psi_k=(c+2)(\phi_{k+1}-\phi_k)
\]
(substitute $c=t+t^{-1}$ and use \eqref{eq:cheb}). For $a=0$ the first identity reads $\phi_{2k+1}=Y\,(\phi_{k+1}-\phi_k)$.

We compare the four evaluations with Corollary~\ref{cor:back1}(ii), namely the third and the second formula there for $\ell=n$ with the first two evaluations, and the fourth and the first formula there, for $\ell=n$ and $\ell=n+1$, with the last two. Since $\binom{k+2}2+\binom k2$ is odd and $\binom{k+2}2+\binom{k+1}2\equiv k+1\pmod 2$, this gives the four identities of Corollary~\ref{rem:C32}, after division by $Y\ne0$. The vanishing statement is the last assertion of Corollary~\ref{cor:back1}(ii).
\end{proof}

\section{Computer checks}\label{ss:checks}
None of the proofs relies on computer calculations. As an independent check, all numbered statements and the main intermediate identities in the proofs were verified by computer in exact arithmetic, over the integers, the rationals, prime fields and polynomial rings. Path numbers and Hankel determinants are computed from their definitions, and determinants in quotient rings by division with remainder or by the three-term recurrence of the residue classes. The scripts (in Python, using only the standard library), their complete output and a table are available as ancillary files. The table lists, for each numbered statement of this paper, the script that checks it and the range that was tested. The scripts also contain negative controls, that is, wrong variants of statements that must fail.

For example, Theorem~\ref{thm:B} was checked in $\Z[x,c]$ for $k\le6$ and $n\le2k+10$, and with $x$ an indeterminate and three values of $c$ for $k\le10$ and $n\le2k+18$. The identity \eqref{eq:C2} was checked with $x$ an indeterminate for $k\le12$ and $|\ell|\le50$, Theorem~\ref{thm:hankelgen} for six families of weights with up to three indeterminates and for random weights with $k\le7$ and $N\le4$, and Proposition~\ref{thm:spectral} for random monic polynomials of degree at most $6$ and, for $P=P_k$, in $\Z[x,y]$ for $k\le5$. The identities of Theorem~\ref{thm:curve} were checked with $x$ an indeterminate for $k\le4$ and at rational points for $k\le6$, Theorem~\ref{thm:typeB} with $x$ an indeterminate for $k\le5$, and Theorem~\ref{thm:return} with $u$ an indeterminate for $k\le8$ (both for fixed values of $c$). The identities involving symplectic and odd orthogonal characters in Remarks~\ref{rem:symplectic} and~\ref{rem:typeB}(4) were checked exactly, with the bialternant formula as the definition of the characters, but only at sample points: random rational points for Remark~\ref{rem:symplectic}, and random rational points and roots in prime fields for Remark~\ref{rem:typeB}(4). Statements about other papers are not checked by computer.

\section*{Declaration of generative AI and AI-assisted technologies in the manuscript preparation process}
The results of this paper, their proofs, the text and the verification scripts were produced by AI agents (Anthropic's Claude models) working under the direction of the author. The author has checked the proofs of Sections~\ref{sec:prelim} and~\ref{sec:HM}, of Proposition~\ref{thm:reduction}, Corollary~\ref{thm:module} and Theorem~\ref{thm:hankelgen} line by line. All numbered statements were checked by exact computer algebra in the ranges listed in the table described in Appendix~\ref{ss:checks}; the identities involving characters in Remarks~\ref{rem:symplectic} and~\ref{rem:typeB}(4) were checked at sample points. The author takes full responsibility for the content of the paper.

\section*{Acknowledgements}
The author thanks Tianxing Yu for sharing his proof of Cigler's conjecture and the manuscript \cite{Yu}, on which this paper builds.

\end{document}